\documentclass[a4paper, reqno, english,11pt]{smfart}

\usepackage{a4wide}

\RequirePackage{mathrsfs} \let\mathcal\mathscr

\usepackage{amsfonts, amsthm, amsmath, amssymb}

\numberwithin{equation}{section}

\renewcommand\AA{\mathbb{A}}
\newcommand\QQ{\mathbb{Q}}
\newcommand\CC{\mathbb{C}}
\newcommand\NN{\mathbb{N}}
\newcommand\RR{\mathbb{R}}
\newcommand\ZZ{\mathbb{Z}}

\newcommand\PP{\mathbb{P}}

\newcommand\cB{\mathcal{B}}

\newcommand\cU{\mathcal{U}}
\newcommand\cV{\mathcal{V}}

\newcommand\cW{\mathcal{W}}
\newcommand\cA{\mathcal{A}}

\renewcommand\u{\mathbf{u}}
\renewcommand\v{\mathbf{v}}
\newcommand\x{\mathbf{x}}
\newcommand\y{\mathbf{y}}
\newcommand\z{\mathbf{z}}
\newcommand\hh{\mathbf{h}}
\newcommand\bN{\mathbf{N}}

\renewcommand\d{\,\mathrm{d}}

\renewcommand{\le}{\leqslant}
\renewcommand{\ge}{\geqslant}
\renewcommand{\leq}{\leqslant}
\renewcommand{\geq}{\geqslant}

\newcommand{\bmu}{\boldsymbol{\mu}}

\newcommand\ve{\varepsilon}
\newcommand\D{\Delta}

\renewcommand\phi{\varphi}

\newcommand\la{\lambda}
\newcommand\al{\alpha}
\newcommand\be{\beta}

\newcommand{\Kring}{\mathfrak{o}_K}
\newcommand{\Lring}{\mathfrak{o}_L}

\DeclareMathOperator{\Mod}{mod}
\renewcommand{\bmod}[1]{\,(\Mod{ #1})}

\newtheorem{theorem}{Theorem}
\newtheorem*{theorem*}{Theorem}

\newtheorem{lemma}{Lemma}

\theoremstyle{definition}
\newtheorem*{ack}{Acknowledgements}

\newcommand{\twosum}[2]{\sum_{\substack{#1\\#2}}}

\newcommand{\Q}{\mathbb{Q}}
\newcommand{\R}{\mathbb{R}}

\renewcommand{\b}[1]{{\bf #1}}
\newcommand{\cl}[1]{{\mathcal #1}}
\newcommand{\Z}{\mathbb{Z}}

\renewcommand{\mod}[1]{\hspace{-2.9mm}\pmod{#1}}

\DeclareMathOperator{\tr}{Tr}
\newcommand{\str}{\widetilde{{\tr}}}

\newcommand\w{\mathbf{w}}

\renewcommand\t{\mathbf{t}}

\newcommand\del{\delta}

\newcommand{\N}{\mathbf{N}}
\DeclareMathOperator{\meas}{meas}

\DeclareMathOperator*{\Osum}{\sum{}^*}
\newcommand{\sumE}{\sideset{}{^{(E)}}\sum}
\DeclareMathOperator{\osum}{\sum{}^*}

\begin{document}

\title{Quadratic polynomials represented by norm forms}

\author{T.D. Browning}
\address{School of Mathematics\\
University of Bristol\\ Bristol\\ BS8 1TW
\\ United Kingdom}
\email{t.d.browning@bristol.ac.uk}
\author{D.R. Heath-Brown}
\address{Mathematical Institute\\ University of Oxford\\ Oxford\\ OX1
  3LB
\\ United Kingdom}
\email{rhb@maths.ox.ac.uk}

\date{\today}

\begin{abstract}
Let $P(t)\in\QQ[t]$ be an irreducible quadratic polynomial and suppose
that $K$ is a quartic extension of $\QQ$ containing the roots of 
$P(t)$.  Let
$\bN_{K/\QQ}(\x)$ be a full norm form for the extension
$K/\QQ$. We show that the variety
\[P(t)=\bN_{K/\QQ}(\x)\neq 0\]
satisfies the Hasse principle and weak approximation.
The proof uses analytic methods.
\end{abstract}


\maketitle
\tableofcontents

\section{Introduction}\label{intro}

Let $k$ be a number field with set of valuations $\Omega_k$. 
Given an algebraic variety $X$ defined over $k$ we have the obvious
inclusions
$$
X(k)\xrightarrow{\Delta} X(\mathbb{A}_k)\subseteq 
\prod_{\nu \in \Omega_k} X(k_\nu),
$$
where $\Delta$ is the diagonal embedding of the set $X(k)$ of
$k$-rational points into the set 
$X(\mathbb{A}_k)$ of ad\`eles of $X$. Moreover 
the set $X(\mathbb{A}_k)$ is empty if and only if $\prod_{\nu \in
  \Omega_k} X(k_\nu)$ is empty and clearly provides a local obstruction to
the existence of $k$-rational points on $X$. Recall that 
a class $\mathcal{X}$ of algebraic varieties $X$ defined over
$k$ is said to satisfy the Hasse principle if $X(k)\neq \emptyset$ whenever
$X(\mathbb{A}_k)\neq \emptyset$. Likewise $\mathcal{X}$
is said to satisfy weak approximation if whenever it is non-empty the image
of $X(k)$ under $\Delta$ is dense in $X(\mathbb{A}_k)$ in the
product of $\nu$-adic topologies.

This paper is concerned
with the Hasse principle and weak approximation for
the class of varieties satisfying the Diophantine equation
\begin{equation} \label{PNZ'}
P(t)=\bN_{K/k}(x_1,\ldots,x_n)\neq 0,
\end{equation}
where $\bN_{K/k}$ is a full norm form for an extension $K/k$ of number
fields, and $P(t)$ is a polynomial over $k$.  Thus if $[K:k]=n$ and
we fix a basis $\{\omega_1,\ldots,\omega_n\}$ for $K$ as a vector space
over $k$, then 
$$
\bN_{K/k}(x_1,\ldots,x_n):=
N_{K/k}(x_1\omega_1+\cdots+x_n\omega_n).
$$ 
Throughout this paper we
will use $\bN_{K/k}$ to denote a norm form, and $N_{K/k}$ to denote
the corresponding field norm.

Progress on this problem has been limited, and we begin by discussing
what is known in the simplest cases.
A crude measure of difficulty is given by the number of distinct roots of $P(t)$  
over an algebraic closure $\bar{k}$. 
When $P(t)$ is a non-zero constant polynomial the Hasse principle for 
\eqref{PNZ'}
is known as the ``Hasse norm principle''.  The validity of the Hasse
norm principle for cyclic extensions $K/k$ was established by
Hasse himself, but for non-cyclic extensions there can be counterexamples.
There is an extensive literature on the subject and it is known, for
example, that the Hasse norm principle holds if the field $K$ has
prime degree over $k$ (Bartels \cite{bartels-81}, for example); or the Galois group
of $N/k$ is dihedral, where $N$ is the normal closure of $K$ over $k$
(Bartels \cite{bartels-81'}); or the extension $K/k$ is Galois and every 
Sylow subgroup of the Galois group is cyclic (Gurak \cite[Corollary 
3.2]{gurak});.

Following the work of 
Colliot-Th\'el\`ene and Sansuc \cite[Proposition 9.1]{CTSa}, 
we also have simple sufficient conditions to ensure that ``weak
approximation for norms'' holds, by which we mean that 
weak approximation holds for
\eqref{PNZ'} when $P(t)$ is a non-zero constant polynomial. 
Let $N$ be the normal closure of $K$ over  $k$. Then weak
approximation holds if either the degree $[K:k]$ is prime, or if
the Galois group of $N/k$ has cyclic Sylow subgroups.
In particular the latter result 
implies that it suffices for $K/k$ to be cyclic. 

The next case to consider is that in which $P(t)=ct^d$ for some $c\in
k^\times$ and some positive degree $d$. In this situation the Hasse
principle and weak approximation may fail.  However \eqref{PNZ'} 
is a principal 
homogeneous space under an algebraic $k$-torus, and 
the work of 
Sansuc \cite{sansuc} 
and  
Voskresenski$\breve{\i}$ \cite{vosk}  
shows that the
Brauer--Manin obstruction is the only obstruction to the Hasse
principle and weak approximation on any smooth projective model of
this variety. 

When $[K:k]=2$ and $P(t)$ has degree $3$ or $4$ then \eqref{PNZ'}
defines a Ch\^atelet surface. The arithmetic of such surfaces is
well-understood.  The Hasse principle and weak approximation may fail,
but it has been shown by
Colliot-Th\'el\`ene, Sansuc and Swinnerton-Dyer \cite{crelle}
that all such failures are explained by the 
Brauer--Manin obstruction.
The same conclusion is  available when $[K:k]=3$ and 
$\deg P(t)\leq 3$, by work of 
Colliot-Th\'el\`ene and Salberger \cite{ct-salb}.  

There have also been investigations into \eqref{PNZ'} when  $P(t)$
factors completely over $k$, with at most two roots. In this case one
may write
$$
P(t)=c(t-a)^u(t-b)^v,
$$
with $a,b,c\in k$ and $u,v\in\NN$.  
It is known that one has the Hasse principle and weak 
approximation whenever the  
Brauer--Manin obstruction is empty, providing that we work over the 
ground field $k=\Q$. This was first proved
under the assumption that $\gcd(u,v,n)=1$,
by Heath-Brown and Skorobogatov \cite{HBS}, a condition that was
subsequently removed by Colliot-Th\'el\`ene, Harari and Skorobogatov 
\cite{CTHS}.  While all the previous work described had been purely 
algebraic, the approach used by
Heath-Brown and Skorobogatov combined analytic machinery, in the form
of the Hardy--Littlewood circle
method, with the previous descent approach to the Brauer--Manin
obstruction.  The circle method can be adapted, with some effort, 
to apply to ground fields other than $k=\Q$.  However, for simplicity,
this possibility was not explored in \cite{HBS}.

Very little is known about other polynomials $P(t)$.
When $P(t)$ is a non-zero separable polynomial with degree at least $2$, 
it is conjectured that the Hasse principle and weak approximation hold 
whenever the 
Brauer--Manin obstruction is empty. 
When $K/k$ is cyclic and Schinzel's Hypothesis is granted, work of  
Colliot-Th\'el\`ene, Skorobogatov and Swinnerton-Dyer
\cite[Theorem~1.1]{98a}
yields a positive answer to this question. 
Note that 
this result is already a special case of earlier work of Colliot-Th\'el\`ene
and Swinnerton-Dyer \cite{CT-94} on pencils of Severi--Brauer
varieties, but this connection is only made
clear in the discussion  \cite[page 10]{98a}.
Note that when $K/k$ is cyclic the Brauer--Manin obstruction is known 
to be empty if $P(t)$ is 
irreducible over $k$. (This follows from Corollary 2.6(c) of 
Colliot-Th\'el\`ene, Harari and Skorobogatov \cite{CTHS}, which shows
that the Brauer group contains only vertical elements when $K/k$ is
cyclic.  However it is not hard to show that the vertical part of the
Brauer group is trivial when $P(t)$ is irreducible over $k$, using 
the remark on page 76 of \cite{CTHS}.)

Finally we mention that there is potential for tackling the case in
which $P(t)$ is an arbitrary polynomial which splits completely over $k$ into $d$ linear factors, 
at
least in the case $k=\QQ$, by using ideas from the work of Green and
Tao \cite{GT} together with the main theorem from Green, Tao and
Ziegler \cite{GTZ}.  In this case the methods of Heath-Brown and
Skorobogatov \cite{HBS} reduce the problem \eqref{PNZ'} to one
involving a system of equations
\[\bN_{K/\QQ}(\x_i)+a_i\bN_{K/\QQ}(\x_0)=c_i\bN_{K/\QQ}(\y_i)\neq 0,
\quad (1\le i\le d-1).\]
In the language of \cite{GT}, this is a system of linear forms of 
finite complexity.  The machinery of Green and Tao, and of Green, Tao and
Ziegler, allows one to handle such systems when the norms are replaced
by primes, and it seems reasonable to hope that a variant of the 
method would allow one to handle the problem above.
This plan was first mentioned to us by Professor Wooley.

It will be apparent from the foregoing survey that the most obvious
open case is that in which $P(t)$ is an irreducible quadratic
over $k$, and this is the goal of the present paper.
Again we shall be dependent on techniques from analytic
number theory which have not been fully developed for ground fields
other than $k=\Q$, so we shall confine attention to this
latter case. With this restriction our goal will be to establish the Hasse 
principle and weak approximation for
\begin{equation} \label{PNZ}
P(t)=\bN_{K/\QQ}(x_1,\ldots,x_n)\not=0,
\end{equation}
under suitable assumptions on the extension $K/\QQ$. 
Let  $|\cdot|_\nu$ denote the $\nu$-adic norm, which we extend to 
vectors by setting 
$|\z|_\nu:=\max_{1\leq i\leq m}|z_i|_\nu, $
if $\z=(z_1,\ldots,z_m)$. When $\nu=\infty$ we will simply write
$|\cdot|_\infty=|\cdot|.$   
With this in mind the following is our main result.

\begin{theorem}\label{thmain}
Let $P(t)\in\QQ[t]$ be an irreducible quadratic polynomial and 
let $K$ be a quartic extension of $\QQ$ containing a root
of $P(t)$. 
Suppose that, for every $\nu\in \Omega_\QQ$, 
we are given a solution $(t^{(\nu)},\x^{(\nu)})\in
\QQ_\nu^{5}$ of \eqref{PNZ}. Let $S\subset \Omega_\QQ$ be any finite 
subset and let $\ve>0$. Then
there is a solution $(t,\x)\in \QQ^{5}$ of \eqref{PNZ} such that 
\begin{equation}
  \label{tap}
  |t-t^{(\nu)}|_\nu<\ve, \quad 
  |\x-\x^{(\nu)}|_\nu<\ve, 
\end{equation}
for every $\nu \in S$.  Thus the Hasse principle and weak
approximation hold for  \eqref{PNZ}. 
\end{theorem}

It is interesting to note that our result is both unconditional and
concerns field extensions $K/\QQ$ which may be non-cyclic. This marks 
a departure from the sort of results achieved in \cite{98a}.
In fact our theorem answers in the affirmative a question posed by 
Colliot-Th\'el\`ene, during the 2005 Bremen workshop
``Rational points on curves --- explicit methods'',
about the Hasse
principle for \eqref{PNZ} in the special case that $K$ is a
biquadratic extension containing a root of $P(t)$ (cf.\  
the questions at the close of \S 2 in the work of  
Colliot-Th\'el\`ene, Harari and Skorobogatov  
\cite{CTHS}).

The proof of Theorem \ref{thmain} relies on techniques from analytic
number theory and is 
inspired by work of Fouvry and Iwaniec \cite{FI}, who proved
that there are infinitely many primes $p$ of the form $a^2+q^2$, with
$q$ also prime.  More generally they showed how to produce primes of
the form $a^2+q^2$ with $q$ from any sufficiently dense set.  Our
argument involves many complexities of detail, but also one major
simplification, since we have only to produce integers in
$N_{K/\QQ}(K^{\times})$, rather than primes.

We can generalise our result mildly,
to include the case in which $P(t)=cQ(t)^u$ for an odd positive integer $u$, 
 where
$Q(t)\in \QQ[t]$ is an irreducible quadratic polynomial. 
This is achieved by establishing a bijection between
solutions  of \eqref{PNZ} and solutions of the corresponding
equation in which $P(t)$ is quadratic. To do this we begin by
choosing $e,f\in\ZZ$ for which $eu+4f=1$.  The equation \eqref{PNZ}
becomes $cQ(t)^u=\bN_{K/\QQ}(\x)$, and raising to the power $e$ we
obtain $c^eQ(t)^{1-4f}=\bN_{K/\QQ}(\x)^e$, whence
$c^eQ(t)=\bN_{K/\QQ}(\x)^eQ(t)^{4f}$. Since $K$ has degree $4$ over
$\QQ$, we deduce that $P_0(t):=c^eQ(t)$ is a norm from $K$ whenever
$P(t)=cQ(t)^u$ is a norm from $K$. The converse deduction is similar.
Thus if $c^eQ(t)=\bN_{K/\QQ}(\x)$ then, raising both sides to the
power $u$ we find that $c^{1-4f}Q(t)^u=\bN_{K/\QQ}(\x)^u$, whence
$cQ(t)^u=c^{4f}\bN_{K/\QQ}(\x)^u$. Thus $P(t)$ is a norm from $K$
whenever $P_0(t)$ is a norm from $K$. It is now immediate from
Theorem \ref{thmain} that we have the Hasse principle and weak
approximation for $P(t)=cQ(t)^u$.   

By a simple change of
variable we may assume that $P(t)=c(1-at^2)$ in Theorem
\ref{thmain},  where $c$ is a non-zero rational and $a$ is a 
square-free integer. 
Let 
$$
L:=\QQ(\sqrt{a}).
$$
The fields to which our theorem applies 
take the shape $K=L(\sqrt{\beta})$, with $\beta\in L$. 
In particular we have $L\subseteq K$ in 
the statement of Theorem \ref{thmain}. It turns out that most of our
argument carries over to an arbitrary degree $n$ extension 
$K$ of $\QQ$ that contains $L$ as a subfield. 
Given  $P(t)=c(1-at^2)$ as above, we suppose that
$(t^{(\nu)},\x^{(\nu)})\in  
\QQ_\nu^{n+1}$  are 
solutions of \eqref{PNZ}, for each  $\nu\in \Omega_\QQ$. 
Then we
want to determine conditions on $K$, beyond the hypothesis $L\subseteq
K$,  such that for any finite set $S\subset \Omega_\QQ$ 
we can 
find a solution $(t,\x)\in \QQ^{n+1}$ of \eqref{PNZ} for which the
weak approximation condition \eqref{tap} holds for each $\nu \in S$.

In pursuing this goal we may assume that
$\{\omega_1,\ldots,\omega_n\}$ is an
integral basis for the ring of integers
$\Kring$, with $\omega_1=1$. 
By the transitivity 
of norms we have $N_{K/\QQ}=N_{L/\QQ}\circ N_{K/L}$, since $L\subseteq
K$.
Hence any norm from $K$ to $\QQ$ is also a norm from $L$ to $\QQ$. 
We will make frequent use of this fact in our work.

Since $(1-at^2)=N_{L/\QQ}(1+t\sqrt{a})$ it  follows from
the hypotheses of the theorem that the equation
$c=N_{L/\QQ}(u+v\sqrt{a})$ can be solved for $u,v\in \QQ_\nu$ for any
$\nu \in \Omega_\QQ$.  The Hasse norm principle therefore implies that there
exists $\delta\in L^\times$ such that 
$$
c=N_{L/\QQ}(\delta)^{-1}.
$$
Thus it will suffice to work with the equation
\begin{equation}
  \label{eq:2.3}
1-at^2=N_{L/\QQ}(\del)\bN_{K/\QQ}(\x)\neq 0,
\end{equation}
rather than \eqref{PNZ},  with square-free $a\in \ZZ$ and 
non-zero $\delta \in L$.
We are then given a finite set $S\subset \Omega_\QQ$ and a solution
$(t^{(\nu)},\x^{(\nu)})\in 
\QQ_\nu^{n+1}$ of this equation  for every $\nu \in \Omega_\QQ$ and we
wish to establish the existence of a 
solution $(t,\x)\in \QQ^{n+1}$ such that \eqref{tap} holds for every
$\nu \in S$.
Our plan is to achieve this by arranging that
\begin{equation}\label{tdef-b}
\bN_{K/\QQ}(\w)(1+t\sqrt{a})=\del\bN_{K/L}(\y)\not=0,
\end{equation}
where
\[
\bN_{K/L}(y_1,\ldots,y_n):=N_{K/L}(y_1\omega_1+\cdots+y_n\omega_n).
\]
Then if $\beta_1=y_1\omega_1+\cdots+ y_n\omega_n$ and
$\beta_2=w_1\omega_1+\cdots+ w_n\omega_n$ we obtain (\ref{eq:2.3}) on taking
$x_1\omega_1+\cdots +x_n\omega_n$ to correspond to the element 
$\beta=\beta_1\beta_2^{-2}$. It will be convenient to write $\x=\y.\w^{-2}$
for the vector $\x$ produced by this construction.

In fact we will establish the following result,
which demonstrates the Hasse principle and weak approximation
for \eqref{tdef-b}, for any number field $K$ of degree $n$ that
contains $L$. 

\begin{theorem}\label{thmain-a}
Let $\delta \in L^\times$ and assume that  $L\subseteq K$. 
Suppose that, for every $\nu\in \Omega_\QQ$, 
we are given a solution $(t^{(\nu)},\w^{(\nu)},\y^{(\nu)})\in
\QQ_\nu^{2n+1}$ of \eqref{tdef-b}. Let $S\subset \Omega_\QQ$ be any
finite  
subset and let $\ve>0$. Then
there is a solution $(t,\w,\y)\in \QQ^{2n+1}$ of \eqref{tdef-b} 
such that \[  |t-t^{(\nu)}|_\nu<\ve, \quad  |\w-\w^{(\nu)}|_\nu<\ve, \quad 
  |\y-\y^{(\nu)}|_\nu<\ve,\]
for every $\nu \in S$.  
\end{theorem}

In \S~\ref{s:deduc} we will show how Theorem \ref{thmain} follows
from this result when $K$ is a quadratic extension
of  $L$. This will be achieved via the following result.

\begin{lemma}\label{newlem}
Let $K$ be a quadratic extension of $L$. Let
$S\subset\Omega_{\QQ}$ be a finite set and let $\ve>0$ be given.
Suppose that the equation
\begin{equation}\label{PNZ1}
c(1-at^2)=\bN_{K/\QQ}(\x)\neq 0
\end{equation}
has solutions $(t^{(\nu)},\x^{(\nu)})$ everywhere locally.  Then there 
exists $\del=\del_{\ve}\in L$ with $c=N_{L/\QQ}(\del)^{-1}$ such that
\begin{equation}\label{tdef-b'}
1+t\sqrt{a}=\del\bN_{K/L}(\x)\not=0
\end{equation}
has solutions $(t_0^{(\nu)},\x_0^{(\nu)})$ everywhere locally,
with
\[|t^{(\nu)}-t_0^{(\nu)}|_{\nu}<\ve,\quad 
|\x^{(\nu)}-\x_0^{(\nu)}|_{\nu}<\ve,
\]
for every $\nu\in S$.
\end{lemma}

We should emphasise here that when we speak of local solutions 
we are thinking of zeros over $\QQ_{\nu}$ of the polynomial,
defined over $\QQ_{\nu}$, which specifies the equation.  In
particular, elements of the completions $L_{\mu}$, for $\mu\mid\nu$, 
do not occur.

It has been suggested to us by Professor Colliot-Th\'el\`ene that the 
open descent method of Colliot-Th\'el\`ene and 
Skorobogatov \cite{ct-skoro} might be used to establish a variant of
Lemma \ref{newlem} in which $K$ is an arbitrary finite extension of $L$.
The proposed lemma would then give the same conclusion as Lemma
\ref{newlem}, but under the assumption that the solutions 
$(t^{(\nu)},\x^{(\nu)})$ of (\ref{PNZ1}) produce an ad\`elic point
orthogonal to the Brauer group of the variety. 
Once combined with 
Theorem \ref{thmain-a}, this should demonstrate that the 
Brauer--Manin obstruction 
to the Hasse principle and to weak approximation is the only one
for \eqref{PNZ}, when $P(t)\in\QQ[t]$ is an irreducible quadratic 
polynomial and 
$K$ is an arbitrary extension of $\QQ$ containing a root of $P(t)$.
However this would still leave open the difficult problem of
calculating the Brauer group, which our route avoids.

By Lemma \ref{newlem}, 
given local solutions of \eqref{PNZ1}, we may produce an equation
\eqref{tdef-b} in which $\w=(1,0,0,\ldots,0)$ and $\y=\x$, and 
which has corresponding local solutions suitably close
to those of \eqref{PNZ1}.  We may then use Theorem \ref{thmain-a} to
produce a global solution of \eqref{tdef-b} close to the given local
solutions. Finally, taking the norm from $L$ to $\QQ$ we obtain a
suitable global solution of \eqref{PNZ1}.  It should be pointed out
that this argument uses the fact that the map from $(\w,\y)$ to
$\x=\y.\w^{-2}$ is continuous for $|\cdot |_{\nu}$ providing that we
avoid a neighbourhood of $\w=\mathbf{0}$.
We stress that the only point in the paper where we use our
assumption 
that $K/L$ is quadratic occurs in the proof of Lemma \ref{newlem}.

We proceed to indicate the initial steps in our treatment of
\eqref{tdef-b} in Theorem \ref{thmain-a}. 
A suitable value of $t\in\QQ$ will exist, providing that
\begin{equation}\label{yeq1}
\tr_{L/\QQ}\big(\delta \bN_{K/L}(\y)\big)=2\bN_{K/\QQ}(\w)\neq 0.
\end{equation}
Moreover, we will then have
\[
t=\frac{1}{2}\tr_{L/\QQ}\left(\frac{
\bN_{K/\QQ}(\w)^{-1}\delta \bN_{K/L}(\y)}{\sqrt{a}}\right).
\]
Thus if we have a solution of (\ref{yeq1}) in which $\y$ and $\w$ are
sufficiently close to $\y^{(\nu)}$ and $\w^{(\nu)}$ it will be
automatic that the corresponding solution $t$ will be close to
$t^{(\nu)}$. It follows that we have only to establish a suitable
Hasse principle and weak approximation result for (\ref{yeq1}).

We must make a further manoeuvre before reaching our fundamental
equation.  As mentioned above, the work of Fouvry and Iwaniec handles
primes.  Using standard machinery, prime numbers are dealt with by
means of ``Type~I sums'' and ``Type II sums''. Of these, Type II sums
involve bilinear forms in which the prime number is replaced by a
product of integers $uv$, which have to lie in suitable ranges.  In
our case we can insist that our norm $\bN_{K/L}(\y)$ is a product
$\bN_{K/L}(\u)\bN_{K/L}(\v)$, thereby eliminating the need to
consider Type I sums.  Indeed, since we can specify the sizes of $\u$
and $\v$, the treatment of the Type II sums will also be simplified
somewhat. Thus instead of attacking \eqref{yeq1} we shall consider 
the Diophantine equation
$$
\tr_{L/\QQ}\big(\delta \bN_{K/L}(\u)\bN_{K/L}(\v)\big)=
2\bN_{K/\QQ}(\w)\neq 0,
$$
with the aim of finding suitably
localised solutions 
\[(\u, \v, \w)\in \ZZ^n\times\ZZ^n\times\ZZ^n.\]

Let $\sigma$ denote the non-trivial automorphism of $L$
and suppose that $\{1,\tau\}$ is a $\ZZ$-basis for $\Lring$, and hence
also a $\QQ$-basis for $L$. 
For technical reasons it will be convenient to replace the
trace $\tr_{L/\QQ}$ by a ``skew-trace'' 
$$
\str(x,y):=\tr_{L/\QQ}(xy^\sigma D_L^{-1})
$$
for $x,y\in L$, where $D_L=\tau-\tau^{\sigma}$.
Thus $(D_L)$ is the different of $L/\QQ$. On writing
\[x=\delta \bN_{K/L}(\u),\quad 
y=\left(\bN_{K/L}(\v)D_L \right)^\sigma\]
our condition becomes
\begin{equation}
  \label{TR2}
\str(x,y)=2\bN_{K/\QQ}(\w)\neq 0.
\end{equation}
We will count suitably restricted solutions of this equation.  If 
$\mathcal{N}$ is the number of such solutions we can write
$$
\mathcal{N}=\sum_{x \in \Lring }\sum_{y \in \Lring} 
\al(x)\be(y)\la\big(\str(x,y)\big).
$$
Here the function $\al(x)$, respectively $\be(y)$, 
counts appropriately restricted representations of $x$ by $\delta
\N_{K/L}(\u)$, respectively of  
$y$ by 
$\left(\bN_{K/L}(\v)D_L \right)^\sigma$.
Moreover $\la(l)$ counts suitably constrained solutions of
$l=2\N_{K/\QQ}(\w)$.

Our expression for $\cl{N}$ can be viewed as a bilinear form.  We
have good techniques for estimating these, going back to the
works of Vinogradov.  However the methods are designed to produce upper
bounds for bilinear forms in which we expect cancellation, while our
problem is to establish an asymptotic formula for an expression in
which all the terms are non-negative.  We shall therefore split
$\alpha(x)$ into two parts $\alpha(x)=\hat{\alpha}(x)+\alpha_0(x)$,
and write 
\[\cl{N}=\cl{M}+\cl{E},\]
where $\cl{M}$ is a main term, and contains the contribution from 
$\hat{\alpha}(x)$, while $\cl{E}$ is an error term, and is the 
corresponding expression involving $\alpha_0(x)$.  
Thus we will need $\hat{\alpha}(x)$ to be a
sufficiently simple function that we can compute $\cl{M}$ directly.
Moreover we will want $\alpha_0(x)$ to produce sufficient cancellation on
average, so that a bilinear form estimation of $\cl{E}$ can be achieved.
The underlying principle here is exactly that which Linnik
\cite{linnik} developed
in his ``dispersion method''.

In \S~\ref{s:approx} we will describe a general procedure for producing an
approximation of the type $\hat{\alpha}(x)$.  In our context $x$ runs
over the ring $\Lring$, but we will begin by presenting the method as
it applies to sequences indexed by $\ZZ$, since we hope this will prove
to be of independent interest. 
With this choice made, 
our proof that the bilinear form
$\mathcal{E}$ 
makes a satisfactory overall contribution to the asymptotic formula is
the subject of \S~\ref{s:bil}.
It is this part of our argument which is based on ideas from the 
work of Fouvry and Iwaniec \cite{FI}, who provide a general
framework for estimating sums of this sort.
Finally, the asymptotic evaluation of $\mathcal{M}$ will
be executed in \S~\ref{s:M} and \S~\ref{s:ssi}.

\begin{ack}
Some of this work was done while the authors were visiting the 
{\em Hausdorff Institute} in Bonn and the {\em Institute
for Advanced Study} in Princeton, and while the second author was
visiting the {\em Mathematical Sciences Research Institute} in
Berkeley. 
The hospitality and financial
support of these bodies is gratefully acknowledged. 
While working on this paper the first author was 
supported by EPSRC grant number
\texttt{EP/E053262/1}.
\end{ack}

\section{Deduction of Theorem \ref{thmain}}\label{s:deduc}

The goal of this section is to prove Lemma \ref{newlem}.  As explained
in the previous section, this is enough to allow the deduction of
Theorem \ref{thmain} from Theorem~\ref{thmain-a}. 
Let $P(t)=c(1-at^2)$,  where $c$ is a non-zero rational and $a$ is a 
square-free integer.   
For the moment let $K$ be an arbitrary number field of degree $n$
containing $L=\QQ(\sqrt{a})$.  In particular $n$ is even.
Fix any $\delta_0 \in L^\times$ such that $c=N_{L/\QQ}(\delta_0)^{-1}$ and let
$S\subset \Omega_\QQ$ be finite. We may 
assume that $S$ contains the archimedean 
valuation, together with any valuations that become ramified in 
$K$, and any valuations $\nu \in \Omega_\QQ$ for which $v_\mu(\delta_0)\neq 0$
for some $\mu\in \Omega_L$ above $\nu$.

We proceed to look for a suitable $\del$, fulfilling the conditions of
Lemma \ref{newlem}, by examining values 
$\del=\delta_0\gamma^\sigma \gamma^{-1}$.
We will see that when $\nu\not\in S$ any $\gamma\in L^{\times}$ 
is acceptable.
Thus our first task is to find a value $\gamma$ which works for the ``bad''
places $\nu\in S$.  For such valuations we claim
that there exist $\gamma_1^{(\nu)},\gamma_2^{(\nu)} \in \QQ_\nu$,
with 
$$
\gamma^{(\nu)}=
\gamma_1^{(\nu)}+\gamma_2^{(\nu)}\sqrt{a}\neq 0,
$$ 
such that 
$$
1+t^{(\nu)}\sqrt{a}=\delta_0{\gamma^{(\nu)}}^\sigma
{\gamma^{(\nu)}}^{-1}\bN_{K/L}(\x^{(\nu)}).
$$
It then suffices to choose $\gamma \in L^\times$ so that $\gamma$ is close to
$\gamma^{(\nu)}$ for each $\nu \in S$, since then we may take
$t_0^{(\nu)}=t^{(\nu)}$ and find a suitable $\x_0^{(\nu)}$ close to
$\x^{(\nu)}$.

To establish the claim we begin by noting that the form $\bN_{K/L}$ 
decomposes as
\begin{equation}\label{nsplit}
\bN_{K/L}=N_1+N_2\sqrt{a},
\end{equation}
over $L$, where $N_1,N_2 \in \QQ[x_1,\ldots,x_n]$ are forms of 
degree $n/2$. Setting $\delta_0=\delta_1+\delta_2\sqrt{a}$
and $\gamma^{(\nu)}=c_1+c_2\sqrt{a}$, and multiplying through by
$\gamma^{(\nu)}$, we see that the equation becomes
\begin{equation}\label{teq}
\begin{split}
(c_1+c_2\sqrt{a})&(1+t^{(\nu)}\sqrt{a})\\
&=(c_1-c_2\sqrt{a})
(\delta_1+\delta_2\sqrt{a})(N_1(\x^{(\nu)})+N_2(\x^{(\nu)})\sqrt{a}).
\end{split}
\end{equation}
Thus our problem is to show the existence of $(c_1,c_2)\in\QQ_{\nu}^2$ 
satisfying this, given
the condition \eqref{eq:2.3}, namely 
\begin{equation}\label{add}
1-a{t^{(\nu)}}^2=
(\delta_1^2-a\delta_2^2)(N_1(\x^{(\nu)})^2-aN_2(\x^{(\nu)})^2)\not=0.
\end{equation}
If we set $A_1=\delta_1N_1(\x^{(\nu)})+a\delta_2N_2(\x^{(\nu)})$ and 
$A_2=\delta_1N_2(\x^{(\nu)})+\delta_2N_1(\x^{(\nu)})$
for convenience, then 
\eqref{teq} becomes a pair of conditions
\begin{equation}\label{sim}
c_1(1-A_1)+c_2(at^{(\nu)}+aA_2)=c_1(t^{(\nu)}-A_2)+c_2(A_1+1)=0.
\end{equation}
We need to find a solution $c_1,c_2$ of these, 
with $c_1^2-ac_2^2\not=0$.  In doing so we may assume \eqref{add}, 
which becomes
\[
1-a{t^{(\nu)}}^2=A_1^2-aA_2^2\not=0.
\]
However the determinant of the system \eqref{sim} is
\begin{align*}
(1-A_1)(A_1+1)-(at^{(\nu)}+aA_2)(t^{(\nu)}-A_2)
&=(1-a{t^{(\nu)}}^2)-(A_1^2-aA_2^2)\\
&=0.
\end{align*}
Moreover, if
$c_1=\pm\sqrt{a}c_2\not=0$ one readily deduces from \eqref{teq} that
$1-a{t^{(\nu)}}^2=0$, which is impossible. This suffices for
the proof of the claim.

In handling the case $\nu\not\in S$ the following lemma will be
useful.  It will be proved at the end of the section.

\begin{lemma}\label{newlem2}
Let $\nu\in\Omega_{\QQ}$ be a finite place, unramified in $K$.  
Suppose that
$$
\beta=b_1+b_2\sqrt{a}\in L^{\times}
$$ 
is a unit in $L_{\mu}$ for each 
place $\mu$ of $L$ above $\nu$. Then $\beta=\bN_{K/L}(\x)$ for some 
$\x\in\QQ_{\nu}^n$, by which we mean that if $N_1$ and $N_2$ are as in
\eqref{nsplit} then $b_1=N_1(\x)$ and $b_2=N_2(\x)$.
\end{lemma}

While Lemma \ref{newlem2} is valid for arbitrary 
extensions $K$ of $L$, we now  
make the assumption that $K$ is a quadratic extension of 
$L$. In particular we have $n=4$ in the above discussion.  

When $\nu\not\in S$ our task is to show that if
\[
f(t):=(1+t\sqrt{a})\gamma(\gamma^{\sigma})^{-1}\del_0^{-1},
\]
then there exists $t_0^{(\nu)}\in\QQ_{\nu}$ such that $f(t_0^{(\nu)})$
is of the form $\bN_{K/L}(\x_0^{(\nu)})$ for some vector 
$\x_0^{(\nu)}\in\QQ_{\nu}^4$.  Since $\nu\not\in S$ 
there are no weak approximation conditions to be satisfied.
We begin by considering the case in which $\nu\not\in S$ is inert in
$L/\QQ$. 
Then $v_\nu(\gamma^\sigma \gamma^{-1})=0$ for any $\gamma$, and
$v_{\nu}(\del_0)=0$ by the choice of $S$, whence
$\gamma(\gamma^{\sigma})^{-1}\del_0^{-1}$ will be a unit in
$L_\nu$. Lemma \ref{newlem2} 
then shows that there exists
$\x_0^{(\nu)}\in \QQ_\nu^4$ such that $f(0)=\bN_{K/L}(\x_0^{(\nu)})$.
It therefore suffices to take $t_0^{(\nu)}=0$ in \eqref{tdef-b'}.

Finally we must deal with the case $\nu \not \in S$, with $\nu$ split
in $L$. 
Suppose $\nu$ splits as $\mu_1$ and $\mu_2=\mu_1^\sigma$ in $L$. 
Write $p$ for the rational prime associated to $\nu$.
Let
\[
v_{\mu_1}(\gamma)=e_1,\quad
v_{\mu_2}(\gamma)=e_2.
\]
Let $\mathfrak{p}_1$ and $\mathfrak{p}_2=\mathfrak{p}_1^\sigma$ be the 
prime ideals associated to $\mu_1$ and $\mu_2$ respectively, so that 
$
(\gamma)=\mathfrak{p}_1^{e_1} \mathfrak{p}_2^{e_2} \mathfrak{g}
$ for some ideal $\mathfrak{g}$ coprime to both $\mathfrak{p}_1$ and
$\mathfrak{p}_2$.  Choose $\alpha_1\in\mathfrak{p}_1\setminus
(\mathfrak{p}_1^2\cup\mathfrak{p}_2)$, and set
$\alpha_2=\alpha_1^\sigma$. 
We write $e=2+|e_1|+|e_2|$ and 
\[
h_1=1+\tfrac12(|e_1|+|e_2|+e_2-e_1),
\quad
h_2=1+\tfrac12(|e_1|+|e_2|+e_1-e_2),
\]
so that $e, h_1$ and $h_2$ are integers, and $e$ is strictly positive.
Since $\alpha_i\in L$ and $K/L$ is quadratic, we have 
$N_{K/L}(\alpha_i)=\alpha_i^2$.
We may then calculate that 
\[
v_{\mu_1}(1+p^{-e}\sqrt{a})=-e, 
\quad v_{\mu_1}(\gamma)=e_1,\quad
v_{\mu_1}(\gamma^\sigma)=e_2,
\]
and
\[v_{\mu_1}\left(N_{K/L}(\alpha_1^{h_1}\alpha_2^{h_2})\right)=
2h_1v_{\mu_1}(\alpha_1)+2h_2v_{\mu_1}(\alpha_2)=2h_1,\]
whence
\[v_{\mu_1}\left(f(p^{-e})N_{K/L}(\alpha_1^{h_1}\alpha_2^{h_2})\right)=
-e+e_1-e_2+2h_1=0.\]
Similarly we have
\[v_{\mu_2}\left(f(p^{-e})N_{K/L}(\alpha_1^{h_1}\alpha_2^{h_2})\right)=0.\]
Thus Lemma
\ref{newlem2} tells us that $f(p^{-e})
N_{K/L}(\alpha_1^{h_1}\alpha_2^{h_2})$ can be written as 
$\bN_{K/L}(\y)$ for some $\y\in\QQ_\nu^4$.  It follows that 
$f(p^{-e})$ takes the form $\bN_{K/L}(\x_0)$ for some $\x_0\in\QQ_\nu^4$, 
whence $t_0^{(\nu)}=p^{-e}$ is acceptable 
in \eqref{tdef-b'}.  This completes the proof of Lemma \ref{newlem}.

It now remains to establish Lemma \ref{newlem2}, for which we return to an
arbitrary extension $K$ of degree $n$ over $\QQ$, which contains 
$L$ as a subfield. 
It will be convenient
to use the notation $i_{\mu}$ for the embedding of $L$ into the
completion $L_{\mu}$, and similarly for valuations of $\QQ$ and
$K$. Thus our hypothesis is that $i_{\mu}(\beta)$ is a unit for every
$\mu\mid\nu$. For each such $\mu$ we choose a place 
$\lambda(\mu)\in\Omega_K$ such that
$\lambda(\mu)\mid\mu$.  Then the extension of local fields
$K_{\lambda(\mu)}/L_{\mu}$ is unramified, whence $i_{\mu}(\beta)$ must
be a norm from $K_{\lambda(\mu)}$ (see Gras \cite[Corollary 1.4.3, part
(ii), page 75]{gras}, for example).  For each $\mu$ above $\nu$ we may
therefore write 
\[
i_{\mu}(\beta)=N_{K_{\lambda(\mu)}/L_{\mu}}(y_{\lambda(\mu)}),
\]
for appropriate elements $y_{\lambda(\mu)}\in K_{\lambda(\mu)}$.  If 
$\lambda\in\Omega_K$ lies above $\mu$ but is different from 
$\lambda(\mu)$ we take $y_{\lambda}=i_{\lambda}(1)$, so that
\[i_{\mu}(1)=N_{K_{\lambda}/L_{\mu}}(y_{\lambda}).\]
We now use weak approximation to find elements $y^{(j)}\in K$ such
that
\[
\left|i_{\lambda}(y^{(j)})-y_{\lambda}\right|_{\lambda}<\frac{1}{j}
\]
for every $\lambda\in\Omega_K$ above $\nu$.  We note in particular that the
sequence $\left(y^{(j)}\right)$ converges with 
respect to each valuation $\lambda$
above $\nu$. We now have
\[
\lim_{j\rightarrow\infty} N_{K_{\lambda}/L_{\mu}}(i_{\lambda}(y^{(j)}))
=
\begin{cases}
i_{\mu}(\beta), & \mbox{if $\lambda=\lambda(\mu)$ for some
$\mu\mid\nu$,}\\ 
i_{\mu}(1), & \mbox{otherwise},
\end{cases}
\]
where the limit is with respect to $|\cdot|_{\mu}$.
We therefore conclude that
\[\prod_{\lambda\mid\mu}N_{K_{\lambda}/L_{\mu}}(i_{\lambda}(y^{(j)}))
\rightarrow i_{\mu}(\beta).\]
However, according to Gras \cite[Proposition 2.2, page 93]{gras} we have
\[
\prod_{\lambda\mid\mu}N_{K_{\lambda}/L_{\mu}}(i_{\lambda}(y^{(j)}))
=i_{\mu}(N_{K/L}(y^{(j)})),
\]
so that
\[
i_{\mu}(N_{K/L}(y^{(j)}))\rightarrow i_{\mu}(\beta).
\]
Since this holds for all $\mu$ above $\nu$ it follows that
\[N_1(\x^{(j)})\rightarrow b_1\;\;\;\mbox{and}\;\;\;
N_2(\x^{(j)})\rightarrow b_2,\]
the convergence being with respect to $|\cdot|_{\nu}$, 
where
\[
y^{(j)}=x_1^{(j)}\omega_1+\cdots+x_n^{(j)}\omega_n.
\]
Finally, since the sequence $y^{(j)}\in K$ converges for every
valuation $\lambda$ above $\nu$, the sequence $\x^{(j)}\in\QQ^n$ must
converge in $\Q_{\nu}$, yielding the required vector $\x\in\Q_{\nu}^n$.

\section{Preliminaries for Theorem \ref{thmain-a}}

We are now ready to begin the proof of Theorem \ref{thmain-a}, which
concerns  the Hasse principle and weak
approximation for the variety  \eqref{tdef-b}.
For the remainder of the paper let $K$ be an arbitrary
number field of degree $n$ that contains $L=\QQ(\sqrt{a})$.
As noted, it will be 
convenient to work with the equivalent variety \eqref{yeq1}.  For ease
of reference we repeat the definition here:-
\begin{equation}\label{yeq2}
\tr_{L/\QQ}\big(\delta \bN_{K/L}(\y)\big)=2\bN_{K/\QQ}(\w)\neq 0.
\end{equation}
We are presented with local solutions
$\y^{(\nu)},\w^{(\nu)}$ for every valuation $\nu$, and wish
to find a global solution which approximates these for every $\nu\in S$.
We claim that it suffices to consider the case
in which $\y^{(\nu)},\w^{(\nu)}$ are integral
for every finite $\nu\in S$.  Indeed, let us suppose that we can solve
\eqref{yeq2} with the local conditions 
\begin{equation}\label{tap31}
     |\y-\y^{(\nu)}|_\nu<\ve, \quad
  |\w-\w^{(\nu)}|_\nu<\ve, 
\end{equation}
providing that
$\y^{(\nu)},\w^{(\nu)}$ are integral for
all finite $\nu\in S$.  Let us suppose further that we have
a general set of values for 
$\y^{(\nu)},\w^{(\nu)}$ satisfying
\eqref{yeq2}.  
Then we may choose an integer $N\in\NN$ such that
$N^2\y^{(\nu)},N\w^{(\nu)}$ are integral for 
all finite $\nu\in S$. 
We note that these values will still satisfy
\eqref{yeq2}.  Then, by our assumption, we can find a global solution
$\y,\w$ of \eqref{yeq2} which satisfies
\[ 
    |\y-N^2\y^{(\nu)}|_\nu<\ve|N^2|_\nu, \quad
  |\w-N\w^{(\nu)}|_\nu<\ve|N|_\nu \]
for all $\nu\in S$.  It follows that $N^{-2}\y,
N^{-1}\w$ 
is a solution of \eqref{yeq2} fulfilling the
condition \eqref{tap31}.  This establishes our claim.

We now use weak approximation for $\ZZ^n$ to produce vectors
$\y^{(M)},\w^{(M)}\in\ZZ^n$ such that
\[  |\y^{(M)}-\y^{(\nu)}|_\nu<\ve, \quad
  |\w^{(M)}-\w^{(\nu)}|_\nu<\ve \]
for all finite $\nu\in S$.  Thus (\ref{tap31}) becomes
\begin{equation}\label{eq:size11} 
\y\equiv \y^{(M)} \bmod{M}, \quad
\w\equiv \w^{(M)} \bmod{M}
\end{equation}
for an appropriate modulus $M\in\NN$.

Having suitably re-interpreted the weak approximation conditions for
the finite places we turn our attention to the infinite place.  Here
we use a similar re-scaling argument to conclude that if
$Y,W\in\NN$ satisfy 
\[W\equiv 1\bmod{M}\;\;\; \mbox{and}\;\;\;
Y=W^2,\]
then a solution $\y,\w$ of \eqref{yeq2} which satisfies both
\eqref{eq:size11} and the $\nu=\infty$ constraints
\begin{equation}\label{eq:size21}
|\y-Y\y^{(\RR)}|<\ve Y, \quad
|\w-W\w^{(\RR)}|<\ve W, 
\end{equation}
gives rise to a solution $Y^{-1}\y,W^{-1}\w$  
of \eqref{yeq2} which meets the original condition \eqref{tap31}.  Since 
$\y^{(\RR)}$ and $\w^{(\RR)}$ cannot vanish we may choose
$\ve$ sufficiently small that neither of $\y$ or $\w$ can be zero in
\eqref{eq:size21}. 

As in \S~\ref{intro} we replace the vector $\y$ by $\u$ and $\v$ to
produce the variety \eqref{TR2}, whose definition it is convenient to 
repeat here:- 
\begin{equation}  \label{TR3}
\mathcal{J}:\;\str\big(\delta \bN_{K/L}(\u),
(\bN_{K/L}(\v)D_L)^\sigma \big)=2\bN_{K/\QQ}(\w)\neq 0. 
\end{equation}
It is clear that if $\u$ and $\v$ are sufficiently close to
$\u^{(\nu)}:=\y^{(\nu)}$ and
\begin{equation}\label{vone}
\v^{(\nu)}:=(1,0,0,\ldots,0)
\end{equation}
in $\QQ_{\nu}$ then $\y$ will be suitably close to $\y^{(\nu)}$.  We
therefore assume that \eqref{TR3} has local solutions
$\u^{(\nu)},\v^{(\nu)},\w^{(\nu)}$ for all places $\nu$,
with $\v^{(\nu)}$ given by \eqref{vone}, and we aim to find a global solution
such that
\begin{equation}
  \label{tap3}
      |\u-\u^{(\nu)}|_\nu<\ve, \quad
  |\v-\v^{(\nu)}|_\nu<\ve, \quad
  |\w-\w^{(\nu)}|_\nu<\ve, 
\end{equation}
for  $\nu \in S$.  Since $\u^{(\nu)}$ and $\v^{(\nu)}$ are integral at
all finite $\nu\in S$ we can re-interpret the corresponding conditions
as congruences
\[\u\equiv \u^{(M)} \bmod{M},\quad \v\equiv \v^{(M)} \bmod{M}\]
with integer vectors $\u^{(M)}$ and $\v^{(M)}$ for which
\begin{equation}\label{vone2}
\v^{(M)}:=(1,0,0,\ldots,0).
\end{equation}
For technical reasons we will move
$\u^{(\RR)}$ in \eqref{TR3} very slightly, and make a corresponding
adjustment in $\w^{(\RR)}$ to compensate, so as to ensure that
\begin{equation}\label{todo}
\frac{\partial \bN_{K/\QQ}(\u^{(\RR)})}{\partial
  u_i}\not=0,\;\;\;(1\le i\le n).
\end{equation}
For the infinite place we replace the parameter $Y$ by two further
values $U$ and $V$ satisfying $UV=Y$ and impose the conditions
\[|\u-U\u^{(\RR)}|<\ve U, \quad|\v-V\v^{(\RR)}|<\ve V\]
instead of $|\y-Y\y^{(\RR)}|<\ve Y$.

We can now summarise our conclusions in the following result. 

\begin{lemma}\label{lem:suffice}
Suppose we are given local solutions of \eqref{TR3} for every valuation 
$\nu$ of $\QQ$, subject to the condition \eqref{vone}.  Let 
$\ve>0$ also be given.  Then there is a modulus $M\in\NN$ having
$|M|_{\nu}<1$ for all finite $\nu\in S$, and 
a solution $(\u^{(M)},\v^{(M)},\w^{(M)})$ of \eqref{TR3} 
over $\ZZ/M\ZZ$ satisfying \eqref{vone2}, 
having the following property.  Let $V, H_0$ 
be integer parameters with 
$H_0\equiv 
V\equiv 1\bmod{M}$.  Let $H=H_0^2$ and suppose that 
$V\ge H\ge 1$.
If
\[U=HV,\quad W=H^{1/2}V,\]
then any solution 
$(\u, \v, \w)\in\ZZ^{3n}$ of \eqref{TR3} satisfying 
\begin{equation}\label{eq:size1}
\begin{split}
&
\u\equiv \u^{(M)} \bmod{M}, \\
&\v\equiv \v^{(M)} \bmod{M}, \\
&\w\equiv \w^{(M)} \bmod{M},
\end{split}
\end{equation}
and
\begin{equation}\label{eq:size2}
\begin{split}
&|\u-U\u^{(\RR)}|<\ve U, \\
&|\v-V\v^{(\RR)}|<\ve V, \\
&|\w-W\w^{(\RR)}|<\ve W, 
\end{split}
\end{equation}
gives rise to a global solution of 
  \eqref{TR3} satisfying \eqref{tap3}.  Moreover, for every finite
  place $\nu\in S$ there is a solution $(\u,\v,\w)\in\ZZ_\nu^{3n}$ of 
\eqref{TR3} satisfying \eqref{eq:size1}.
\end{lemma}

It might help the reader at this point to say more about the r\^ole of
the parameters $H$ and $V$.   We shall think of $H$ 
as being a small fixed power of $V$.
When we estimate 
error terms in our analysis we cannot afford to lose 
any power $V^{\theta}$ of $V$, unless $\theta$ can be taken
arbitrarily small.  On the other hand there
will be certain points in our argument where
we will lose factors of 
$V^{\eta}$ with arbitrary small
$\eta>0$. This will not matter since we will make a key saving
which is a power of $H$, so that there is a net gain overall.

With this in mind, many of our estimates will involve factors of the
type $V^{O(\eta)}$. 
These involve the standard
convention that there are implicit order constants for each
occurrence of the $O(\cdot)$ notation, which need not be the same on
each occasion.  Since we are taking the degree $n$ of $K$ to be fixed,
we will allow these implicit order constants to depend on
$n$. Recalling that
$H\le V$ we may replace terms involving any combinations of 
$H^{\eta},U^{\eta}$ and $W^{\eta}$ by $V^{O(\eta)}$.
The number $\eta$ will be a sufficiently small
positive constant, which will be fixed throughout the proof.  We could
have chosen to specify its value at the outset of the argument, but we
feel it is more instructive merely to impose the condition that $\eta$ is
sufficiently small, at various points in the proof.

We are now ready to cast our problem in terms of a bilinear form.  If
$R$ is any ring it will be convenient to write $\mathcal{J}(R)$ for
the set of solutions $(\u,\v,\w)$ of the equation
\eqref{TR3} in which each of the vectors has coordinates lying in $R$. 
In the light of Lemma \ref{lem:suffice} we introduce the
counting function
\[
N(H,V):=
\#\{(\u,\v,\w)\in \mathcal{J}(\ZZ): \mbox{\eqref{eq:size1} and
    \eqref{eq:size2} hold}\}.
\]
Henceforth we will allow the constants implied by the notations
$\ll$, $\gg$ and $O(\cdot)$ to depend on 
$$
\u^{(\RR)},\v^{(\RR)},\w^{(\RR)},
\u^{(M)},\v^{(M)}, \w^{(M)}, M, \delta, L, K,\;\;\mbox{and}\;\; \ve,
$$
which are to be regarded as fixed once and for all.

In our work we will restrict the values over which $\v$ runs by 
stipulating that if $\bN_{K/L}(\v)=\bN_1(\v)+\bN_2(\v)\tau$ then
\begin{equation}\label{Econd}
\gcd(\bN_1(\v) ,\bN_2(\v))=1.
\end{equation}
In particular it follows that
\begin{equation}\label{Ebnd}
\gcd(a , b)=O(1)\quad\mbox{for}\quad  
(\bN_{K/L}(\v)D_L)^{\sigma}=a+b\tau.
\end{equation}

We are now ready to specify the sets over which we will sum.  In
the case of the variable $\u$ there is a technical point to be dealt
with in \S~\ref{app2}.  For the time being we give ourselves 
independent linear forms $L_1(\u),\ldots, L_n(\u)$ whose r\^ole will
become clear later.  Let $G$ be a further parameter, tending to
infinity with $V$, which we
assume is in the range $1\leq G \leq H$.
We then define the regions
\begin{equation}
  \label{eq:Abox}
\begin{split}
\cU&:=\left\{\u \in \RR^n: \max_{1\leq i \leq
    n}|L_i(\u)-UL_i(\u^{(\RR)})|<G^{-1}U\right\},\\ 
\cV&:=\left\{\v \in \RR^{n}: 
|\v-V\v^{(\RR)}|<G^{-1}V
\right\},\\
\cW&:=\{\w \in \RR^n: |\w-W\w^{(\RR)}|<G^{-1}W\}.
\end{split}
\end{equation}
In order to interpret our counting function $N(H,V)$ as a bilinear form,
we let 
\begin{equation}
  \label{eq:AL}
\al(x):=
\#
\{
\u\in \cU\cap\ZZ^{n}: 
\u\equiv \u^{(M)} \bmod{M},~ \delta \bN_{K/L}(\u)=x
\}
\end{equation}
and
\begin{equation}
  \label{eq:BE}
\be(y):=\#\left\{
\v\in \cV\cap\ZZ^n:  
\begin{array}{l}
\v\equiv \v^{(M)} \bmod{M},\\
\mbox{\eqref{Econd} holds and $(\bN_{K/L}(\v)D_L)^{\sigma}=y$}
\end{array}
\right\},
\end{equation}
for $x,y \in \Lring$.   Lastly we define the function
\begin{equation}
  \label{eq:LA}
  \la(l):=
\#
\left\{
\w\in \cW\cap\ZZ^{n}:  
\w\equiv \w^{(M)} \bmod{M},~ 2\bN_{K/\QQ}(\w)=l
\right\},
\end{equation}
on $\ZZ$.  Notice that $\u^{(\RR)},\v^{(\RR)}$ and $\w^{(\RR)}$ 
are all non-zero, whence $\u,\v$ and $\w$ will be 
non-zero throughout $\cU,\cV$ and $\cW$, if $G$ is large enough. 
It follows in particular that $\alpha,\beta$ and 
$\lambda$ are supported on non-zero $x,y\in\Lring$ and $l\in\ZZ$. 

We now define the bilinear form
$$
\cl{N}(G, H,V):=
\sum_{x \in \Lring }\sum_{y \in \Lring} \al(x)\be(y)\la\big(\str(x,y)\big).
$$
It is easy to check that
$
N(H,V)\geq \cl{N}(G,H,V)
$ 
whenever $G\gg \varepsilon^{-1}$, where the implied constant is
allowed to depend
on the linear forms $L_1,\ldots,L_n$, a convention that we adhere to
for the remainder of the paper. 
It now suffices to demonstrate that $ \cl{N}(G,H,V)>0$ for large
values of $G$.  Note that we will ultimately take 
$G=\log V.$

Although we have framed  $\cl{N}(G,H,V)$ as a bilinear form,
it is not an upper bound for $\cl{N}(G,H,V)$ that we seek but an asymptotic
formula as $G\rightarrow \infty$. As indicated in the introduction we
will begin by extracting a main term from our expression for
$\cl{N}(G,H,V)$. Instrumental in this will be finding a decomposition
$\alpha(x)=\hat{\alpha}(x)+\alpha_0(x)$, for an appropriate
approximation $\hat{\alpha}(x)$ to $\al(x)$. 
We will then write
\begin{equation}
  \label{eq:me}
\cl{N}(G,H,V)=\cl{M}(G,H,V)+\cl{E}(G,H,V),
\end{equation}
where 
\[\cl{M}(G,H,V):=\sum_{x \in \Lring }\sum_{y \in \Lring} 
\hat{\al}(x)\be(y)\la\big(\str(x,y)\big)\]
is regarded as the main term and 
\[\cl{E}(G,H,V):=\sum_{x \in \Lring }\sum_{y \in \Lring} 
\al_0(x)\be(y)\la\big(\str(x,y)\big)\]
is the error term.  
The handling of  $\cl{E}(G,H,V)$ will be executed in 
\S~\ref{s:bil} and the estimation of 
$\cl{M}(G,H,V)$ will be the subject of \S~\ref{s:M} and \S~\ref{s:ssi}.

Our treatment of $\cl{E}(G,H,V)$ requires bounds for
\[\sum_{\substack{x\in R\\ x\equiv x_0\bmod{h}}}\al_0(x)\]
uniformly for small moduli $h$ and square regions $R$.  
Thus our
approximation $\hat{\al}$ will have to be such that $\al_0$ averages
to zero over all congruence classes to small moduli.  This will be
achieved via a quite general procedure described in the next section.
Our estimate for $\cl{E}(G,H,V)$ will also require bounds for 
\[\sum_{x \in \Lring} |\al_0(x)|^2,\quad
\sum_{y \in \Lring} |\be(y)|^2\quad \mbox{and}\quad
\sum_{l \in \ZZ} |\la(l)|^2,\] 
for which we have the following result.
\begin{lemma}\label{lem:upper-abl}
For any  $\eta>0$ we have
$$
\sum_{x \in \Lring} |\al(x)|^2\ll_\eta U^{n+\eta},\quad
\sum_{y \in \Lring} |\be(y)|^2\ll_\eta V^{n+\eta},
$$
and 
$$
\sum_{l\in \ZZ}|\la(l)|^2 \ll_\eta W^{n+\eta}.
$$
\end{lemma}

Since $|\al_0(x)|^2\le 2|\al(x)|^2+2|\hat{\al}(x)|^2$ we will also require
a bound for 
\[\sum_{x \in \Lring} |\hat{\al}(x)|^2.\]
This will be established in \S~\ref{app2}. 
We conclude this section by proving Lemma~\ref{lem:upper-abl}.
We will discuss the upper bound for the case of
$\be$, the remaining estimates being dealt with similarly. 
Let 
$h(\v)=(\bN_{K/L}(\v)D_L)^{\sigma}$. 
We shall show that if $\v'$ is given, then there are
$O_{\eta}(V^{\eta})$ choices of $\v$ for which
$h(\v)=h(\v')$. We set $\varpi=h(\v')$, so that
$\varpi\in\Lring\setminus \{0\}$.  
If
$\rho=v_1\omega_1+\cdots+v_n\omega_n\in\Kring$, then $\rho\mid \varpi$ in
$\Kring$, and the conjugates of $\rho$ are all $O(V)$ in absolute 
magnitude.  The number of
admissible elements $\rho$ is therefore 
\[
\ll_{\xi}V^{\xi}|N_{L/\QQ}(\varpi)|^{\xi}
\ll_{\xi}V^{\xi}(V^n)^{\xi}= V^{(1+n)\xi}
\]
for any $\xi>0$.  If we now take $\xi=\eta/(1+n)$ it follows that the number of 
$\v$ corresponding to a given
$\v'$ is $O_{\eta}(V^{\eta})$, as required.

The fact that a given value $\varpi=h(\v)$ 
is attained
$O(V^{\eta})$ times will be used at various points in the rest of the
paper without further comment.  Similarly,  we shall  use estimates
$O(V^{\eta})$ for
the number of representations of $\varpi$ as $\bN_{K/L}(\u)$ with
$\u\in\cU$, or of $l$ as $2\bN_{K/\QQ}(\w)$ with $\w\in\cW$.

\section{A general approximation principle} \label{s:approx}

Our goal now is to split
$\alpha(x)$ into two parts $\alpha(x)=\hat{\alpha}(x)+\alpha_0(x)$,
where $\hat{\alpha}(x)$ is a
sufficiently simple function that we can compute $\cl{M}(G,H,V)$
directly. 
Moreover we will want $\alpha_0(x)$ to produce sufficient cancellation on
average, so that a bilinear form estimation of $\cl{E}(G,H,V)$ can be achieved.

We start by describing a general procedure for producing an
approximation of the type $\hat{\alpha}(x)$.  In our context $x$ runs
over the ring $\Lring$, but we will begin by presenting the method as
it applies to sequences indexed by $\Z$, since we hope this will prove
to be of independent interest.
The underlying ideas are perhaps not new. In particular there are 
certain similar features in recent work of Br\"udern \cite{brudern}.  
However we have not found anything in the literature which exactly 
meets our needs. 

We begin by supposing that we are given a  sequence $k(1),\ldots,k(N)$
of complex 
numbers.  We aim to show how to approximate $k(n)$
locally by a function $\hat{k}(n)$.  By this we mean that for any congruence
class $a\bmod{q}$ the sum
\[S(a,q):=\twosum{n\le N}{n\equiv a\bmod{q}}k(n)\]
will be approximated by
\[\hat{S}(a,q):=\twosum{n\le N}{n\equiv a\bmod{q}}\hat{k}(n),\]
at least for small values of $q$.  Naturally we can do this by merely
setting $\hat{k}(n)=k(n)$, but we seek a function $\hat{k}(n)$ 
which is defined in
terms of the density of the sequence $k(n)$ in congruence classes to
small moduli.

As an example of what we have in mind, consider the sequence
$k(n)=\Lambda(n)$, where $\Lambda(n)$ is the von Mangoldt function.  If
we choose any constant $A\ge 3$ and set
\[\Lambda_Q(n):=
\sum_{q\le Q}\frac{\mu(q)}{\phi(q)}\;\twosum{a\bmod{q}}{\gcd(a,q)=1}e_q(an),\]
with $Q=(\log x)^A$, then
\begin{equation}
  \label{eq:gold}
\twosum{n\le x}{n\equiv b\bmod{h}}\Lambda(n)=
\twosum{n\le x}{n\equiv b\bmod{h}}\Lambda_Q(n)+O_A(x(\log x)^{-A/2})
\end{equation}
uniformly for $h\le Q$, for all residue classes $b\bmod{h}$.  This
result follows from Heath-Brown \cite[Lemma~1]{RMI}.

We first introduce our fundamental hypotheses.  We assume that we have
an arithmetic function $\rho(a,q)$ and a ``smooth'' function $\omega(n)$, 
together with a bound $E\ge 1$ 
such that
\begin{equation}\label{E}
\left|S(a,q)-\rho(a,q)S\right|\le E,
\end{equation}
with
\[S=\sum_{n\le N}\omega(n),\]
for all residue classes $a$ to moduli $q\le Q$.  We observe that
\[\twosum{a\bmod{rs}}{a\equiv b\bmod{r}}S(a,rs)=S(b,r)\]
for all $b,r,s$, and we therefore impose the natural condition that
\begin{equation}\label{comb}
\twosum{a\bmod{rs}}{a\equiv b\bmod{r}}\rho(a,rs)=\rho(b,r)
\end{equation}
for all $b,r,s$.  Since we can always re-scale the functions $\rho$ and
$\omega$ in \eqref{E} there is no loss in generality in assuming that
\begin{equation}\label{one}
\rho(0,1)=1.
\end{equation}
Although it is not necessary in general, it will
prove convenient to assume that
\[\rho(a,q)\in\R,\;\;\;\rho(a,q)\ge 0\]
for all pairs $a,q$.

We will also require a smoothness condition on the
function $\omega(n)$, which we formulate as the bound
\begin{equation}\label{W}
\left|\twosum{n\le N}{n\equiv a\bmod{q}}\omega(n)-q^{-1}S\right|\le W,
\end{equation}
for all residue classes $a$ to moduli $q\le Q^2$.
\bigskip

Our choice for $\hat{k}(n)$ will be motivated by the treatment of
major arcs in the circle method.
If we consider the exponential sum
\[\Sigma(\alpha):=\sum_{n\le N}k(n)e(\alpha n),\]
then when $\alpha$ is close to $a/q$ one would use the major-arc
approximation
\[\left\{\sum_{c\bmod{q}}\rho(c,q)e_q(ac)\right\}
\left\{\sum_{n\le N}\omega(n)e((\alpha-a/q)n)\right\}.\]
When $\alpha$ is not close to $a/q$ the above expression tends to
be small.  Hence it is reasonable to approximate $\Sigma(\alpha)$ by
\[\sum_{q\le Q}\;\twosum{a\bmod{q}}{\gcd(a,q)=1}
\left\{\sum_{c\bmod{q}}\rho(c,q)e_q(ac)\right\}
\left\{\sum_{n\le N}\omega(n)e((\alpha-a/q)n)\right\}\]
for all real $\alpha$.  Picking out the coefficient of $e(\alpha n)$
in the above expression we are therefore led to suggest the choice
\begin{equation}\label{hatdef}
\hat{k}(n):=\omega(n)\sum_{q\le Q}\;\twosum{a\bmod{q}}{\gcd(a,q)=1}
e_q(-an)\sum_{c\bmod{q}}\rho(c,q)e_q(ac).
\end{equation}

We proceed to investigate the sum $\hat{S}(b,h)$ with $h\le Q$.  We
have
\[\hat{S}(b,h)=\twosum{n\le N}{n\equiv b\bmod{h}}\omega(n)
\sum_{q\le Q}\;\sum_{c\bmod{q}}\rho(c,q)\twosum{a\bmod{q}}{\gcd(a,q)=1}
e_q(a(c-n)).\]
The final sum over $a\bmod{q}$ is a Ramanujan sum, for which the standard
evaluation as
\[
\sum_{d\mid \gcd(q,c-n)} d\mu(q/d)
\]
produces
\[\hat{S}(b,h)=\sum_{q\le Q}\;\sum_{c\bmod{q}}\rho(c,q)\;\sum_{d\mid
  q}
d\mu(q/d)
\sum_{\substack{n\le N\\ n\equiv b\bmod{h}\\ n\equiv c\bmod{d}}}\omega(n).\]
The simultaneous congruences $n\equiv b\bmod{h}$ and $n\equiv c\bmod{d}$
are only soluble if $\gcd(d,h)$ divides $ c-b$,  
in which case there is a unique solution
$e\bmod{[d,h]}$ say.  Thus if $\gcd(d,h)\mid c-b$ our hypothesis \eqref{W}
shows that
\[\hat{S}(b,h)=\cl{M}+\cl{E},\]
with
\[\cl{M}=S\sum_{q\le Q}\;\sum_{c\bmod{q}}\rho(c,q)\twosum{d\mid
  q}{\gcd(d,h)\mid c-b}
\frac{d\mu(q/d)}{[d,h]}\]
and
\[|\cl{E}|\le W\sum_{q\le Q}\;\sum_{c\bmod{q}}\rho(c,q)\sum_{d\mid
  q}d|\mu(q/d)|.\] 
In view of \eqref{comb} and \eqref{one} we have
\[\sum_{c\bmod{q}}\rho(c,q)=\rho(0,1)=1,\]
whence the crude bound
\[\sum_{d\mid q}d|\mu(q/d)|\le q^2\]
yields $|\cl{E}|\le WQ^3.$

Turning to the main term $\cl{M}$ we observe in general that
\begin{align*}
\twosum{d\mid q}{\gcd(d,h)\mid k}\frac{d\mu(q/d)}{[d,h]}
&=
h^{-1}\twosum{d\mid q}{\gcd(d,h)\mid k}\mu(q/d)(d,h)\\
&=h^{-1}\sum_{d\mid q}\mu(q/d)\sum_{e\mid d,h,k}e\sum_{f\mid d/e,h/e}\mu(f)\\
&=h^{-1}\sum_{e\mid q,h,k}e\sum_{f\mid q/e,h/e}\mu(f)
\twosum{d\mid q}{ef\mid d}\mu(q/d)\\
&=h^{-1}\sum_{e\mid q,h,k}e\sum_{f\mid q/e,h/e}\mu(f)
\sum_{g\mid q/(ef)}\mu\left(\frac{q/(ef)}{g}\right).
\end{align*}
The final sum vanishes unless $ef=q$, in which case we must have
$q\mid h$.  It then follows that
\begin{align*}
\twosum{d\mid q}{\gcd(d,h)\mid k}\frac{d\mu(q/d)}{[d,h]}
&=h^{-1}\sum_{e\mid q,k}e\mu(q/e)\\
&=h^{-1}\twosum{a\bmod{q}}{\gcd(a,q)=1}e_q(ak).
\end{align*}
Inserting this result into our formula for the main term $\cl{M}$, we see
that
\[\cl{M}=Sh^{-1}
\sum_{q\mid h}\sum_{c\bmod{q}}\rho(c,q)\twosum{a\bmod{q}}{
\gcd(a,q)=1}e_q(a(c-b)).\] 
Using \eqref{comb} we have
\[\rho(c,q)=\twosum{d\bmod{h}}{d\equiv c\bmod{q}}\rho(d,h),\]
whence
\[\cl{M}=Sh^{-1}\sum_{d\bmod{h}}\rho(d,h)\sum_{q\mid h}
\twosum{a\bmod{q}}{\gcd(a,q)=1}e_q(a(d-b)).\]
As $q$ runs over divisors of $h$, and $a$ runs over residue classes
coprime to $q$, the fractions $a/q$ run over the entire set 
\[\left\{\frac{n}{h}:\, 0\le n<h\right\}.\]
We therefore deduce that
\[\cl{M}=Sh^{-1}\sum_{d\bmod{h}}\rho(d,h)\sum_{n\bmod{h}}e_h(n(d-b))
=S\rho(b,h),\]
since the summation over $n$ produces the value $h$ when $d\equiv
b\bmod{h}$, and the value $0$ otherwise.

In view of \eqref{E} we may now summarise our results as follows.
\begin{lemma}\label{Mainapp}
With the above assumptions we have
\[|S(b,h)-\hat{S}(b,h)|\le WQ^3+E\]
for all $h\le Q$ and all residue classes $b$ modulo $h$.
\end{lemma}

Thus $\hat{k}(n)$ approximates $k(n)$ well, in congruence classes to small
moduli. It may be instructive to consider the effect of this procedure
on the sequence $k(n)=\Lambda(n)$ that we discussed earlier. Taking
$$
\rho(a,q)=\begin{cases}
1/\phi(q), & \mbox{if $\gcd(a,q)=1$,}\\
0, &\mbox{otherwise,}
\end{cases}
$$
we readily deduce from the Siegel--Walfisz theorem that for any $A\geq
1$ there exists a constant $C_A>0$ such that 
$$
\sum_{\substack{n\leq N\\ n\equiv a\mod{q}}}\Lambda(n)=\rho(a,q)N
+O\left(N\exp(-C_A (\log N)^{1/2})\right),
$$
uniformly for $q\leq (\log N)^A$. On taking $\omega(n)=1$
we see that $E\ll N\exp(-C_A (\log N)^{1/2})$ 
is admissible in \eqref{E}. Since
$W\ll 1$ in \eqref{W}, the approximation in \eqref{eq:gold} is a trivial
consequence of Lemma \ref{Mainapp} with 
\begin{align*}
\Lambda_Q(n)
&=
\omega(n)\sum_{q\le Q}\;\twosum{a\bmod{q}}{\gcd(a,q)=1}
e_q(-an)\sum_{\substack{b\bmod{q}\\
    \gcd(b,q)=1}}\frac{e_q(ab)}{\phi(q)}\\
&=
\sum_{q\le
  Q}\frac{\mu(q)}{\phi(q)}\;\twosum{a\bmod{q}}{\gcd(a,q)=1}e_q(an),
\end{align*}
for $q\leq (\log N)^A$.

\bigskip

In applications it may be important to know about the size of $\hat{k}(n)$,
and we therefore investigate the mean square
\[\Sigma:=\sum_{n\le N}|\hat{k}(n)|^2.\]
If we write
\[c_{a,q}:=\sum_{b\bmod{q}}\rho(b,q)e_q(ab)\]
then
\[\hat{k}(n)=\omega(n)\sum_{q\le
  Q}\twosum{a\bmod{q}}{\gcd(a,q)=1}c_{a,q}e_q(-an).\] 
Thus, if we assume that
\[
|\omega(n)|\le \omega_0
\]
for all $n\in\NN$, then the dual large sieve produces
\begin{align*}
\Sigma&\le\omega_0^2\sum_{n\le N}\left|\sum_{q\le Q}
\twosum{a\bmod{q}}{\gcd(a,q)=1}c_{a,q}e_q(-an)\right|^2\\
&\le\omega_0^2(N+Q^2)\sum_{q\le Q}\twosum{a\bmod{q}}{\gcd(a,q)=1}|c_{a,q}|^2.
\end{align*}
In view of \eqref{E} we have
\[|Sc_{a,q}-T(a,q)|\le qE,\]
with
\[T(a,q):=\sum_{b\bmod{q}}S(b,q)e_q(ab)=\sum_{n\le N}k(n)e_q(an).\]
It follows that $|Sc_{a,q}|^2\le 2|T(a,q)|^2+2q^2E^2$, and hence that
\[|S|^2\Sigma\le 2\omega_0^2(N+Q^2)\left(E^2Q^4+\sum_{q\le Q}
\twosum{a\bmod{q}}{\gcd(a,q)=1}|T(a,q)|^2\right).\]
We now apply the standard large sieve to deduce that
\[|S|^2\Sigma\le 2\omega_0^2(N+Q^2)
\left(E^2Q^4+(N+Q^2)\sum_{n\le N}|k(n)|^2\right).\]
We express this result formally in the following lemma.

\begin{lemma}\label{L2}
We have
\[\sum_{n\le N}|\hat{k}(n)|^2\le 2|S|^{-2}\omega_0^2(N+Q^2)
\left(E^2Q^4+(N+Q^2)\sum_{n\le N}|k(n)|^2\right).\]
Moreover
\[\sum_{n\le N}|\hat{k}(n)|^2\ll
\left(\frac{\omega_0N}{|S|}\right)^2\left(N+\sum_{n\le N}|k(n)|^2\right),\]
providing that $EQ^2\le N$.
\end{lemma}

Thus under suitable circumstances the $L^2$-norm of $\hat{k}(n)$ will
have  order of magnitude bounded by the $L^2$-norm of $k(n)$.

\bigskip

Having described the situation for sequences $k(n)$ indexed by $\Z$ we
return to our original problem, in which we have a sequence
$\alpha(x)$ with $x$ running over $\Lring$. We will describe the
situation as it applies to a quite general function $\alpha$, and only
later, in \S~\ref{app2}, 
restrict to the function \eqref{eq:AL}.  We shall consider sums in
which $x$ runs over a region $R$, say, and lies in a congruence class
$y\bmod{q}$, where $y\in\Lring$ and $q\in\NN$.  We therefore set
\begin{equation}\label{Ssetup}
S(y,q):=\twosum{x\in R}{x\equiv y\bmod{q}}\alpha(x)
\end{equation}
and
\[\hat{S}(y,q):=\twosum{x\in R}{x\equiv y\bmod{q}}\hat{\alpha}(x),\]
and assume that we have functions $\rho(y,q)$ and $\omega(x)$ such
that $\rho(y,q)$ is non-negative and
\begin{equation}\label{E2}
\left|S(y,q)-\rho(y,q)S\right|\le E
\end{equation}
for $q\le Q$, with $E\ge 1$ and
\[S:=\sum_{x\in R}\omega(x).\]
As before we will assume that
\begin{equation}\label{comb2}
\twosum{y\bmod{rs}}{y\equiv z\bmod{r}}\rho(y,rs)=\rho(z,r)
\end{equation}
and
\begin{equation}\label{one2}
\rho(0,1)=1.
\end{equation}
Finally, we shall suppose that
\begin{equation}\label{W2}
\left|\twosum{x\in R}{x\equiv y\bmod{q}}\omega(x)-q^{-2}S\right|
\le W,
\end{equation}
for all residue classes $y$ to moduli $q\le Q^2$.

In order to define $\hat{\alpha}(x)$ we will need an analogue of the
coprimality 
condition $\gcd(a,q)=1$ which occurs in \eqref{hatdef}.  It turns out that
the correct generalisation is to require that there is no non-trivial
common divisor of $y$ and $q$ in $\NN$.  Of course this is not the same
as requiring $\gcd(y,q)=1$ in $\Lring$.  We write 
\[\Osum_{y\bmod{q}}\]
to denote a sum in which
$y\in\Lring$ runs over residue classes modulo $q$ such that $y$ and
$q$ have no non-trivial common divisor in $\NN$.  

We also need an
analogue for the exponential function $e_q(a)$.  Recall that 
$\{1,\tau\}$ is a
$\Z$-basis for $\Lring$, and hence also a $\Q$-basis for $L$.  If
$x=a+b\tau\in L$, with $a,b\in\Q$, and if $q\in\NN$, we define
\[e^{(L)}(x):=e(b)\;\;\;\mbox{and}\;\;\;e_q^{(L)}(x):=e_q(b).\]
These exponential functions have the property that
if $y\in\Lring$ then
\[
\sum_{x\bmod{q}}e_q^{(L)}(xy)=
\begin{cases}
q^2, & \mbox{if $q\mid y$,}\\
0, & \mbox{if $q\nmid y$.}
\end{cases}
\]
Thus, for the analogue of the Ramanujan sum we have
\begin{equation}
  \label{eq:gen-ram}
\Osum_{x\bmod{q}}e_q^{(L)}(xy)=\sum_{d\mid \gcd(q,y)}d^2\mu(q/d).  
\end{equation}

We are now ready to specify our approximation to $\alpha(x)$.  We set
\begin{equation}\label{hataldef}
\hat{\alpha}(x):=\omega(x)\sum_{q\le Q}\Osum_{y\bmod{q}}
e_q^{(L)}(-yx)\sum_{z\bmod{q}}\rho(z,q)e_q^{(L)}(yz).
\end{equation}
An argument precisely analogous to that used for Lemma \ref{Mainapp}
then produces the following result.
\begin{lemma}\label{Mainapp2}
With the above assumptions we have
\[|S(z,h)-\hat{S}(z,h)|\le WQ^4+E\]
for all $h\le Q$ and all residue classes $z$ modulo $h$.
\end{lemma}

To produce an analogue of Lemma \ref{L2} we shall assume that
\[R=\{a+b\tau\in \Lring:\,|a-a_0|,|b-b_0|< N\}\]
for certain $a_0,b_0$.
Then we will have a large sieve inequality for $\Lring$, taking the
form
\begin{equation}\label{ls2}
\sum_{q\le Q}\Osum_{y\bmod{q}}
\left|\sum_{x\in R}c_xe_q^{(L)}(xy)\right|^2\le
(\sqrt{2N}+Q)^4\sum_{x\in R}|c_x|^2.
\end{equation}
This follows from the two dimensional large sieve 
of Huxley \cite[Theorem~1]{hux}.    Under the condition that
$Q^2\le N$, which we now impose, we will then have $(\sqrt{2N}+Q)^4\ll N^2$.
We proceed to argue as before
to deduce from the dual of the above estimate that
\[\sum_{x\in R}|\hat{\alpha}(x)|^2\ll\omega_0^2N^2
\sum_{q\le Q}\Osum_{y\bmod{q}}|c_{y,q}|^2,\]
with
\[c_{y,q}=\sum_{z\bmod{q}}\rho(z,q)e_q^{(L)}(yz).\]
This time the estimate \eqref{E2} yields
\[|Sc_{y,q}-T(y,q)|\le q^2E,\]
with
\[T(y,q)=\sum_{x\in R}\alpha(x)e_q^{(L)}(xy).\]
Continuing as before we then obtain
\[|S|^2\sum_{x\in R}|\hat{\alpha}(x)|^2\ll\omega_0^2N^2
\left(E^2Q^7+N^2\sum_{x\in R}|\alpha(x)|^2\right).\]
This gives us the following conclusion.
\begin{lemma}\label{L22}
Suppose that
\[R=\{a+b\tau\in \Lring:\,|a-a_0|,|b-b_0|< N\}\] 
for certain $a_0,b_0$, and assume that there is a constant $\omega_0$ such that
$|\omega(x)|\le \omega_0$ for all $x\in R$.  Then
\[\sum_{x\in R}|\hat{\alpha}(x)|^2\ll
\left(\frac{\omega_0N^2}{|S|}\right)^2
\left(N^2+\sum_{x\in R}|\alpha(x)|^2\right),\]
providing that $Q^2\le N$ and $E^2Q^7\le N^4$.
\end{lemma}

\section{The functions $\hat{\al}(x)$ and $\al_0(x)$}\label{app2}

In this section we will apply Lemmas \ref{Mainapp2} and \ref{L22} to
the function \eqref{eq:AL}.  It will be important for us to produce
functions $\hat{\alpha}$ and $\omega$ which depend only on the set
$\cU$, and not on the set $R$, since we require results which are
uniform in $R$.

We will write $x\in\Lring$ in the form 
$x=x_1+ x_2\tau$ with $x_1,x_2\in\ZZ$.  We shall also
write $\delta \bN_{K/L}(\u)$ as
$\bN_1(\u)+\bN_2(\u)\tau$. The reader should observe that this
notation does not coincide with that used temporarily in
\eqref{nsplit}, nor that used in our discussion of \eqref{Econd}.
We let $R$ be a square in the $(x_1,x_2)$-plane, with sides
parallel to the $x_1$ and $x_2$ axes, and side-length $\rho
\ll U^{n/2}$. This corresponds, by abuse of notation, to the square $R$
used in defining the sums $S(y,q)$ in \eqref{Ssetup}. Extending this
abuse of notation we shall allow $\delta \bN_{K/L}(\u)$ to denote the 
ordered pair $(\bN_1(\u),\bN_2(\u))$.

We will specify the function $\omega(x)=\omega(x_1,x_2)$ and verify 
its properties at the
end of this section.  For the time being we content ourselves with
describing the key features as follows.

\begin{lemma}\label{omegaprops}
There is a continuously differentiable
function $\omega:\RR^2\rightarrow [0,\infty)$, depending on
$\cU$, for which $\omega(\x+\hh)-\omega(\x)\ll U^{-n/2}|\hh|$ for all
$\x,\hh\in\RR^2$, and such that
\begin{equation}\label{omint}
\int_R\omega(x_1,x_2)\d x_1\d x_2=M^n\meas \{\u\in\cU:
\delta\bN_{K/L}(\u)\in R\},  
\end{equation}
for every square $R$ as above. 
Furthermore $\omega$ is supported on a disc of radius $O(U^{n/2})$ and
satisfies $\omega(\x)\ll 1$ throughout this disc.  Moreover there is a
disc of radius $\gg G^{-1}U^{n/2}$ around the point 
$U^{n/2}\delta \bN_{K/L}(\u^{(\RR)})$ 
on which we have $\omega(\x)\gg G^{2-n}$.
\end{lemma}

Although it is a real measure which occurs on the right-hand-side of
\eqref{omint} what occurs naturally for us is the corresponding
cardinality
\[S_0:=
\#\{\u\in \cU\cap\ZZ^n:\delta\bN_{K/L}(\u)\in R\}.\]
We proceed to establish a relation between the two.  For each 
$\u\in \cU\cap\ZZ^n$ we let 
\[S(\u):=\{\y\in\RR^n: u_i\le y_i<u_i+1,~ (1\le i\le n)\}\] 
and
\[\cU^{(-)}:=\bigcup\{S(\u):\u\in\ZZ^n, ~S(\u)\subseteq\cU\}.\] 
Thus $\cU^{(-)}\subseteq\cU$, and the number of integer vectors in $\cU$
but not $\cU^{(-)}$ is $O(U^{n-1})$.  In particular,
\begin{equation}\label{AA}
S_0=
\#\{\u\in \cU^{(-)}\cap\ZZ^n: \delta\bN_{K/L}(\u)\in
R\}+O(U^{n-1}).
\end{equation}
Now if $\u\in \cU^{(-)}\cap\ZZ^n$ and $\y\in S(\u)$ then
$\delta\bN_{K/L}(\y)=\delta\bN_{K/L}(\u)+O(U^{n/2-1})$.  Thus there are squares
$R_1$ and $R_2$ of sides $\rho_1$ and $\rho_2$ respectively, such that
\[|\rho_1-\rho|\ll U^{n/2-1}\;\;\;\mbox{and}\;\;\;
|\rho_2-\rho|\ll U^{n/2-1},\]
and for which
\[\delta
\bN_{K/L}(\y)\in R_1\;\Rightarrow \;\delta \bN_{K/L}(\u)\in R
\;\Rightarrow \;\delta\bN_{K/L}(\y)\in R_2\]
whenever $\y\in S(\u)$.  We deduce that
\begin{align*}
\#\{\u\in \cU^{(-)}\cap\ZZ^n:\, &\delta\bN_{K/L}(\u)\in R\}\\
&\ge
\meas\{\y\in\cU^{(-)}: \delta\bN_{K/L}(\y)\in R_1\}
\end{align*}
and
\begin{align*}
\#\{\u\in \cU^{(-)}\cap\ZZ^n:\, &\delta\bN_{K/L}(\u)\in R\}\\
&\le
\meas\{\y\in\cU^{(-)}: \delta\bN_{K/L}(\y)\in R_2\}.
\end{align*}
However $\cU^{(-)}$ is contained in $\cU$ and differs from it 
by a set of measure $O(U^{n-1})$, whence
\begin{align*}
\meas\{\y\in\cU^{(-)}:\, &\delta\bN_{K/L}(\y)\in R_1\}\\
&\ge
\meas\{\y\in\cU: \delta\bN_{K/L}(\y)\in R_1\}+O(U^{n-1})
\end{align*}
and
\[
\meas\{\y\in\cU^{(-)}: \delta\bN_{K/L}(\y)\in R_2\}
\le \meas \{\y\in\cU: \delta\bN_{K/L}(\y)\in R_2\}.
\]
According to \eqref{omint} we may then deduce that
\[\#\{\u\in \cU^{(-)}\cap\ZZ^n:\delta\bN_{K/L}(\u)\in R\} 
\ge M^{-n}\int_{R_1}\omega(x_1,x_2)\d x_1\d x_2+O(U^{n-1})\] 
and
\[\#\{\u\in \cU^{(-)}\cap\ZZ^n:\delta\bN_{K/L}(\u)\in R\}
\le M^{-n}\int_{R_2}\omega(x_1,x_2)\d x_1\d x_2. 
\]
However the squares $R_1$ and $R_2$ each differ from $R$ by a
set of measure $O(U^{n-1})$, and furthermore $\omega\ll 1$.  Thus
\[\#\{\u\in \cU^{(-)}\cap\ZZ^n:\delta\bN_{K/L}(\u)\in R\}=
M^{-n} 
\int_{R}\omega(x_1,x_2)\d x_1\d x_2+O(U^{n-1}).\]
Finally, if we compare this with \eqref{AA}, we deduce that
\begin{equation}\label{BB}
S_0=M^{-n} 
\int_{R}\omega(x_1,x_2)\d x_1\d x_2+O(U^{n-1}).
\end{equation}

We next consider
\[\twosum{\x\in R}{\x\equiv \y\bmod{q}}\omega(\x).\]
To each point $\x$ counted in the above sum we associate the square
$R(\x)=[x_1,x_1+q)\times [x_2,x_2+q)$.  If $\t\in R(\x)$ then
$\omega(\t)-\omega(\x)\ll qU^{-n/2}$ by Lemma~\ref{omegaprops}, whence
\[\int_{R(\x)}\omega(\t)\d t_1\d t_2=q^2\omega(\x)+O(q^3U^{-n/2}).\]
We sum this for points $\x\in R$ with $\x\equiv \y\bmod{q}$.  Since
there are $O(q^{-2}U^n)$ such points if $q\le U^{n/2}$ we deduce that
\[\twosum{\x\in R}{\x\equiv \y\bmod{q}}\omega(\x)=
q^{-2}\twosum{\x\in R}{\x\equiv \y\bmod{q}}
\int_{R(\x)}\omega(\t)\d t_1\d t_2+O(qU^{n/2}).\]
The union of the squares $R(\x)$  will be a square
$R'$ say, whose sides are within a distance $q$ of the sides of $R$.
Thus $R$ and $R'$ differ by a set of measure $O(qU^{n/2})$, if $q\le
U^{n/2}$. Since $\omega(\t)\ll 1$ for all $\t$ this produces
\begin{align*}
\twosum{\x\in R}{\x\equiv \y\bmod{q}}\omega(\x)&=
q^{-2}\int_{R'}\omega(\t)\d t_1\d t_2+O(qU^{n/2})\\
&= q^{-2}\int_{R}\omega(\t)\d t_1\d t_2+O(qU^{n/2}).
\end{align*}
Comparing this with \eqref{BB} we therefore obtain
\[M^{-n}
\twosum{\x\in R}{\x\equiv\y\bmod{q}}\omega(\x)=q^{-2}S_0+O(U^{n-1}),\]
providing that $n\ge 4$ and $q\le U$.
In particular, when $q=1$ we obtain
\begin{equation}\label{S0S}
S=M^{n}S_0+O(U^{n-1}),
\end{equation}
so that \eqref{W2} holds with $W\ll U^{n-1}$ and any $Q\le U^{1/2}$.   

We next consider the condition \eqref{E2}.  In view of the definition
\eqref{eq:AL} we see that
\[S(y,q)=\sum_\w\#\{\u\in \cU\cap\ZZ^n:\u\equiv\w\bmod{[M,q]},\,
\delta\bN_{K/L}(\u)\in R\},\]
where $\w$ runs over vectors modulo $[M,q]$ for which
$\w\equiv\u^{(M)}\bmod{M}$ and $\delta\bN_{K/L}(\w)\equiv y\bmod{q}$.  For a
general modulus $r\le U$ we proceed to compare the sizes of the sets
\[T(r,\x):=\{\u\in \cU\cap\ZZ^n:\u\equiv\x\bmod{r},\,
\delta\bN_{K/L}(\u)\in R\}\]
for the values $\x=\w$ and $\mathbf{0}$. For $\u=(u_1,\ldots,u_n)$ let
\[\u^*=(r[u_1/r],\ldots,r[u_n/r]).\]
If $\u\in T(r,\w)$ then $\u^*$ will
belong to $T(r,\mathbf{0})$ unless either $\u$ is within a distance
$O(r)$ of the boundary of $\cU$, or $\delta\bN_{K/L}(\u)$ is within a distance
$O(rU^{n/2-1})$  of the boundary of $R$.  A pair $(N_1,N_2)$ arises at
most $O_{\eta}(U^{\eta})$ times as a value of $\delta\bN_{K/\QQ}(\u)$ with
$\u\in\cU\cap\ZZ^n$, and hence it follows that 
\begin{equation}\label{new21}
\#T(r,\mathbf{0})=\# T(r,\mathbf{w})+O_{\eta}(rU^{n-1+\eta})
\end{equation}
for any $\eta>0$.  If we now sum for all $\w\bmod{r}$ we deduce that
\[r^n\#T(r,\mathbf{0})=\#\{\u\in \cU\cap\ZZ^n:\delta\bN_{K/L}(\u)\in R\}
+O_{\eta}(r^{n+1}U^{n-1+\eta})\]
and hence that
\[r^n\#T(r,\mathbf{0})=S_0+O_{\eta}(r^{n+1}U^{n-1+\eta}).\]
Substituting this back into \eqref{new21} we deduce that 
\[\# T(r,\w)=r^{-n}S_0+O_{\eta}(rU^{n-1+\eta}).\]
We therefore see that
\[S(y,q)=M^{-n}\rho(y,q)S_0+O_{\eta}(q^{n+1}U^{n-1+\eta}),\]
with
\begin{equation}
  \label{eq:vetch}
  \rho(y,q):=
\frac{M^n}{[M,q]^{n}}\#\left\{\mathbf{s}\bmod{[M,q]}:~
\begin{array}{l}
\mathbf{s}\equiv\u^{(M)}\bmod{M},\\
\delta\bN_{K/L}(\mathbf{s})\equiv y\bmod{q}
\end{array}\right\}.
\end{equation}
The conditions \eqref{comb2} and \eqref{one2} are now readily checked,
and we see that \eqref{E2} follows from \eqref{S0S} with
$E\ll_{\eta}Q^{n+1}U^{n-1+\eta}$. 

We are now ready to apply Lemma \ref{Mainapp2}, which produces the
following result.

\begin{lemma}\label{Mainapp3}
Let $R$ be a square with side $\rho\ge 1$ satisfying $\rho\ll
U^{n/2}$.  Then if $\eta$ is
any positive constant we have
\[\twosum{x\in R}{x\equiv y\bmod{q}}\alpha_0(x)
\ll_{\eta}Q^{n+1}U^{n-1+\eta},\]
for all $q\le Q\le U^{1/2}$ and all $y\bmod{q}$. 
\end{lemma}

In our application the square $R$ will vary and so it is 
crucial that Lemma \ref{Mainapp3}  is uniform in squares of side length 
$O(U^{n/2})$. 
As to Lemma \ref{L22}, we will only need a result for the $L^{2}$-norm
taken over all $x$.  Hence we choose $R$ to be a square centred at the
origin and with side
$\rho=cU^{n/2}$, in which the constant $c$ is taken sufficiently 
large that $x\in R$ whenever $\hat{\alpha}(x)\not=0$.
We may clearly take
$\omega_0\ll 1$ by Lemma \ref{omegaprops} and 
\[\sum_{x\in R}|\alpha(x)|^2\ll_{\eta}U^{n+\eta}\]
by Lemma \ref{lem:upper-abl}.
Our remaining task is thus to estimate $S$ from below. However our
choice of $R$ ensures that $S_0=\#(\cU\cap\ZZ^n)\gg U^nG^{-n}$, whence 
\eqref{S0S} yields $S\gg U^nG^{-n}$, assuming that $G\le U^{1/(n+1)}$, 
say.  Consequently we deduce from Lemma
\ref{L22} the following bound.
\begin{lemma}\label{L23}
For any constant $\eta>0$ we have
\[
\sum_{x\in\Lring}|\hat{\alpha}(x)|^2\ll_{\eta} U^{n+\eta}G^{O(1)}
\]
providing that $Q^{n+5}\le U^{1-\eta}$ and $G\leq U^{1/(n+1)}$.
Moreover we then have
\[
\sum_{x\in\Lring}|\alpha_0(x)|^2\ll_{\eta} U^{n+\eta}G^{O(1)}.
\]
\end{lemma}
The final estimate follows immediately from Lemma \ref{lem:upper-abl},
since
\[  \sum_{x\in\Lring} |\al_0(x)|^2\ll
  \sum_{x\in\Lring} |\al(x)|^2+ \sum_{x\in\Lring}
  |\hat{\al}(x)|^2.\]
\bigskip

We end this section by proving Lemma \ref{omegaprops}. 
We first show 
that the map from $\CC^n$ to $\CC^2$ given by
$\u\mapsto(\bN_1(\u),\bN_2(\u))$ is non-singular at any point for
which $\bN_{K/\QQ}(\u)\not=0$, by which we mean
that $\nabla\bN_1(\u)$ and $\nabla\bN_2(\u)$ are only proportional if
$\bN_{K/\QQ}(\u)=0$. Clearly it suffices to do the same for the map 
\[
\u\mapsto(\bN_1(\u)+\tau\bN_2(\u),\bN_1(\u)+\tau^{\sigma}\bN_2(\u)).
\]
However $\bN_1(\u)+\tau\bN_2(\u)$ (respectively,
$\bN_1(\u)+\tau^\sigma\bN_2(\u)$) is 
a product of $\delta$ (respectively, $\delta^\sigma$) 
with 
certain conjugates of
$u_1\omega_1+\cdots+u_n\omega_n$. 
Thus we can write our map in the form 
\[
\u\mapsto(\Lambda_1(\u)\cdots \Lambda_{n/2}(\u)\,,\,
\Lambda_{n/2+1}(\u)\cdots \Lambda_n(\u)),\] 
with suitable linear forms $\Lambda_i$, which will be
linearly independent. Hence if we set $v_i=\Lambda_i(\u)$ for $1\le i\le n$ 
it will suffice to show that the map 
\[\v\mapsto(v_1\cdots v_{n/2}\,,\,v_{n/2+1}\cdots v_n)\]
is non-singular whenever $v_1\cdots v_n\not=0$.
This however is trivial.

Since $\bN_{K/\QQ}(\u^{(\RR)})\not=0$, we deduce that there are
indices $i$ and $j$ such that
\[\frac{\partial\bN_1(\u^{(\RR)})}{\partial u_i}
\frac{\partial\bN_2(\u^{(\RR)})}{\partial u_j}> 
\frac{\partial\bN_1(\u^{(\RR)})}{\partial u_j}
\frac{\partial\bN_2(\u^{(\RR)})}{\partial u_i}.\] 
We suppose for notational simplicity that $i=1$ and $j=2$, the other
cases being similar.
By continuity we then have
\begin{equation}\label{om1}
J(\u):=\frac{\partial\bN_1(\u)}{\partial u_1}
\frac{\partial\bN_2(\u)}{\partial u_2}- 
\frac{\partial\bN_1(\u)}{\partial u_2}
\frac{\partial\bN_2(\u)}{\partial u_1}\gg U^{n-2}
\end{equation}
throughout $\cU$, if $G$ is large enough. We now write $\v=(u_1,u_2)$
and $\u=(\v,\w)$, where $\w=(u_3,\ldots,u_n)$.  We will find it
convenient to number the entries in $\w$ as $(w_3,\ldots,w_n)$.
By the Implicit Function Theorem, if 
\[\bN_1(\v,\w)=x_1,\;\;\; \bN_2(\v,\w)=x_2\]
with $(\v,\w)\in\cU$, we can express $\v$ in terms of $\x$ and $\w$ as
$\v=\v(\x,\w)$. 

We may now calculate that
\begin{align*}
\meas\{\u\in\cU &: \delta\bN_{K/L}(\u)\in R\}\\
&= \int_{\x\in R}\int_{\w\in W(\x)}J(\v(\x,\w),\w)^{-1}\d w_3\cdots\d
w_{n}\, \d x_1\d x_2,
\end{align*}
where 
\[W(\x):=\{\w\in\RR^{n-2}: (\v(\x,\w),\w)\in\cU\}.\]
We therefore define
\begin{equation}\label{eq:ouromega}
\omega(\x):=M^n\int_{\w\in W(\x)}J(\v(\x,\w),\w)^{-1}\d w_3\cdots\d w_{n},
\end{equation}
so that \eqref{omint} immediately follows.  It is clear from the
definition of $W(\x)$ that $\omega$ is supported on the set of values
$\delta\bN_{K/L}(\u)$ as $\u$ runs over $\cU$.  Thus the support is
contained in a disc of radius $O(U^{n/2})$ as claimed in the lemma. Moreover
\begin{equation}\label{501}
\meas\{W(\x)\}\le 
\meas\{\w\in\RR^{n-2}: (\v,\w)\in\cU\;\;\mbox{for some}\;\v\}
\ll U^{n-2},
\end{equation}
whence \eqref{om1} yields $\omega(\x)\ll 1$ as required.

The function $\v(\x,\w)$ will be continuously
differentiable with respect to $\x$ and $\w$.  Moreover it will be
weighted-homogeneous of degree 1 in $\w$ and degree $2/n$ in $\x$.
Since $\partial x_k/\partial
v_l\ll U^{n/2-1}$ for $1\le k,l\le 2$ we deduce from \eqref{om1} that
\begin{equation}\label{pv}
\frac{\partial v_k}{\partial x_l}\ll J(\u)^{-1}U^{n/2-1}\ll U^{1-n/2},
\;\;\; (1\le k,l\le 2).
\end{equation} 
We observe from \eqref{todo} that $\partial
\bN_{K/\QQ}(\u^{(\RR)})/\partial w_3$ does not vanish, whence 
\[\partial\bN_j(\u^{(\RR)})/\partial w_3\not=0\]
for at least one of $j=1$ or
$j=2$.  Moreover
\[\frac{\partial \bN_j(\v(\x,\w),\w)}{\partial w_3}=
\frac{\partial x_j}{\partial w_3}=0\]
for $j=1,2$, whence
\[
\frac{\partial \bN_j}{\partial v_1}\frac{\partial v_1}{\partial w_3}+
\frac{\partial \bN_j}{\partial v_2}\frac{\partial v_2}{\partial w_3}
+\frac{\partial \bN_j}{\partial w_3}=0,\;\;\;(j=1,2).
\]
Thus at least one of $\partial v_1/\partial w_3$ and
$\partial v_2/\partial w_3$ is non-vanishing at $\u^{(\RR)}$.  We
suppose that $\partial v_1/\partial w_3\not=0$, the alternative 
case being similar.  By continuity we then have 
$|\partial v_1/\partial w_3|\gg 1$ for $|\u-\u^{(\RR)}|\ll G^{-1}$, if
$G$ is large enough.  It follows that
\begin{equation}\label{om3}
\left|\frac{\partial v_1}{\partial w_3}\right|\gg 1
\end{equation}
throughout $\cU$, since the partial derivative is homogeneous in $\w$,
of weight zero.

We can now investigate $\omega(\x+\hh)-\omega(\x)$.  On $W(\x+\hh)\cap
W(\x)$ we have
\begin{align*}
J(\v(\x+\hh,\w),\w)^{-1}- &J(\v(\x,\w),\w)^{-1}\\
&\ll 
U^{4-2n}|J(\v(\x+\hh,\w),\w)-J(\v(\x,\w),\w)|,
\end{align*}
by \eqref{om1}.  Moreover
\[J(\v(\x+\hh,\w),\w)-J(\v(\x,\w),\w)\ll U^{n-3}
|\v(\x+\hh,\w)-\v(\x,\w)|,\]
since $J(\u)$ is a form in $\u$ of degree $n-2$.  It then follows
from \eqref{pv} that
\[J(\v(\x+\hh,\w),\w)^{-1}-J(\v(\x,\w),\w)^{-1}\ll U^{2-3n/2}|\hh|\]
on $W(\x+\hh)\cap W(\x)$.  We therefore see from \eqref{501} that
the contribution to $\omega(\x+\hh)-\omega(\x)$ from the set 
$W(\x+\hh)\cap W(\x)$ will be $O(U^{-n/2}|\hh|)$. This is
satisfactory.

On the remaining
range we merely use the bound $J^{-1}\ll U^{2-n}$. It therefore
suffices to estimate the measure of the set of $\w\in\RR^{n-2}$ for
which either $(\v(\x+\hh,\w),\w)\in\cU$ or $(\v(\x,\w),\w)\in\cU$, but
not both.  By substituting $\x'=\x+\hh$ and $\hh'=-\hh$ if necessary,
we may suppose that $(\v(\x,\w),\w)\in\cU$ and 
$(\v(\x+\hh,\w),\w)\not\in\cU$.  In view of \eqref{pv} this means that 
$(\v(\x,\w),\w)$ lies in $\cU$ at
a distance $O(U^{1-n/2}|\hh|)$ from the boundary of $\cU$.  Let
$\cU_{\x}(\hh)$ denote the set of such points $\w$. As yet we have not
completely specified the set $\cU$ in \eqref{eq:Abox}, 
and it is now time to do so. For indices $i\ge 3$ we merely choose
$L_i(\u)=w_i$.  Then if $(\v(\x,\w),\w)$ lies in $\cU$ at a distance
$O(U^{1-n/2}|\hh|)$ from the edge defined by $L_i$ we see that the
corresponding $w_i$ runs over an interval of length
$O(U^{1-n/2}|\hh|)$, so that the contribution to ${\rm meas}(\cU_{\x}(\hh))$
is $\ll U^{1-n/2}|\hh|U^{n-3}=U^{n/2-2}|\hh|$.  We take the remaining
linear forms $L_i$ to be $v_1$ and $v_1+\lambda v_2$.  Here $\lambda$ is
a non-zero constant chosen sufficiently small that
\[\left|\frac{\partial (v_1+\lambda v_2)}{\partial w_3}\right|\gg 1\]
throughout $\cU$.  In view of \eqref{om3} such a choice will be
possible, since we have $\partial v_2/\partial w_3\ll 1$.
Suppose now that $w_4,w_5,\ldots,w_{n}$ are fixed. Then 
if $(\v(\x,\w),\w)$ lies in $\cU$ at a distance
$O(U^{1-n/2}|\hh|)$ from the edge defined by the linear form $L_i=v_1$ 
we see that $w_3$ is confined to an interval of length
$O(U^{1-n/2}|\hh|)$, and similarly for the edge defined by
$L_i=v_1+\lambda v_2$.  It follows that the contribution to 
${\rm meas}(\cU_{\x}(\hh))$ is $O(U^{n/2-2}|\hh|)$ in these cases too. Since
$J^{-1}\ll U^{2-n}$ we deduce that $\omega(\x+\hh)-\omega(\x)\ll
U^{-n/2}|\hh|$ as required.

It remains to establish the lower bound for $\omega(\x)$.  It is clear
that $J(\u)$ is homogeneous of degree $n-2$, whence $J(\u)\ll U^{n-2}$
for all relevant $\u$.  In view of \eqref{eq:ouromega} it therefore
suffices to show that $\meas\{W(\x)\}\gg G^{2-n}U^{n-2}$ on a suitable
set of values $\x$.
Now 
\[\left(\v(\delta\bN_{K/L}(\u^{(\RR)}),u^{(\RR)}_3,\ldots,u^{(\RR)}_n)\,,
\,u^{(\RR)}_3,\ldots,u^{(\RR)}_n\right)=\u^{(\RR)}.\]
Moreover if $|\x|\ll U^{n/2}$ and $|\w|\ll U$ then
\begin{align*}
\Big|\v(\x,\w)
&-\v\left(U^{n/2}\delta
\bN_{K/L}(\u^{(\RR)}),Uu^{(\RR)}_3,\ldots,Uu^{(\RR)}_n\right)
\Big| 
\\
&\ll  U\left\{|U^{-n/2}\x-\delta\bN_{K/L}(\u^{(\RR)})|+
\max_{3\le i\le n}|U^{-1}w_i-u^{(\RR)}_i|\right\},
\end{align*}
by the homogeneity properties of $\v(\x,\w)$. We therefore see 
from the definition \eqref{eq:Abox} of $\cU$ that 
there is a small constant $c>0$ such that $(\v(\x,\w),\w)\in\cU$
whenever we have $|w_i-Uu^{(\RR)}_i|\le cG^{-1}U$ for $3\le i\le n$ and 
\[|\x-U^{n/2}\delta\bN_{K/L}(\u^{(\RR)})|\le cG^{-1}U^{n/2}.\]
We therefore have $\meas\{W(\x)\}\gg G^{2-n}U^{n-2}$ in the above
region, as required.

\section{A large sieve bound for $\alpha_0(x)$}

From Lemma \ref{Mainapp3} we know that $\al_0(x)$ is evenly distributed
in all congruence classes for moduli $q\le Q$.  
The condition that $Q\le U^{1/2}$ causes no problems.  However a more
serious constraint on the size of $Q$ comes from the fact that
we cannot handle $\hat{\al}$ if $Q$ is too large.  

The goal of the
present section is to show that in fact $\al_0$ is evenly distribution
for 
``almost all'' congruence classes, for much larger values of $q$.
This will enable us to get equidistribution for large moduli, on
average, while keeping $Q$ sufficiently small that $\hat{\al}$ can be
adequately handled.

Our equidistribution result will be achieved by a quite general 
large sieve argument,
motivated by (but, we believe, simpler than) that used by Fouvry and
Iwaniec \cite[\S~10]{FI}. 
A related procedure is given by Iwaniec and Kowalski 
\cite[Theorem~17.5]{IK}.  However we have found it more 
convenient to use additive characters directly, rather than to 
switch to multiplicative ones as they do. The reader might 
care to note that a somewhat different argument, also of a 
very general kind, appears in Heath-Brown \cite[\S~2]{RMI}.

We assume that we have a function $\al_0(x)$ defined on $\Lring$,
satisfying 
\begin{equation}\label{alass}
\left|\twosum{x\in R}{x\equiv y\bmod{q}}\al_0(x)\right|\le W_0,
\end{equation}
for all $q\le Q$ and $y\in\Lring$, where $R$ is a square of side 
$2N$ as in Lemma \ref{L22}.   For the duration of this section let us write 
$$
\Sigma(S;q,y):=\twosum{x\in S}{x\equiv y\bmod{q}}\al_0(x),
$$
for any square $S$. We proceed to consider
\[
S_1(Q_0):=\sum_{q\le Q_0}q^2\sum_{y\bmod{q}}
\left|\Sigma(R;q,y)\right|^2.
\]
If
\[K(t):=\sum_{x\in R}\al_0(x)e^{(L)}(tx)\]
for $t\in L$ then
\[\twosum{x\in R}{x\equiv y\bmod{q}}\al_0(x)=
q^{-2}\sum_{b\bmod{q}}e_q^{(L)}(-by)K(b/q),\]
and we deduce that
\[S_1(Q_0)=\sum_{q\le Q_0}\;\sum_{b\bmod{q}}|K(b/q)|^2.\]
We now write the fraction $b/q$ in lowest terms as $c/h$, say. Then
\begin{align*}
S_1(Q_0)&=\sum_{h\le Q_0}\;\Osum_{c\bmod{h}}|K(c/h)|^2
\#\{q\le Q_0: h\mid q\}\\
&\le  Q_0 \sum_{h\le Q_0}h^{-1}\Osum_{c\bmod{h}}|K(c/h)|^2.
\end{align*}
The contribution from terms $h>Q$ is at most
\[Q_0Q^{-1}\sum_{h\le Q_0}\;\Osum_{c\bmod{h}}|K(c/h)|^2
\le Q_0Q^{-1}(\sqrt{2N}+Q_0)^4\sum_{x\in R}|\al_0(x)|^2,\]
by the two dimensional large sieve in the form \eqref{ls2}.

For the remaining terms with $h\le Q$ we observe that
\begin{align*}
\Osum_{c\bmod{h}}|K(c/h)|^2 \le\sum_{c\bmod{h}}|K(c/h)|^2
&= h^2\sum_{y\bmod{h}}\left|\Sigma(R;q,y)\right|^2\\ 
&\le  h^4W_0^2,
\end{align*}
by our assumption \eqref{alass}.  Thus terms with $h\le Q$ contribute
at most $Q_0Q^4W_0^2$ to $S_1(Q_0)$.  This enables us to conclude as follows.

\begin{lemma}\label{LS2}
Under the assumption \eqref{alass} we have
\[\sum_{q\le Q_0}q^2\sum_{y\bmod{q}}
\left|\twosum{x\in R}{x\equiv y\bmod{q}}\al_0(x)\right|^2
\ll Q_0Q^4W_0^2+\frac{Q_0(N^2+Q_0^4)}{Q}\sum_{x\in R}|\al_0(x)|^2.\]
\end{lemma}

We will require a form of this estimate in which we have a maximum
over different squares $R$.
We assume that $\al_0(x)$ is supported on a set
\[S=\{a+b\tau\in \Lring:\, (a,b)\in (-N,N]^2\}\]
and proceed to cover this with $K^2$ smaller squares, each contained
in $S$ and of the type
\[R_i=\{a+b\tau\in \Lring:\, (a-a_i\,,\, b-b_i)\in (-N/K, N/K]^2\},\]
for appropriate pairs $(a_i,b_i)$.  Here $K\le N$ is a positive
integer parameter which we will specify in due course.
Now any square $R\subseteq
S$, with sides aligned with the coordinate axes, will include a 
union of certain of the squares $R_i$, outside which there are only 
$O(N^2K^{-1})$ points of $\Lring$.  If we now require that 
\begin{equation}\label{al0bnd}
|\al_0(x)|\le A_0
\end{equation}
for all $x$, then
\[\left|\Sigma(R;q,y)\right|\leq O(N^2K^{-1}A_0)+
\sum_{i\le K^2}\left|\Sigma(R_i;q,y)\right|.\]
Thus if
\[S_2(Q_0):=\sum_{q\le Q_0}q^2\sum_{y\bmod{q}}\max_R
\left|
\Sigma(R;q,y)\right|^2,\]
then
\[S_2(Q_0)\ll N^4K^{-2}A_0^2Q_0^5+K^2\sum_{i\le K^2}
\sum_{q\le Q_0}q^2\sum_{y\bmod{q}}\left|\Sigma(R_i;q,y)\right|^2.\]
Lemma \ref{LS2} then yields
\begin{align*}
\sum_{q\le Q_0}q^2&\sum_{y\bmod{q}}
\left|
\Sigma(R_i;q,y)\right|^2
\ll 
Q_0Q^4W_0^2+\frac{Q_0(N^2K^{-2}+Q_0^4)}{Q}\sum_{x\in R_i}|\al_0(x)|^2.
\end{align*}
Since 
\[\sum_{i\le K^2}\sum_{x\in R_i}|\al_0(x)|^2\leq \sum_{x\in R}|\al_0(x)|^2,\]
we deduce that
\[S_2(Q_0)\ll N^4K^{-2}A_0^2Q_0^5+K^4Q_0Q^4W_0^2+
\frac{K^2Q_0(N^2K^{-2}+Q_0^4)}{Q}\sum_{x\in S}|\al_0(x)|^2.\]
We now choose $K=A_0Q_0^3$, which yields the following conclusion.
\begin{lemma}\label{LS3}
Under the assumptions \eqref{alass} and \eqref{al0bnd}
we have
\begin{align*}
\sum_{q\le Q_0}q^2
&\sum_{y\bmod{q}}\max_R
\left|\twosum{x\in R}{x\equiv y\bmod{q}}\al_0(x)\right|^2\\
&\ll N^4Q_0^{-1}+A_0^4Q_0^{13}Q^4W_0^2
+Q_0Q^{-1}N^2\sum_{x\in  S}|\al_0(x)|^2,
\end{align*}
providing that $A_0Q_0^5\le N$.
\end{lemma}

We apply this last result to our particular situation.  In view of
the conditions in Lemmas~\ref{Mainapp3} and \ref{L23} we take 
$N$ of order $U^{n/2}$, 
$Q^{n+5}\leq U^{1-\eta}$
and $G\leq U^{1/(n+1)}$, so that the value
$W_0=Q^{n+1}U^{n-1+\eta}$ is admissible. In order to estimate
$\al_0(x)$ we note that $\al(x)\ll_{\eta} U^{\eta}$ for any $\eta>0$.
Moreover \eqref{hataldef} yields $|\hat{\al}(x)|\le Q^3\omega(x)$,
since \eqref{comb2} and \eqref{one2} imply
\begin{equation}
  \label{eq:solar}
\sum_{z\bmod{q}}\rho(z,q)= 1.
\end{equation}
Since $\omega(x)\ll 1$ by Lemma \ref{omegaprops}, we will certainly
have $\al_0(x)\ll_{\eta} Q^3U^{\eta}$.  Taking $A_0\ll_\eta
Q^3U^\eta$, the right hand 
side in Lemma \ref{LS3} is therefore
\[\ll_{\eta} U^{2n}Q_0^{-1}+Q_0^{13}Q^{2n+18}U^{2n-2+O(\eta)}+ 
Q_0Q^{-1}U^{2n+\eta}G^{O(1)},\] 
by Lemma \ref{L23}. This enables us to conclude as follows.

\begin{lemma}\label{LS4}
If $Q\le Q_0\le U^{1/(n+16)}$ and $G\leq U^{1/(n+1)}$, then
\[\sum_{q\le Q_0}q^2\sum_{y\bmod{q}}\max_R
\left|\twosum{x\in R}{x\equiv y\bmod{q}}\al_0(x)\right|^2
\ll_{\eta} Q_0Q^{-1}U^{2n+O(\eta)}G^{O(1)},\]
for any $\eta>0$.
\end{lemma}

We end by establishing a trivial bound for the above sum, which provides
an instructive comparison.  If $q\le Q_0\le U^{n/2}$ then 
\begin{align*}
\left|\sum_{\substack{x\in R \\ x\equiv y\bmod{q}}}\al_0(x)\right|^2
&\ll
\#\{x\in R:\, x\equiv y\bmod{q}\} 
\sum_{\substack{x\in R \\ x\equiv y\bmod{q}}}|\al_0(x)|^2\\
&\ll
U^nq^{-2}\sum_{\substack{x\in R \\ x\equiv y\bmod{q}}}|\al_0(x)|^2,
\end{align*}
whence 
\begin{equation}\label{502}
\begin{split}
\sum_{y\bmod{q}}\max_R
\left|\twosum{x\in R}{x\equiv y\bmod{q}}\al_0(x)\right|^2
&\ll U^nq^{-2}\sum_{y\bmod{q}} \;
\sum_{x\equiv y\bmod{q}}|\al_0(x)|^2\\
&=
U^nq^{-2}\sum_{x\in \Lring}|\al_0(x)|^2.
\end{split}
\end{equation}
We therefore obtain the trivial bound
\[\sum_{q\le Q_0}q^2\sum_{y\bmod{q}}\max_R
\left|\twosum{x\in R}{x\equiv y\bmod{q}}\al_0(x)\right|^2
\ll_{\eta}Q_0U^{2n+\eta}G^{O(1)},\]
via Lemma \ref{L23}. Thus Lemma \ref{LS4} provides a saving which is a power of
$Q$, providing that $Q$ is larger than a suitable power of $U^{\eta}G$.

\section{Bilinear forms in dimension $2$}\label{s:bil}

The estimation of bilinear forms is one of the cornerstones of
analytic number theory and can be traced back to the work of
Vinogradov. Given finite sequences $u_m,v_n \in \CC$ and 
a matrix $\mathbf{A}=(a_{m,n})$ of complex numbers, 
the essential problem is to estimate the double
sum
$$
\u^T \mathbf{A}\v=\sum_{m}\sum_n u_mv_n a_{m,n}.
$$
The standard procedure 
is to use Cauchy's inequality to remove the dependence
on $\u=(u_m)$. This gives 
$$
|\u^T \mathbf{A}\v|\leq \left(\sum_{m}|u_m|^2\right)^{1/2}
\left(\sum_m \left|\sum_n 
v_n a_{m,n}\right|^2\right)^{1/2}.
$$
For the second term on the right one 
expands the square and reverses the order of
summation, so as to use suitable information about the sum
$$
\sum_{m}a_{m,n_1}\overline{a_{m,n_2}}.
$$
It may happen that the sum itself is small.  In other cases one can
give an asymptotic evaluation with a main term 
$M(n_1,n_2)$ say.  One
may then complete the analysis via an estimation of the sum
$
\sum_{n_1}\sum_{n_2} v_{n_1}\overline{v_{n_2}}M(n_1,n_2).
$

In our work we will require an analogue of this procedure for
sequences indexed by elements of $\ZZ^2$, rather than by
$\ZZ$.  In this endeavour we are inspired by 
the analytic machinery developed for the Gaussian integers 
$\ZZ[i]$ by Fouvry and Iwaniec \cite[\S~9]{FI}.  
We will provide a completely self-contained account of the method,
which has the 
advantage of being slightly more general in scope.

Let $\mathbf{M}\in{\rm SL}_2(\ZZ)$ and let
$\al,\be : \ZZ^2 \rightarrow \CC$ 
be arbitrary functions with finite $L^2$-norms 
$$
\| \al \|_2 := \left(\sum_{\x\in \ZZ^2}
|\al(\x)|^2\right)^{1/2}, \quad
\| \be \|_2 := \left(\sum_{\y\in \ZZ^2}
|\be(\y)|^2\right)^{1/2}.
$$
Furthermore, let 
$\la:\ZZ\rightarrow \CC$ be such that 
$$
\| \la \|_2 := \left(\sum_{l\in \ZZ} |\la(l)|^2\right)^{1/2}
$$
is also finite. Our objective
is to estimate the double sum
$$
S(\al,\be;\la):=\sum_{\x \in \ZZ^2}\sumE_{\y\in
  \ZZ^2} \al(\x)\be(\y)
\la\big(\x^T\mathbf{M}\y\big),
$$
where $\sum\!{}^{(E)}$ indicates that the summation over
$\y$ is restricted by the condition  $\gcd(y_1,y_2)\le E$.

In the analysis of this section the implied constants in our
$O(\cdot)$ and $\ll$ notation will be allowed to depend implicitly on the
coefficients of the matrix $\mathbf{M}$.  However they will be uniform
in all other parameters.
It will be convenient to introduce a norm on 
elements $\x\in\RR^2$ via $|\x|=\max\{|x_1|,|x_2|\}$. We should
observe that this particular choice of norm is important, since it
will be significant for us that the unit ball is a square.

Let $A,B\geq 1$ with $A\geq B$. We will think of $A$ and $B$ as being
large, with $A/B$ also large, but of considerably smaller order than
$B$.   We will suppose that
$\al,\be$ are 
supported on the sets
\begin{equation}\label{eq:supp}
\cA
=\{\x \in \ZZ^2: |\x|\leq A\}, \quad 
\cB=\{\y \in \ZZ^2: B\leq |\y| \leq 2B\},
\end{equation}
respectively.  Let  $M_0$ denote the maximum of the moduli of the
entries of $\mathbf{M}$, and let
$$
\widetilde{\mathbf{M}}
=\left(
\begin{array}{cc}
 0 & 1\\ 
-1 & 0
\end{array}
\right)\mathbf{M}.
$$
We now define
\[\widetilde{B}=
\{\z\in\ZZ^2:(2M_0)^{-1}B\leq |\z|\leq 4M_0B\}
\]
and we
observe that $\widetilde{\mathbf{M}}\y\in\widetilde{B}$ whenever $\y\in\cB$.

An application of  Cauchy's inequality now gives
$$
|S(\al,\be;\la)|^2\leq \| \be \|_2^2 \cdot \| \la \|_2^2 
\sumE_{\y\in \cB}
 \sum_{l\in \ZZ}
\left|\sum_{\substack{\x \in \ZZ^2\\ \x^T\mathbf{M}\y=l}
}\al(\x)\right|^2.
$$
Enlarging the range of summation for $\y$ we
therefore deduce that
\begin{equation}\label{eq:*1}
|S(\al,\be;\la)|\leq \| \be \|_2 \cdot \| \la \|_2 
\cdot T(\al)^{1/2},
\end{equation}
where
$$
T(\al):=
\sumE_{\substack{\y\in \ZZ^2\\ 
\widetilde{\mathbf{M}}\y\in\widetilde{B}}}
\sum_{l\in
  \ZZ}\left|\sum_{\substack{\x \in \ZZ^2\\ \x^T\mathbf{M}\y=l}
}\al(\x)\right|^2.
$$
Opening up the inner sum we obtain
\begin{align*}
T(\al)=\sumE_{\substack{\y\in \ZZ^2\\ 
\widetilde{\mathbf{M}}\y\in\widetilde{B}}}
\sum_{\substack{\w \in \ZZ^2\\ \w^T\mathbf{M}\y=0}}(\al*\al)(\w), 
\end{align*}
with
$$
(\al*\al)(\w):=
\sum_{\substack{\x,\x'\in \ZZ^2 \\ \x-\x'=\w}} \al(\x)\overline{\al}(\x').
$$

The vectors $\y$ take the form $e\u$ where 
$\u=(u_1,u_2)\in\ZZ^2$ is primitive, and where $e\le E$ is 
a positive integer.  Since 
$\det \mathbf{M}=1$, we see that
the general solution of the linear equation $\w^T\mathbf{M}\u=0$ is
$$
\w=c\widetilde{\mathbf{M}}\u
$$
for $c\in \ZZ$. Hence 
\begin{align*}
T(\al)=\sum_{e\le E}\;
\Osum_{\substack{\u\in \ZZ^2\\ e\widetilde{\mathbf{M}}\u\in\widetilde{B}}}\;
\sum_{c\in \ZZ}(\al*\al)\big(c\widetilde{\mathbf{M}}\u\big)=
\sum_{e\le E}\;\Osum_{e\z\in\widetilde{B}}\;\sum_{c\in \ZZ}(\al*\al)(c\z),
\end{align*}
where $\sum_{\z}^*$ denotes summation for primitive 
vectors $\z=(z_1,z_2)$.  We note that Cauchy's inequality yields 
$|(\al*\al)(\w)|\leq \| \al \|_2^2$ for any $\w\in\ZZ^2$.
Thus the contribution to the above sum from terms with $c=0$ is 
$O(\| \al \|_2^2 B^2 )$, whence 
\begin{equation}
  \label{eq:*2}
T(\al)=T_0(\al)+O(\| \al \|_2^2 B^2 ),
\end{equation}
with 
$$
T_0(\al):=\sum_{e\le E}\;\Osum_{e\z\in \widetilde{B}}\;
\sum_{c\in\ZZ\setminus\{0\}}(\al*\al)(c\z).
$$
We now use the M\"obius function to pick out the condition 
$\gcd(z_1,z_2)=1$, giving
\[T_0(\al)=\sum_{e\le E}\sum_{c\in\ZZ\setminus\{0\}}
\sum_{b=1}^{\infty}\mu(b)T_0(\al;b,c,e),\]
where
\[T_0(\al;b,c,e):=\sum_{\substack{ec^{-1}\w \in \widetilde{B}\\ bc\mid\w}}
(\al*\al)(\w).\]

For $ec^{-1}\w\in\widetilde{B}$ we deduce that 
\begin{equation}\label{zb}
(2M_0)^{-1}e^{-1}cB\le|\w|\le 4M_0e^{-1}cB.
\end{equation}
Moreover
we will have $(\al*\al)(\w)=0$ unless 
$|\w|\le 2A$.  Thus we can restrict the summation over $c$ to the range
\[0<|c|\le 4M_0 AB^{-1}E=C,\]
say.  

We may now handle certain ranges of $b$ and $c$ trivially.  To this
end we give ourselves a parameter $K\ge 1$ which we will choose in due
course. In view of \eqref{zb} the number of available vectors 
$\w\equiv0\bmod{bc}$,
corresponding to a particular triple $b,c,e$ is $O(B^2b^{-2}e^{-2})$.  
Since  $|(\al*\al)(\w)|\leq \| \al \|_2^2$ we deduce that the contribution
to $T_0(\al)$ from $b>K$ is 
$$
\ll \sum_{e\le E}\sum_{c\leq C}\sum_{b>K}\| \al \|_2^2B^2b^{-2}e^{-2}
\ll \| \al \|_2^2ABEK^{-1}.
$$
Similarly the contribution from $|c|\le CK^{-1}$ is  
$$
\ll \sum_{e\le E}\sum_{|c|\leq C/K}\sum_{b}  
\| \al \|_2^2B^2b^{-2}e^{-2}
\ll \| \al \|_2^2ABEK^{-1}.
$$

It follows that
\[T_0(\al)=\sum_{e\le E}\sum_{b\leq K}\mu(b)\sum_{C/K<|c|\leq C} 
T_0(\al;b,c,e)+O(\| \al \|_2^2ABEK^{-1}).\]
Combining this with
\eqref{eq:*1} and \eqref{eq:*2} therefore leads to the
conclusion that 
$$
S(\al,\be;\la) \ll  \| \be \|_2 \cdot \| \la\|_2 
\left\{\| \al \|_2 B+\| \al \|_2 (ABEK^{-1})^{1/2}+T_1(\al)^{1/2}\right\}
$$
for any $A\ge B\ge 1$ and $K\ge 1$, with
\[T_1(\al):=
\sum_{e\le E}\;\sum_{b\le K}\;\sum_{C/K<|c|\le C} |T_0(\al;b,c,e)|.\]

We now open up the convolution $\al*\al$ to obtain
$$
T_0(\al;b,c,e)
=\sum_{\x'\in \ZZ^2}\overline{\al}(\x')
\sum_{\substack{\x\in R \\ \x\equiv \x' \bmod{bc}}} \al(\x),
$$
where 
\[R=\{\x\in \ZZ^2: (2M_0)^{-1}e^{-1}cB\leq |\x-\x'|\leq 4M_0e^{-1}cB\}.\]
Thus
\[|T_0(\al;b,c,e)|\le \sum_{\x'\in \ZZ^2}|\al(\x')|\max_R
\left|\sum_{\substack{\x\in R \\ \x\equiv \x' \bmod{bc}}} \al(\x)\right|,\]
where $R$ now runs over all squares with sides aligned to the axes. It
follows from Cauchy's inequality that
\begin{align*}
|T_0(\al;b,c,e)|^2&\le  \|\al\|_2^2\sum_{\x'\in\cA}\max_R
\left|\sum_{\substack{\x\in R \\ \x\equiv \x' \bmod{bc}}} \al(\x)\right|^2\\
&\ll \|\al\|_2^2(1+A^2(bc)^{-2})T_2(\al;|bc|), 
\end{align*}
with
\begin{equation}\label{T2d}
T_2(\al;q):=\sum_{\u\bmod{q}}\max_R
\left|\sum_{\substack{\x\in R \\ \x\equiv \u\bmod{q}}}\al(\x)\right|^2.
\end{equation}
If we insist that $E\le B$ then $1+A|bc|^{-1}\le 1+AKC^{-1}\ll BKE^{-1}$,
whence
\[T_1(\al)\ll E\sum_{b\le K}\sum_{C/K<c\le C}\|\al\|_2BKE^{-1}
T_2(\al;|bc|)^{1/2}.\]
Hence, if $\tau(q)$ denotes the usual divisor function, and we set
\begin{equation}\label{T3d}
T_3:=\sum_{q\le CK}q^2T_2(\al;q),
\end{equation}
a further application of Cauchy's inequality yields
\begin{align*}
T_1(\al)&\ll BK\|\al\|_2\left\{\sum_{C/K<q\le CK}\tau(q)^2q^{-2}\right\}^{1/2}
T_3^{1/2}\\
&\ll  BK\|\al\|_2\left\{KC^{-1}(\log C)^3\right\}^{1/2}T_3^{1/2}\\
&\ll  (BK)^{3/2}(AE)^{-1/2}(\log A)^{3/2}\|\al\|_2T_3^{1/2}.
\end{align*}
We may now conclude as follows.

\begin{lemma}\label{lem:bil}
Let $A\ge B\geq E\ge 1$ and define $C=4M_0AB^{-1}E$.
Let $\al,\be : \ZZ^2 \rightarrow \CC$  be
functions supported on the sets \eqref{eq:supp}, and
let $\la:\ZZ\rightarrow\CC$. Define $T_2(\al;q)$ and $T_3$ as 
in \eqref{T2d} and \eqref{T3d}.
Then for any $K\ge 1$ we have 
\[
S(\al, \be;\la)
\ll \|\be \|_2 \cdot  \|\la\|_2 
\left(\|\al\|_2 B +\| \al \|_2
  (ABEK^{-1})^{1/2}+T_1(\al)^{1/2}\right),
\]
with
\[T_1(\al)\ll (BK)^{3/2}(AE)^{-1/2}(\log A)^2\|\al\|_2T_3^{1/2}.\]
\end{lemma}

Before preceding further it may be helpful to comment on the above estimate.
For the purposes of this illustration we shall suppose that $K$ is
chosen so that 
\begin{equation}\label{kc}
\max\{E^2,\log^2 A\}\le K\le A/B.
\end{equation}
If we merely estimate $S(\al, \be;\la)$ via Cauchy's inequality,
using the fact that $\x^T\mathbf{M}\y=l$ has $O(AB\log A)$ solutions $\x,\y$  
for each given $l$, 
we are led to a trivial bound 
\[
S(\al, \be;\la)\ll\|\al\|_2\cdot \|\be\|_2\cdot 
\|\la\|_{2}\sqrt{AB \log A}.
\] 
Hence the first term in
the above lemma saves at least $\sqrt{A/B}$. Similarly the second
term saves at least $\sqrt{K/E}$.  Both these are at least
$K^{1/4}$.  To analyse the third
term we use the argument in \eqref{502} to deduce that  
\[T_2(\al;q)\ll A^2q^{-2}\|\al\|_2^2\] 
for $q\ll A$.  We therefore obtain the trivial bound  
\begin{align*} 
T_1(\al)&\ll  
\|\al\|_2^2 A^{1/2}B^{3/2}C^{1/2}K^2E^{-1/2}\log^2 A\\  
&\ll  \|\al\|_2^2 ABK^2\log^2 A\\
&\ll \|\al\|_2^2 ABK^3.  
\end{align*}
For comparison, in order for the third term in 
Lemma \ref{lem:bil} to produce a 
comparable saving to that in the first two terms, 
one would wish to replace
the above by 
\[T_1(\al)\ll \|\al\|_2^2 ABK^{-1},\]
say.  The type of saving we require is exactly that given by Lemma
\ref{LS4}, providing that we work with $S(\al_0,\be;\la)$.

\bigskip

We now show how Lemma \ref{lem:bil} can be applied to estimate
\[\cl{E}(G,H,V)=\sum_{x \in \Lring }\sum_{y \in \Lring} 
\al_0(x)\be(y)\la\big(\str(x,y)\big),\]
when $G=\log V$. 
We have yet to specify the parameter $Q$ used in the definition of
$\hat{\al}$, but we shall do this shortly. 
Under suitable hypotheses 
we will see in \S~\ref{s:M} and \S~\ref{s:ssi}  
that
the main term  $\cl{M}(G,H,V)$ has order at least
$(\log V)^{1-3n}H^nV^{2n}$.  Our goal
is to show that $\cl{E}(G,H,V)$ is smaller than this by at least
a power of $H$.

If $x=x_1+x_2\tau$ and $y=y_1+y_2\tau$, then $\str(x,y)=\x^T\mathbf{M}\y$,
where
\[\mathbf{M}=\left(\begin{array}{cc}0 & -1\\ 
1 & 0 \end{array}\right).\]
By abuse of notation we will write $\al_0(\x)=\al_0(x)$ 
and $\be(\y)=\be(y)$.
From our definitions of $\al, \hat{\al}$ and $\omega$ in
\eqref{eq:AL}, \eqref{hataldef} and \eqref{eq:ouromega} we see that if
$\al_0(x)\not=0$ then $x=\delta\bN_{K/L}(\u)$  for some $\u\in\cU$.
Moreover, from \eqref{eq:Abox} we deduce that
$\delta\bN_{K/L}(\u)=U^{n/2}\delta\bN_{K/L}(\u^{(\RR)})+o(U^{n/2})$ as 
$G\rightarrow\infty$.  Since
$\delta\bN_{K/L}(\u^{(\RR)})\not=0$ it follows that $\al(\x)$ is
supported on a suitable set $|\x|\le A$ with $A=c_1 U^{n/2}$ for a
certain constant $c_1>0$.  We may
analyse the support of $\be$ similarly.  Since
$\bN_{K/L}(\v^{(\RR)})\not=0$ we deduce that $\be(\y)$ is
supported on a set $B\le|\y|\le 2B$ with $B=c_2 V^{n/2}$ for a
suitable positive constant $c_2$.  Moreover, $\beta(\y)$ is also
supported on vectors $\y$ with $\gcd(y_1,y_2)\ll 1$, by
\eqref{Ebnd}.  We may therefore take $E$ of order $1$.

Estimates for $\|\be\|_2$ and $\|\la\|_2$ are given by Lemma 
\ref{lem:upper-abl}, while Lemma \ref{L23} handles $\|\al_0\|_2$. 
Inserting these
bounds into Lemma \ref{lem:bil} and adopting the assumption
\eqref{kc}, together with the bounds
$Q^{n+5}\le U^{1-\eta}$ and $G\leq U^{1/(n+1)}$,
we then see that
\[B+(ABEK^{-1})^{1/2}\ll U^{n/4}V^{n/4}K^{-1/2}\]
and
\[
(BK)^{3/2}(AE)^{-1/2}(\log A)^2\ll U^{-n/4}V^{3n/4}K^{5/2},
\]
so that
\begin{align*}
\cl{E}(G,H,V)\ll_{\eta}~& V^{(n+\eta)/2}W^{(n+\eta)/2}
\Big(U^{(n+\eta)/2}\cdot U^{n/4}V^{n/4}K^{-1/2}\\
&+ U^{-n/8}V^{3n/8}K^{5/4}\cdot
U^{(n+\eta)/4}T_3^{1/4}\Big)G^{O(1)}.
\end{align*}

We now apply Lemma \ref{LS4}, along with the bounds $C\ll AB^{-1}E$
and $G=O_\eta(V^\eta)$,
to deduce that
\[T_3=\sum_{q\le CK}q^2T_2(\al;q)
\ll_{\eta} CKQ^{-1}U^{2n}V^{O(\eta)}\ll_{\eta}
KQ^{-1}U^{5n/2}V^{-n/2+O(\eta)},\] 
providing that 
\begin{equation}\label{kc'}
Q\le CK\le U^{1/(n+16)}.
\end{equation}
Note that when $G=\log
V$ the condition $G\leq U^{1/(n+1)}$ is automatic for large
$V$. Moreover, if $Q\le U^{1/(n+16)}$ then we automatically have
$Q^{n+5}\le U^{1-\eta}$, if $\eta$ is small enough.

Thus, on recalling that $U=HV$ and $W=H^{1/2}V$, we deduce that
\begin{align*}
\cl{E}(G,H,V)&\ll_{\eta}H^{n/4}V^{n+O(\eta)}
\left(H^{3n/4}V^{n}K^{-1/2}
+H^{3n/4}V^{n}K^{3/2}Q^{-1/4}\right)\\
&= H^{n}V^{2n+O(\eta)}(K^{-1/2}+K^{3/2}Q^{-1/4}).
\end{align*}
It is now apparent that we should take $K=Q^{1/8}$, which will satisfy
\eqref{kc} and \eqref{kc'} if $Q\ll H^{4n/7}$ and 
$(\log V)^{16}\ll Q\ll H^{-4n}U^{8/(n+16)}$. We conclude as follows.

\begin{lemma}\label{lem:treat_E}
Let $G=\log V$ and assume that 
\[(\log V)^{16}\ll Q\ll  
\min\{H^{n/2}, H^{-5n}U^{8/(n+16)}\}.\]
 Then we have 
$$
\cl{E}(G,H,V)\ll_{\eta} Q^{-1/16}H^{n}V^{2n+O(\eta)}.
$$
\end{lemma}

\section{Estimation of the main term}\label{s:M}

The purpose of this section and the next 
is to produce a satisfactory estimate for the 
sum
$$
\cl{M}=\cl{M}(G,H,V)=
\sum_{x \in \Lring }\sum_{y \in \Lring}
\hat{\al}(x)\be(y)\la\big(\str(x,y)\big),
$$
as $G\rightarrow \infty$, 
where $\hat{\al}$ is the approximation to $\al$ which we constructed in
\eqref{hataldef}. 
Our goal will be to demonstrate that $\cl{M}$ has
order $G^{-2n }H^nV^{2n}$, which simple heuristics suggest to be
the expected size of $\cl{N}(G,H,V)$.  In fact we will fall a little
short of this, showing that 
$\cl{M}\gg G^{1-3n}H^nV^{2n}$ when suitable constraints are placed
on the parameters $Q,G,H$ and $V$.

Let 
\begin{equation}
\begin{split}\label{eq:woodruff}
\mathcal{V}_1
&:=\{\v\in \cV\cap\ZZ^{n}: 
\v\equiv \v^{(M)} \bmod{M}, ~
\mbox{\eqref{Econd} holds}\},\\
\mathcal{W}_1&:=\{\w\in \cW\cap\ZZ^{n}: 
\w\equiv \w^{(M)} \bmod{M}\},
\end{split}
\end{equation}
where $\cV, \cW$ are given by \eqref{eq:Abox}. 
Opening up the functions $\beta$ and $\la$ from \eqref{eq:BE} and
\eqref{eq:LA}, respectively, we find that
$$
\cl{M}=
\sum_{\v\in \cV_1}\sum_{\w\in \cW_1}
\sum_{\substack{x \in \Lring \\
\tr_{L/\QQ}(x\bN_{K/L}(\v))=
2\bN_{K/\QQ}(\w)}}
\hspace{-0.2cm}
\hat{\al}(x)=
\sum_{\v\in \cV_1}\sum_{\w\in \cW_1}
\cl{M}(\v,\w),
$$
say.
Recall that $D_L=\tau-\tau^\sigma$ where $\{1,\tau\}$ is a 
$\ZZ$-basis for $\Lring$.
It will be
convenient to make the observation that 
$$
2\tr_{L/\QQ}(\tau^2)-\tr_{L/\QQ}(\tau)^2=D_L^2,
$$
which is a non-zero integer. 
By enlarging the weak approximation set $S$ in Theorem \ref{thmain-a},
if necessary, we may clearly assume that $M$ contains any prime
divisors of $2D_L^2$.

We now set in motion our analysis of $\cl{M}(\v,\w)$, for given $\v\in
\cV_1$ and $\w\in \cW_1$.
Suppose that $\bN_{K/L}(\v)$ decomposes as 
$\bN_1(\v)+\bN_2(\v)\tau$, for suitable forms $\bN_1,\bN_2$ of degree
$n/2$. Then \eqref{Econd} demands that that
$\bN_1(\v)$ and $\bN_2(\v)$ be coprime.
 Moreover, in view of \eqref{vone2} and
our convention that $\omega_1=1$ in 
the integral basis $\{\omega_1, \ldots, \omega_n\}$
for $K$ over $\QQ$, it is clear that 
$$
\bN_1(\v)\equiv 1 \bmod{M}, \quad 
\bN_2(\v)\equiv 0 \bmod{M}.
$$
Let us write $x=x_1+x_2\tau$. 
The constraint in the $x$
summation is equivalent to 
$$
a_1x_1+a_2x_2=b,
$$
for integers
$a_1,a_2,b$ such that 
\begin{equation}\label{eq:aib}
\begin{split}
a_1&=\tr_{L/\QQ}(\bN_{K/L}(\v))=2\bN_1(\v)+\tr_{L/\QQ}(\tau)\bN_2(\v), \\
a_2&=
\tr_{L/\QQ}(\tau\bN_{K/L}(\v))=
\tr_{L/\QQ}(\tau)\bN_1(\v)+\tr_{L/\QQ}(\tau^2)\bN_2(\v),\\
b&=2\bN_{K/\QQ}(\w).
\end{split}
\end{equation}
In order to analyse $\cl{M}(\v,\w)$ we will need to recall the
expression for $\hat{\alpha}(x)$.  Let the function
$\omega(x)=\omega(x_1,x_2)$ be given by \eqref{eq:ouromega}, with 
key properties as in Lemma~\ref{omegaprops}. 
Then for any parameter $Q\geq 1$ and $x\in \Lring$ 
we have 
$$
\hat{\alpha}(x)=\omega(x)\sum_{q\le Q}\; \Osum_{t\bmod{q}}
e_q^{(L)}(-tx)\sum_{z\bmod{q}}\rho(z,q)e_q^{(L)}(tz),
$$
where  $\rho(z,q)$ is given by \eqref{eq:vetch} and 
the notation $\osum$ means that the sum is taken over
$t=t_1+t_2\tau
$, with $t_1,t_2 \in \ZZ/q\ZZ$ such that
$\gcd(q,t_1,t_2)=1$.
Define 
\begin{equation}\label{eq:c_q}
c_q(t):=\sum_{z \bmod{q}}\rho(z,q) 
e_q^{(L)}(tz),
\end{equation}
for $t\bmod{q}$. 
It is clear from \eqref{eq:solar} that $|c_q(t)|\leq 1$. 
Inserting our expression for $\hat{\al}(x)$ and breaking the inner sum
into congruence classes modulo $q$, we see that 
$$
\cl{M}(\v,\w)= 
\sum_{q\le Q}\;\Osum_{t\bmod{q}}c_q(t)
\sum_{r\bmod{q}}e_q^{(L)}(-tr)\cl{L}_r(\v,\w),
$$
where
$$
\cl{L}_r(\v,\w):=\sum_{\substack{a_1x_1+a_2x_2=b \\
x\equiv r\bmod{q}}}\omega(x).
$$

We proceed to investigate the summation conditions on $x$.  We
first prove that $a_1$ and $a_2$ are not both zero, that 
$b$ is non-zero, and that
\begin{equation}\label{eq:kappa}
\gcd(a_1,a_2)=2^\kappa, \quad \kappa:= \begin{cases}
1, & \mbox{if $2\mid \tr_{L/\QQ}(\tau)$,}\\
0, & \mbox{if $2\nmid \tr_{L/\QQ}(\tau)$.}
\end{cases}
\end{equation}
In particular $\gcd(a_1,a_2)\mid b$. 

To establish the claim we observe
that $\gcd(a_1,a_2)$ is a common divisor of 
$$
2a_2-\tr_{L/\QQ}(\tau)a_1=D_L^2\bN_2(\v)
$$
and 
$$
\tr_{L/\QQ}(\tau^2)a_1-\tr_{L/\QQ}(\tau)a_2=D_L^2\bN_1(\v).
$$
However $\v$ and $\w$ are in the regions \eqref{eq:Abox} 
with $\bN_{K/L}(\v^{(\RR)})$ and
$\bN_{K/\QQ}(\w^{(\RR)})$ non-zero.  Thus if $G$ is large enough we will have
$\bN_{K/L}(\v)$ and $\bN_{K/\QQ}(\w)$ non-zero.  Hence $b$ will be non-zero and
similarly, since $D_L\not=0$, the numbers $a_1$ and $a_2$ cannot both be zero.
We also see that
$$
\gcd(a_1,a_2)\mid D_L^2\gcd(\bN_1(\v),\bN_2(\v))=D_L^2.
$$ 
Since
$\bN_1(\v)\equiv 1 \bmod{p}$ and 
$\bN_2(\v)\equiv 0 \bmod{p}$ for any prime divisor $p$ of $2D_L^2$, the
claim readily follows.

We now write $a_i'=2^{-\kappa}a_i$ for $i=1,2$, and $b'=2^{-\kappa}b$,
so that the first condition on $x$ becomes $a_1'x_1+a_2'x_2=b'$.  Thus
the sum will be empty unless $a_1'r_1+a_2'r_2\equiv b'\bmod{q}$, as we
now assume.  For future use we note that this is equivalent to demanding that 
\begin{equation}\label{eq:rcond}
2^{-\kappa}\tr_{L/\QQ}(r\bN_{K/L}(\v))\equiv 
2^{1-\kappa}\bN_{K/\QQ}(\w) \bmod{q}.
\end{equation}
We may therefore write $a_1'r_1+a_2'r_2=b'+qk$ for some
$k\in\ZZ$.  Then, if we take $x=r+qy$, the summation conditions on $x$
translate into the requirement that $a_1'y_1+a_2'y_2=-k$.
Since $a_1'$ is coprime to $a_2'$  we can find integers  
$\overline{a_1},\overline{a_2}$ such that 
$a_1'\overline{a_1}+a_2'\overline{a_2}=1$.  We then have
$a_1'(y_1+k\overline{a_1})+a_2'(y_2+k\overline{a_2})=0$, 
so that
$y_1+k\overline{a_1}=-ma_2'$ and $y_2+k\overline{a_2}=ma_1'$ for some
integer $m$.  It follows that the summation conditions
$a_1x_1+a_2x_2=b$ and $x\equiv r\bmod{q}$ are satisfied if and only
if $x_1$ and $x_2$ take the forms $x_1=r_1-\overline{a_1}kq-a_2'qm$ and
$x_2=r_2-\overline{a_2}kq+a_1'qm$ respectively. Our conclusion is
therefore that
$$
\mathcal{L}_r(\v,\w)=\sum_{m\in \ZZ}f(m),
$$
if (\ref{eq:rcond}) holds, where 
$$
f(m)
:=\omega\big(r_1-\overline{a_1}kq-a_2'qm,
r_2-\overline{a_2}kq+a_1'qm\big). 
$$

We would now like to replace the discrete summation over $m$ by a
continuous integral.
For this a relatively crude approach is available to us through 
Lemma~\ref{omegaprops}. Thus if $v\in [0,1]$ it follows that
$$
f(m+v)-f(m) \ll U^{-n/2}q\max\{|a_1'|,|a_2'|\}.
$$
Moreover Lemma 
\ref{omegaprops} tells us that $\omega$ is supported on a disc of
radius $O(U^{n/2})$, whence $f$ is supported on an interval with
length $O(U^{n/2}/(q\max\{|a_1'|,|a_2'|\}))$. Recalling that
$U\geq V$ in Lemma \ref{lem:suffice} and 
$\max\{|a_1'|,|a_2'|\}\ll V^{n/2}$, this therefore produces the conclusion
$$
\left|\int_{-\infty}^\infty f(m)\d m - \sum_{m\in \ZZ}f(m)\right|\ll
1+q U^{-n/2}\max\{|a_1'|,|a_2'|\}\ll q.
$$
Assuming that $a_2\neq 0$ the change of variables 
$m=q^{-1}(-2^\kappa x+r_1/a_2'-\overline{a_1}kq/a_2')$ now yields
$$
\mathcal{L}_r(\v,\w)=
\frac{2^\kappa}{q}I(\v,\w)+O(q ),
$$
for $r$ satisfying \eqref{eq:rcond}, 
with 
\begin{equation}
  \label{eq:defJ}
I(\v,\w):=
\int_{-\infty}^\infty
\omega(a_2x,-a_1x+b/a_2)\d x.
\end{equation}
Here $a_1,a_2,b$ depend on $\v$ and $\w$ and are given by \eqref{eq:aib}.
If $a_2=0$ we  reverse the r\^oles of $a_1$ and $a_2$ to 
produce an integral involving
$\omega(-a_2 x+b/a_1,a_1 x)$.

Recall the estimate $|c_q(t)|\leq 1$ that we recorded above.
Inserting our estimate for $\mathcal{L}_r(\v,\w)$ into that 
for $\cl{M}(\v,\w)$, it now follows that 
\begin{align*}
\cl{M}(\v,\w)= 2^\kappa I(\v,\w)
\sum_{q\le Q}\frac{1}{q}\Osum_{t\bmod{q}}c_q(t)
\sum_{\substack{r\bmod{q}\\ \mbox{\scriptsize{\eqref{eq:rcond}
        holds}}}} e_q^{(L)}(-tr)
+O\left(Q^6\right).
\end{align*}
We now sum both sides over all $\v \in \cV_1$ and $\w\in \cW_1.$
On observing that $\#\cV_1=O(G^{-n}V^n)$ and $\#\cW_1=O(G^{-n}W^n)$,
we see that the overall contribution from the error term is 
\begin{align*}
&\ll Q^6G^{-2n}V^nW^n \ll Q^6H^{n/2}V^{2n}.
\end{align*}
This will be satisfactory if $Q$ is sufficiently small compared to  $H$.
Our work so far has shown that 
\begin{equation}\label{addn}
\cl{M}= 2^\kappa \sum_{\v\in \cV_1}\sum_{\w\in \cW_1}I(\v,\w)
\sum_{q\le Q}\frac{1}{q}\cl{C}+O\left(Q^6 H^{n/2}V^{2n}\right),
\end{equation}
where $\cV_1,\cW_1$ are given by \eqref{eq:woodruff}, and
\[\cl{C}:=\Osum_{t\bmod{q}}c_q(t)
\sum_{\substack{r\bmod{q}\\ \mbox{\scriptsize{\eqref{eq:rcond}
        holds}}}} e_q^{(L)}(-tr).\]

Opening up \eqref{eq:c_q}, we find that
\[\cl{C}=\sum_{\substack{r\bmod{q}\\ 
\mbox{\scriptsize{\eqref{eq:rcond} holds}}}} 
\sum_{z\bmod{q}} \rho(z,q)\Osum_{t\bmod{q}} e_q^{(L)}\left(t(z-r)\right).\]
We have seen that the condition \eqref{eq:rcond} can be written
$a_1'r_1+a_2'r_2\equiv b'\bmod{q}$, for non-zero integers
$a_1',a_2',b'$ such that $\gcd(a_1',a_2' )=1$.
The inner sum is a Ramanujan sum and thus it follows from
\eqref{eq:gen-ram}, combined with \eqref{comb2}, that 
\begin{align*}
\cl{C}
&=\sum_{u\mid q}u^2 \mu(q/u)
\sum_{\substack{r \bmod{q}\\
a_1'r_1+a_2'r_2\equiv b'\bmod{q}}}\;
\sum_{\substack{z\bmod{q}\\ z\equiv r \bmod{u}}} \rho(z,q)\\
&=\sum_{u\mid q}u^2 \mu(q/u)
\sum_{\substack{r \bmod{q}\\
a_1'r_1+a_2'r_2\equiv b'\bmod{q}}}
\rho(r,u).
\end{align*}
Let $q=uv$. Writing $r=r'+ur''$ for $r'\bmod{u}$ and $r''\bmod{v}$ we
see that the sum over $r$ is equal to
$$
\sum_{\substack{r' \bmod{u}\\
a_1'r_1'+a_2'r_2'\equiv b'\bmod{u}}}
\rho(r',u)
\#\left\{
r''\bmod{v}: 
a_1'r_1''+a_2'r_2''\equiv g(r') \bmod{v}
\right\},
$$
where $g(r')=u^{-1}(
b'-a_1'r_1'-a_2'r_2')$.
Since 
$a_1'$ and $a_2'$ are coprime 
it follows that there are precisely $v$
possibilities for $r''.$ The definition \eqref{eq:vetch} of
$\rho(r',u)$ therefore reveals that 
$$
\cl{C}=M^n q\sum_{u\mid q}\frac{u \mu(q/u)}{[M,u]^n} 
\#\left\{\b{s}\bmod{[M,u]}:~
\begin{array}{l}
\b{s}\equiv\u^{(M)}\bmod{M},\\
2^{-\kappa}F(\v;\w;\b{s}) \equiv 0\bmod{u}
\end{array}\right\},
$$
where 
$$
F(\v;\w;\b{s})
:= \tr_{L/\QQ}(\delta\bN_{K/L}(\b{s})\bN_{K/L}(\b{v}))-2\bN_{K/\QQ}(\w).
$$
This is the polynomial that underpins the variety introduced 
in \eqref{TR3}.

We proceed to insert this expression for $\cl{C}$ into \eqref{addn}.
We will use the M\"obius function to 
remove the coprimality condition 
\eqref{Econd}, which is implicit in the definition \eqref{eq:woodruff}
of $\cV_1$. We therefore arrive at the estimate
\begin{equation}\label{eq:twayblade}
\begin{split}
\cl{M}=~& 2^\kappa M^n 
\sum_{k=1}^\infty \mu(k)\sum_{q\le Q}\sum_{u\mid q}\frac{u \mu(q/u)}{[M,u]^n} 
\sum_{\substack{\b{s}\bmod{[M,u]}\\
\b{s}\equiv\u^{(M)}\bmod{M}}} \\
&\times \sum_{\substack{(\v,\w)\in \cV_2\times \cW_1\\
2^{-\kappa}F(\v;\w;\b{s}) \equiv 0\bmod{u}}}I(\v,\w)
+O\left(Q^6H^{n/2}V^{2n}\right),
\end{split}
\end{equation}
where $I(\v,\w)$ is given by \eqref{eq:defJ} and 
$\cV_2$ is defined as for $\cV_1$, but with the condition 
\eqref{Econd} replaced by $k \mid \bN_{K/L}(\v)$. This  latter condition
is taken componentwise as $k\mid \bN_1(\v)$ and $k\mid \bN_2(\v)$, in
the usual way.  

The sum over $k$ is empty when $k\gg V^{n/2}$ so that we only need
consider values $k\ll V^{n/2}$.  However we
need to reduce this range further.  Since
$\bN_{K/L}(\v^{(\RR)})\neq 0$ it follows from \eqref{eq:aib} and
the definition of $\mathcal{V}$ that 
\[\max\{|a_1|, |a_2|\}\gg V^{n/2}.\] 
The integration in \eqref{eq:defJ} is therefore 
over an
interval of length $\ll V^{-n/2}U^{n/2}=H^{n/2}$, by Lemma~\ref{omegaprops}.
A second application of Lemma~\ref{omegaprops} to bound the size
of $\omega$ now yields
\begin{equation}
  \label{eq:speedwell}
I(\v,\w) \ll H^{n/2}. 
\end{equation}
Hence there is an absolute constant $c>0$
such that the overall contribution to $\cl{M}$ from terms with $k
>K$ is 
\begin{align*}
&\ll H^{n/2}
\sum_{K<k\ll V^{n/2}} 
\sum_{q\le Q}
\sum_{u\mid q} u 
\sum_{\substack{\w\in \cW_1}}
M_k(cV)\\
&\ll_\eta Q^{2+\eta} H^{n/2}W^n
\sum_{K<k\ll V^{n/2}} 
M_k(cV),
\end{align*}
where
\[
M_k(X):= \#\{\v \in \ZZ^n:
|\v|\leq X, ~k \mid \bN_{K/L}(\v)
\}.
\]
Clearly $k^2 \mid \bN_{K/\QQ}(\v)$ whenever 
$k \mid \bN_{K/L}(\v)$, whence
$$
M_k(X)\leq  \sum_{\substack{ N\ll X^n\\ k^2\mid N}} \#\{\v \in \ZZ^n:
|\v|\leq X, ~\bN_{K/\QQ}(\v)=N\}\ll_\eta \left(1+\frac{X^n}{k^2} \right)X^\eta.
$$
The overall contribution to the main term from $k
>K$ is therefore seen to be 
\begin{align*}
&\ll_\eta Q^{2+\eta} H^{n/2}W^nV^{n+\eta}
\sum_{K<k\ll V^{n/2}} 
\frac{1}{k^2}\\
&\ll_\eta
 \frac{Q^{2+\eta} H^{n}V^{2n+\eta}}{K}.
\end{align*}
This allows us to truncate
the summation over $k$ in \eqref{eq:twayblade} 
with acceptable
error.  

We now replace the sum over $\b{s}$ in 
 \eqref{eq:twayblade} by one in which the variable runs
modulo $[M,u,k]$.  This has the effect of multiplying it by
$[M, u,k]^{-n}[M,u]^n$.
We may therefore summarise our findings in the following result.

\begin{lemma}\label{lem:stage1} 
Let $K\geq 1$ and let $\kappa$ be
given by \eqref{eq:kappa}. 
Then we have 
\begin{align*}
\cl{M}=~& 2^\kappa M^n
\sum_{k\leq K}
\mu(k)\sum_{q\le Q}\sum_{u\mid q}\frac{u \mu(q/u)}{\Delta^n} 
\sum_{\substack{\b{s}\bmod{\Delta}\\
\b{s}\equiv\u^{(M)}\bmod{M}}} 
\mathcal{K}_{k,u}(\b{s})\\
&+
O_\eta\left(H^{n}V^{2n+\eta}\left\{\frac{Q^6}{H^{n/2}}+
 \frac{Q^{2+\eta} }{K}\right\}\right),
\end{align*}
where $\Delta:=[M,u,k]$ and 
\begin{align*}
\mathcal{K}_{k,u}(\b{s})
:=
\sum_{\substack{(\v,\w)\in \cV_2\times \cW_1\\
2^{-\kappa}F(\v;\w;\b{s}) \equiv 0\bmod{ u}}}I(\v,\w).
\end{align*}
Here $\cW_1$ is given by
\eqref{eq:woodruff},
 $I$ is given by \eqref{eq:defJ} and 
$\cV_2$ 
is defined as for $\cV_1$ in \eqref{eq:woodruff}, but with 
\eqref{Econd} replaced by $k \mid \bN_{K/L}(\v)$. 
 \end{lemma}

The next phase of the argument concerns a detailed analysis of 
$\mathcal{K}_{k,u}(\b{s})$. 
It is natural to break the sum over $\v$ and $\w$ into 
congruence classes modulo
$\D=[M, u,k]$. Let 
\[\mathcal{I}(\Delta)=
\mathcal{I}(\b{p},\b{q},\b{s};\Delta):=
\sum_{\substack{\w\in \cW\cap \ZZ^n\\ \w\equiv \b{q}\bmod{\Delta}}}
\sum_{\substack{\v\in \cV\cap\ZZ^{n}\\
\v\equiv \b{p}\bmod{\Delta}}} 
I(\v,\w).\]
Then we have 
\begin{equation}
  \label{eq:toast}
\mathcal{K}_{k,u}(\b{s})
=\sum_{\b{p},\b{q}}\mathcal{I}(\Delta),
\end{equation}
where the sum is over 
$(\b{p},\b{q})\in (\ZZ/\Delta \ZZ)^{2n}$ for which 
\begin{equation}
  \label{eq:daisy}
(\b{p},\b{q},\b{s})\equiv (\v^{(M)},\w^{(M)},\u^{(M)}) \bmod{M}
 \end{equation}
and 
\begin{equation}\label{eq:violet}
\begin{split}
& 2^{-\kappa}F(\b{p};\b{q};\b{s})\equiv 0\bmod{u},\quad
 k \mid \bN_{K/L}(\b{p}).
\end{split}
\end{equation}
Perhaps the most obvious way to deal with $\mathcal{I}(\Delta)$
is to approximate the sums over $\w$ and
$\v$ by integrals. 
A simple change of variables would then permit us to extract the
dependence of 
$\mathcal{I}(\Delta)$ on $\b{p},\b{q},\b{s}$ and
$\Delta$. Instead of this it turns
out that we can manage with a relatively crude direct
comparison of $\mathcal{I}(\Delta)$ with
$\mathcal{I}(1)$, as we proceed to show.

It follows from our conditions of summation that $\Delta\ll k u \leq
KQ$.  
It will be convenient  to make the additional hypothesis
\begin{equation}
  \label{eq:hyp1}
  KQ\leq V,
\end{equation}
which implies in particular that $\Delta\ll V$. 
In view of \eqref{eq:defJ} we have
\begin{align*}
\mathcal{I}(\Delta)
&= \int_{-\infty}^\infty
\sum_{\substack{\w\in \cW\cap \ZZ^n\\ \w\equiv \b{q}\bmod{\Delta}}}
\sum_{\substack{\v\in \cV\cap\ZZ^{n}\\
\v\equiv \b{p}\bmod{\Delta}}} 
f(\v,\w) \d x,
\end{align*}
where if $a_i=a_i(\v)$ and $b=b(\w)$ are given by \eqref{eq:aib} then
$$
f(\v,\w):=
\omega\left(a_2(\v)x,
-a_1(\v)x+\frac{b(\w)}{a_2(\v)}\right).
$$
We will denote by $T(\b{p},\b{q};\Delta)$ the integrand that appears
in this expression for $\mathcal{I}(\Delta)$. Using an
approach based on the proof of Lemma \ref{Mainapp3} in \S~\ref{app2}
we compare
$T(\b{p},\b{q};\Delta)$ with 
$T(\b{0},\b{0};\Delta)$.
For this we will assume without loss of generality that 
$\max_{i=1,2} |a_i(\v)|=|a_2(\v)|$
in the sum over $\v$.  The alternative
possibility is accommodated by a simple change of variables in the
integral over $x$.  Hence the definition of $\cV$ ensures that 
$a_2(\v)$ has order of magnitude 
$V^{n/2}$ and by 
Lemma~\ref{omegaprops} we have $x\ll H^{n/2}$
in  $\mathcal{I}(\Delta)$. 

Under the change of variables 
$\w=\w'+\b{q}$  we see that 
$$
\left|\frac{b(\w)}{a_2(\v)}-
\frac{b(\w')}{a_2(\v)}\right|
\ll \frac{\Delta|\w'|^{n-1}}{|a_2(\v)|}
\ll \frac{\Delta W^{n-1}}{V^{n/2}},
$$
since $|\w'|\ll \Delta+W\ll W$ by \eqref{eq:hyp1}.
We may therefore conclude from Lemma~\ref{omegaprops} that 
\begin{align*}
T(&\b{p},\b{q};\Delta)
=
\sum_{\substack{\v\in \cV\cap\ZZ^{n}\\
\v\equiv \b{p}\bmod{\Delta}}} 
\sum_{\substack{\w'+\b{q}\in \cW\cap \ZZ^n\\ \w'\equiv \b{0}\bmod{\Delta}}}
\left\{
f(\v,\w') + O\left(\frac{\Delta W^{n-1}}{U^{n/2}V^{n/2}} 
\right)\right\}.
\end{align*}
Note that the number of $\v$ appearing in the outer sum is $O(\D^{-n}V^n)$.
For $\w'$ such that $\w'+\b{q}\in \cW$ it is clear that 
$\w'\in \cW$ unless $\w'$ is within a distance $\Delta$ of the boundary of
$\cW$. Invoking Lemma \ref{omegaprops} to deduce that $f$
is bounded absolutely, and recalling that $U=HV$ and $W=H^{1/2}V$, 
it easily follows that
\begin{align*}
\left|
T(\b{p},\b{q};\Delta)-
T(\b{p},\b{0};\Delta)\right|
&\ll (\Delta^{-1}V)^{n}(\Delta^{-1}W)^{n-1}+
\frac{\Delta^{-2n+1}V^{n/2}W^{2n-1}}{ U^{n/2}} \\
&\ll  \Delta^{-2n+1} H^{(n-1)/2}V^{2n-1},
\end{align*}
by \eqref{eq:hyp1}.

We now repeat the above process by 
considering the effect of a change of variables 
$\v=\v'+\b{p}$ in  
$T(\b{p},\b{0};\Delta)$. 
Recalling 
that $\Delta\ll KQ \leq V$, by \eqref{eq:hyp1}, we obtain
\begin{align*}
\left|a_i(\v)-a_i(\v')\right| 
&\ll \Delta |\v'|^{n/2-1} \ll \Delta V^{n/2-1},
\end{align*}
for $i=1,2$.  In particular it follows that
\begin{align*}
\left|
\frac{b(\w)}{a_2(\v)}-\frac{b(\w)}{a_2(\v')}
\right|
&\ll \frac{\Delta V^{n/2-1} W^n}{
|a_2(\v)a_2(\v')|}\\
&\ll \frac{\Delta V^{n/2-1} W^n}{V^{n}}\\
&= \Delta H^{n/2}V^{n/2-1},
\end{align*}
for any $\w \in \cl{W}$.
An application of Lemma \ref{omegaprops} reveals that 
$$
T(\b{p},\b{0};\Delta)=
\sum_{\substack{\w\in \cW\cap \ZZ^n\\ \w\equiv \b{0}\bmod{\Delta}}}
\sum_{\substack{\v'+\b{p}\in \cV\cap\ZZ^{n}\\
\v'\equiv \b{0}\bmod{\Delta}}} 
\left\{f(\v',\w) 
+ O\left( \frac{\Delta H^{n/2}V^{n/2-1}}{U^{n/2}}\right)\right\}.
$$
To control the error term we note that the total number of available
$\v',\w$ in the sums is 
$$
\ll \frac{W^n}{\Delta^n} \cdot  \left(\frac{\Delta+V}{\Delta}+1\right)^n
\ll \Delta^{-2n}V^nW^n.
$$
We conclude that
\begin{align*}
\left|T(\b{p},\b{0};\Delta)
-
T(\b{0},\b{0};\Delta)\right|
&\ll (\Delta^{-1}W)^n (\Delta^{-1}V)^{n-1} + 
\frac{\Delta^{-2n+1}H^{n/2}V^{3n/2-1}W^n}{U^{n/2}}\\
&\ll \Delta^{-2n+1}H^{n/2}V^{2n-1},
\end{align*}
whence
$$
T(\b{p},\b{q};\Delta)=
T(\b{0},\b{0};\Delta)+O\left(
\Delta^{-2n+1}H^{n/2}V^{2n-1}\right).
$$

Summing the latter estimate over $\b{p},\b{q}$ modulo $\Delta$ gives  
$$
T(\b{0},\b{0};1)=
\Delta^{2n}T(\b{0},\b{0};\Delta)+O\left(
\Delta H^{n/2}V^{2n-1}\right),
$$
which we substitute back in to get
$$
T(\b{p},\b{q};\Delta)=
\Delta^{-2n}T(\b{0},\b{0};1)+O\left(
\Delta^{-2n+1}H^{n/2}V^{2n-1}\right).
$$
Now the outer integral in our expression for 
$\mathcal{I}(\Delta)$
is over an interval of length $O(H^{n/2})$. 
We have therefore arrived at the following result.

\begin{lemma}\label{lem:J}
Assume that $KQ\leq V$. Then we have 
$$
\mathcal{I}(\Delta)=
\Delta^{-2n}
\mathcal{I}(1)+O(
\Delta^{-2n+1}H^{n}V^{2n-1}).
$$
\end{lemma}

This result will allow us to separate out what is in effect the 
``singular integral'' associated to our counting problem. 
In the notation of \eqref{eq:defJ} it is given by 
\begin{equation}\label{eq:sig-inf}
\begin{split}
\sigma_\infty(G,H,V)
&:=\mathcal{I}(1)
=
 \sum_{\substack{\w\in \cW\cap \ZZ^n}}
\sum_{\substack{\v\in \cV\cap\ZZ^{n}}}
I(\v,\w).
\end{split}
\end{equation}
It follows from \eqref{eq:speedwell} that
\begin{equation}
  \label{eq:upper-sig}
\sigma_\infty(G,H,V)\ll (G^{-1}V)^n\cdot (G^{-1}W)^n \cdot H^{n/2}=
G^{-2n}H^n V^{2n}.
\end{equation}
It is interesting to compare the present situation with the 
singular integrals arising 
from typical applications of the Hardy--Littlewood circle
method. These are expressed as volumes that reflect the real density
of solutions. It transpires that we will be able to provide a
lower bound for $\sigma_\infty(G,H,V)$ which essentially matches the upper
bound \eqref{eq:upper-sig} without first approximating the sum by an
integral. Nonetheless crucial use will be made of the fact that 
$\sigma_\infty(G,H,V)$ features a sum over points close to a
non-singular real point on the variety \eqref{TR3}.

For given $k$ and $u$, let
\begin{equation}
  \label{eq:saxifrage}
  N_{M}(k,u):=\#\left\{
\b{p},\b{q},\b{s} \bmod{\Delta}
: 
\mbox{\eqref{eq:daisy}, \eqref{eq:violet} hold}
 \right\},
\end{equation}
where 
$\D=[M,u,k]$.
We claim that
\begin{equation}
  \label{eq:bugle}
  N_{M}(k,u)\ll_\xi \frac{\Delta^{3n-1+\xi}}{k},
\end{equation}
for any $\xi>0$.
To see this we write 
$A=\tr_{L/\QQ}\left(\delta \bN_{K/L}(\b{p})\bN_{K/L}(\b{s})\right)$ 
for fixed $\b{p},\b{s}$ modulo $\D$.  Then the number
of $\b{q}\bmod{\D}$ contributing to $N_{M}(k,u)$ is 
at most 
\begin{align*}
&\sum_{\substack{N\ll \D^n\\ N\equiv A \bmod{u}}} 
\#\{\b{q}\in[1,\Delta]^n: 2\bN_{K/\QQ}(\b{q})=N\}
\ll_\xi \frac{\D^{n+\xi/2}}{u}.
\end{align*}
Moreover there are $\D^n$ possibilities 
for $\b{s}\bmod{\D}$ and the number of available 
$\b{p} \bmod{\D}$ is at most
\begin{equation}\label{sfo1}
\sum_{\substack{N\ll \D^n\\ N\equiv 0 \bmod{k^2}}} 
\#\{\b{p}\in[1,\D]^n: \bN_{K/\QQ}(\b{p})=N\}
\ll_\xi \frac{\D^{n+\xi/2}}{k^2}.
\end{equation}
It therefore follows that 
\begin{align*}
N_{M}(k,u)&\ll_\xi \frac{\D^{n+\xi/2}}{u} \cdot \D^n
\cdot \frac{\D^{n+\xi/2}}{k^2} = \frac{\D^{3n+\xi}}{uk^2}.
\end{align*}
Noting that $uk \gg
\Delta$, we easily arrive at \eqref{eq:bugle}.

The extraction of $\sigma_\infty=\sigma_\infty(G,H,V)$ now follows on 
combining Lemmas~\ref{lem:stage1} and \ref{lem:J} with
\eqref{eq:toast} and \eqref{eq:bugle}. Putting these together, and 
assuming that $KQ\leq V$,  we therefore deduce that
   \begin{align*}
\cl{M}=~& \sigma_\infty 2^\kappa M^n
\sum_{k \leq K} \mu(k)
\sum_{q\le Q}
\sum_{u\mid q}\frac{ u \mu(q/u)}{\Delta^{3n}}N_{M}(k,u) \\
& + O_\eta\left(
H^{n}
V^{2n+\eta}\left\{ \frac{Q^6}{H^{n/2}}+
 \frac{Q^{2+\eta} }{K} + \frac{Q^{2+\eta}K^\eta }{V}\right\}
\right),
\end{align*}
where $ N_{M}(k,u)$ is given by \eqref{eq:saxifrage}.
Since $KQ\leq V$ it is clear that the error term is 
$$
\ll_\eta
H^{n}
V^{2n+2\eta}\left\{ \frac{Q^6}{H^{n/2}}+
 \frac{Q^{2} }{K}\right\}.
$$

We now show that the summation over $k$ can be
extended to infinity  with acceptable error. 
It follows from  \eqref{eq:upper-sig} and  \eqref{eq:bugle}
that
the error in so doing is 
\begin{align*}
&\ll_\eta G^{-2n}H^n V^{2n}
\sum_{k>K} \sum_{q\le Q}\sum_{u\mid q}\frac{u}{\Delta^{3n}}
\cdot \frac{\Delta^{3n-1+\eta/2}}{k}
\ll_\eta \frac{Q^{1+\eta}H^nV^{2n}}{K^{1-\eta}}.
\end{align*}
This error term is subsumed by that above when $KQ\leq V$.
Choosing 
$$
K=\frac{H^{n/2}}{Q},
$$
the hypothesis \eqref{eq:hyp1} becomes 
 $H^{n/2}\leq V$.  
Putting everything together we have therefore established the
following result. 

\begin{lemma}\label{lem:onion}
Assume that $H^{n/2}\leq V$. 
Then we have 
$$
\cl{M}=2^\kappa M^n\sigma_\infty \mathfrak{S}(Q)
+ O_\eta\left(Q^6 H^{n/2} V^{2n+2\eta}\right),
$$
where $\sigma_\infty=\sigma_\infty(G,H,V)$ is given by \eqref{eq:sig-inf}  and 
$$
\mathfrak{S}(Q):=\sum_{q\leq Q}\sum_{k=1}^\infty \mu(k)
\sum_{u\mid q}\frac{u \mu(q/u)}{\Delta^{3n}}N_{M}(k,u),
$$
with $\D=[M,u,k]$ and $N_{M}(k,u)$ given by \eqref{eq:saxifrage}.
\end{lemma}

\section{The singular series and integral}\label{s:ssi}

It remains to produce satisfactory lower bounds for the quantities 
$\mathfrak{S}(Q)$ and $\sigma_\infty$.  For the former our strategy
will be to show that as $Q$ tends to infinity the sum $\mathfrak{S}(Q)$
converges to a limit $\mathfrak{S}$ which is a product of local
factors.  We will then show that each of these factors is positive.
In view of (\ref{eq:bugle}) the sum
\[f(q):=\sum_{k=1}^\infty \mu(k)
\sum_{u\mid q}\frac{u \mu(q/u)}{\Delta^{3n}}N_{M}(k,u)\]
is absolutely convergent.  We rearrange it as
\[f(q)=\sum_{u\mid q}u\mu(q/u)\sum_{k=1}^\infty
\mu(k)\D^{-3n}N_{M}(k,u)\]
and write $k=k_1 k_2$ where $k_1\mid uM$ and $\gcd(k_2,uM)=1$.  Then
$[M,u,k]=[M,u]k_2$ and 
\[f(q)=\sum_{u\mid q}u\mu(q/u)\sum_{k_1\mid
  uM}\frac{\mu(k_1)}{[M,u]^{3n}} 
\sum_{k_2}\frac{\mu(k_2)}{k_2^{3n}}N_{M}(k_1k_2,u).\]
The function $N_{M}(k_1k_2,u)$ factors as
\[N_{M}(k_1k_2,u)=N_{M}(k_1,u)\#\left\{\b{p},\b{q},\b{s}\bmod{k_2}
:  k_2 \mid \bN_{K/L}(\b{p}) \right\},\]
where $\Delta=[M,u]$ in $N_M(k_1,u)$. 
We write
\[R(k):=\#\left\{\b{p}\bmod{k}: k\mid \bN_{K/L}(\b{p}) \right\}\]
and observe that 
\begin{equation}\label{RB}
R(k)\ll_\xi k^{n-2+\xi}\qquad\mbox{and}\qquad R(k)<k^n \mbox{ for }
k>1,
\end{equation}
by the argument used for (\ref{sfo1}), and the fact that
$\bN_{K/\QQ}(1,0,0,\ldots,0)=1$.  It follows that
\begin{align*}
\sum_{k_2}\frac{\mu(k_2)}{k_2^{3n}}N_{M}(k_1k_2,u)&=
N_{M}(k_1,u)\sum_{k_2}\frac{\mu(k_2)}{k_2^{n}}R(k_2)\\
&=c\prod_{p\mid uM}\left(1-\frac{R(p)}{p^n}\right)^{-1}N_{M}(k_1,u),
\end{align*}
with
\begin{equation}\label{cprod}
c=\prod_p\left(1-\frac{R(p)}{p^n}\right)>0.
\end{equation}
We therefore see that $f(q)=cf_0(q,M)$ where 
\begin{equation}\label{fex}
f_0(q,M)=\sum_{u\mid q}u\mu(q/u)\sum_{k\mid uM}
\frac{\mu(k)}{[M,u]^{3n}}
\prod_{p\mid uM}\left(1-\frac{R(p)}{p^n}\right)^{-1}N_{M}(k,u)
\end{equation}
is a multiplicative function in the two variables $q$ and $M$. We write
$q=\prod p^{\alpha}$ and $M=\prod p^\mu$, where $\mu=v_p(M)$ vanishes
for all but finitely many primes.  Then
\begin{equation}\label{abs}
\sum_{q=1}^\infty |f(q)|=c\prod_p \sum_{\alpha=0}^\infty
|f_0(p^\alpha,p^\mu)|.
\end{equation}
Providing that this product converges we will be able to deduce that
\[\lim_{Q\rightarrow\infty}\mathfrak{S}(Q)=\mathfrak{S}\]
exists.  Moreover we will have
$$
\mathfrak{S}=c\prod_p\sigma_p
$$
with
\begin{equation}\label{fsum1}
\sigma_p=\sum_{\alpha=0}^\infty f_0(p^\alpha,p^\mu),
\end{equation}
and in order to prove that $\mathfrak{S}>0$ it will suffice to show
that $\sigma_p>0$ for every prime $p$.

Our treatment of the singular series will depend on the following two
lemmas, which we will prove later in this section.  It will be
convenient to write $\x$ for the vector
$(\b{p},\b{q},\b{s})$ and $\x_0$ for
$(\b{v}^{(M)},\b{w}^{(M)},\b{u}^{(M)})$.  Moreover we shall write
$F(\x)=2^{-\kappa}F(\b{p};\b{q};\b{s})$ and
$\b{N}(\x)=\b{N}_{K/L}(\b{p})$. 
\begin{lemma}\label{Mlem}
Suppose that $p^{\mu}\| M$. For any $\beta\ge \max\{\mu,1\}$ define
$$
M(p^\beta,p^\mu):=\#\{\x\bmod{p^\beta}:\,\x\equiv\x_0\bmod{p^\mu},\, 
p^\beta\mid F(\x),\, p\nmid\N(\x)\}.
$$
Then we will have
\begin{equation}\label{5/2}
M(p,1)=p^{3n-1}+O(p^{3n-5/2}).
\end{equation}
\end{lemma}
\begin{lemma}\label{sscon}
For any prime $p$ and $\al \geq 2$ we have 
\[f_0(p^\alpha,p^\mu)\ll (2\alpha+1)^{3n}p^{\alpha-3[\alpha/2]-2}.\]
\end{lemma}

We now prove that the product (\ref{abs}) converges. Lemma \ref{sscon}
yields
\[\sum_{\alpha=2}^{\infty}|f_0(p^\alpha,p^\mu)|\ll
\sum_{\alpha=2}^{\infty}(2\alpha+1)^{3n}p^{\alpha-3[\alpha/2]-2}\ll
p^{-2}.\]
We note also that $\mu$ is non-zero only for the primes in $S$, and
that $f_0(1,1)=1$ Thus a bound $f_0(p,1)\ll p^{-3/2}$ will suffice 
to establish absolute convergence.  In general we have
\begin{align*}
\sum_{\alpha\le\beta}&f_0(p^\alpha,p^\mu)\\
&=\sum_{q\mid p^\beta}
\sum_{u\mid q}u\mu(q/u)\sum_{k\mid up^\mu}
\frac{\mu(k)}{[p^\mu,u]^{3n}}
\prod_{p\mid up^\mu}\left(1-\frac{R(p)}{p^n}\right)^{-1}N_{p^\mu}(k,u)\\
&=\sum_{u\mid p^\beta}u\sum_{k\mid up^\mu}
\frac{\mu(k)}{[p^\mu,u]^{3n}}
\prod_{p\mid up^\mu}\left(1-\frac{R(p)}{p^n}\right)^{-1}N_{p^\mu}(k,u)
\sum_{q\mid p^\beta:\,u\mid q}\mu(q/u),
\end{align*}
in view of (\ref{fex}). The final sum over $q$ vanishes unless 
$u=p^\beta$, and if
$\beta\ge\max\{\mu,1\}$ our expression reduces to
\begin{equation}\label{fsum2}
\begin{split}
\sum_{\alpha\le\beta}f_0(p^\alpha,p^\mu)&=p^{-(3n-1)\beta}
\left(1-\frac{R(p)}{p^n}\right)^{-1}\sum_{k\mid
  p}\mu(k)N_{p^\mu}(k,p^{\beta})\\
&=p^{-(3n-1)\beta}\left(1-\frac{R(p)}{p^n}\right)^{-1}M(p^\beta,p^\mu).
\end{split}
\end{equation}
Taking $\beta=1$ and $\mu=0$ we obtain
\begin{align*}
f_0(p,1)
&=
p^{1-3n}\left(1-\frac{R(p)}{p^n}\right)^{-1}\left\{M(p,1)-p^{3n-1}
\left(1-\frac{R(p)}{p^n}\right)\right\}.
\end{align*}
The required estimate $f_0(p,1)\ll p^{-3/2}$ now follows from Lemma
\ref{Mlem}, since  $R(p)/p^n\ll p^{-3/2}$ by \eqref{RB}. This
completes the proof of absolute convergence.

We turn now to proof that $\sigma_p>0$ for every prime $p$. By
(\ref{fsum2}) we have
\[\sigma_p=
\lim_{\beta\rightarrow\infty}\sum_{\alpha\le\beta}f_0(p^\alpha,p^\mu)
=\left(1-\frac{R(p)}{p^n}\right)^{-1}
\lim_{\beta\rightarrow\infty}p^{-(3n-1)\beta}M(p^\beta,p^\mu).\]
Thus $\sigma_p>0$ providing that
$M(p^\beta,p^\mu)\gg_{p,\mu}p^{(3n-1)\beta}$ as
$\beta\rightarrow\infty$. 

Using Hensel's lemma it will follow that $\sigma_p>0$ if
the variety 
\begin{align*}
F(\v;\w;\u)
&=
\tr_{L/\QQ}(\delta\bN_{K/L}(\b{u})\bN_{K/L}(\b{v}))
-2\bN_{K/\QQ}(\w)=0,
\end{align*}
introduced
in \eqref{TR3}, 
has a non-singular point over $\ZZ_p$  
which satisfies the constraints 
\begin{equation}\label{twocond}
p\nmid \bN_{K/L}(\b{v})\quad\mbox{and}\quad
(\b{v},\w,\u)\equiv (\v^{(M)},\w^{(M)},\u^{(M)}) \bmod{M}.
\end{equation}
The function $M(p,1)$ counts points over $\mathbb{F}_p$ which lie on
the above variety and satisfy the constraints (\ref{twocond}), whether
they are non-singular or not. The number of singular points will be
$O(p^{3n-2})$, whence the estimate (\ref{5/2}) shows that there will
be a suitable non-singular point providing that $p\ge p_0$, say.
We can arrange that the set $S$
includes all primes $p<p_0$, so that Lemma \ref{lem:suffice} ensures the
existence of $p$-adic integer solutions with 
$\v\equiv (1,0,\ldots,0)\bmod{p}$ and $\bN_{K/\QQ}(\w)\not=0$. It
follows that $p\nmid \bN_{K/L}(\b{v})$.  Finally we prove, by 
contradiction, that any such point must be non-singular. For otherwise
we would have
$$
\nabla_\w F(\v;\w;\u)=-2\nabla_\w \bN_{K/\QQ}(\w)=\bf{0}
$$
where $\nabla_\w=(\frac{\partial}{\partial w_1},
\ldots,\frac{\partial}{\partial w_n})$. It would then follow from 
Euler's identity that
$$
\w.\nabla_\w \bN_{K/\QQ}(\w)=n \bN_{K/\QQ}(\w)=0,
$$ 
a contradiction. Thus we have suitable non-singular points for every
prime $p$, and hence it follows that $\sigma_p>0$ for all primes $p$. 

We remark that, by combining (\ref{cprod}), (\ref{fsum1}) and
(\ref{fsum2}), we have
\[\mathfrak{S}=\prod_p\sigma_p^*,\]
with
\[\sigma_p^*:=\lim_{\beta\rightarrow\infty}
p^{-(3n-1)\beta}M(p^\beta,p^\mu).\]
Thus $\mathfrak{S}$ is a standard product of local densities.

It remains to establish Lemmas \ref{Mlem} and \ref{sscon}.  
To handle Lemma \ref{Mlem} we observe that the variety defined by
$F(\x)=0$ takes the simple form
\[X_1\cdots X_n+X_{n+1}\cdots X_{2n}+X_{2n+1}\cdots X_{3n}=0\]
over $\overline{\mathbb{F}_p}$.  This makes it clear that
we have a hypersurface of projective dimension $3n-2$, whose
singular locus has projective dimension $3n-7$.  Here we use the fact
that the singular locus consists of points where two or more
coordinates vanish from each of the sets $\{X_1,\ldots,X_n\}$, 
$\{X_{n+1},\ldots,X_{2n}\}$, and  $\{X_{2n+1},\ldots,X_{3n}\}$.
According to the result of Hooley \cite{hooley} the number of
projective points modulo $p$ differs from $(p^{3n-1}-1)/(p-1)$ by an
amount $O(p^{3n-4})$, whence the number of points in
$\mathbb{A}^n(\mathbb{F}_p)$ is $p^{3n-1}+O(p^{3n-3})$.  It follows
that
\[M(p,1)=p^{3n-1}+O(p^{3n-3})-\#\{\x\bmod{p}:\,p\mid F(\x),\, p\mid\N(\x)\}.\]
In the set on the right we have 
\[p\mid\tr_{L/\QQ}(\delta\bN_{K/L}(\b{p})\bN_{K/L}(\b{q}))
-2\bN_{K/\QQ}(\b{s})\]
and
\[p\mid\bN_{K/L}(\b{p}),\]
from which it follows that $p\mid 2 \bN_{K/\QQ}(\b{s})$.  It follows, by
the argument leading to (\ref{sfo1}) that the number
of possible $\b{p}\bmod{p}$ is $O_{\xi}(p^{n-2+\xi})$ and that the number
of possible $\b{s}\bmod{p}$ is $O_{\xi}(p^{n-1+\xi})$.  We then see that
\[\#\{\x\bmod{p}:\,p\mid F(\x),\,
p\mid\N(\x)\}\ll_{\xi}p^{3n-3+2\xi}\]
for any fixed $\xi>0$, which gives us the required bound (\ref{5/2}).

Turning to the proof of Lemma \ref{sscon} we will begin by supposing
that 
\begin{equation}
  \label{eq:fir}
\alpha\ge \max\{2\mu-1,2\}.
\end{equation}
It follows that $\alpha>\max\{\mu,1\}$, whence (\ref{fsum2}) yields
\begin{equation}\label{f1}
f_0(p^\alpha,p^\mu)=\frac{p^{-(3n-1)\alpha}}{1-R(p)p^{-n}}\left(  
M(p^\alpha,p^\mu)-p^{3n-1}M(p^{\alpha-1},p^\mu)\right). 
\end{equation}
It will be appropriate to observe at this point that $1-R(p)/p^n\gg 1$, which
follows from (\ref{RB}).

We proceed to compare $M(p^e,p^\mu)$ with $M(p^{e+1},p^\mu)$,
using Hensel lifting.  For $t<e$ we define
\[S_t(p^e,p^\mu):=\{\x\bmod{p^e}:\,\x\equiv\x_0\bmod{p^\mu},\, 
p^e\mid F(\x),\, p\nmid\N(\x),\,p^t\|\nabla F(\x)\}.\]
When $t<e/2$ and $t\le e-\max\{\mu,1\}$ one sees that if 
$\x\in S_t(p^e,p^\mu)$ then $\x+p^{e-t}\y\in S_t(p^e,p^\mu)$ for 
all $\y\bmod{p^t}$. Thus $S_t(p^e,p^\mu)$ is composed of cosets modulo
$p^{e-t}$. Moreover $\x+p^{e-t}\y$ will be in $S_t(p^{e+1},p^\mu)$ for
exactly $p^{3n-1}$ choices of $\y\bmod{p}$.  It follows that each
coset modulo $p^{e-t}$ in $S_t(p^e,p^\mu)$ lifts to exactly $p^{3n-1}$ 
cosets modulo $p^{e+1-t}$ in $S_t(p^{e+1},p^\mu)$, and hence that
\begin{equation}\label{compare}
\#S_t(p^{e+1},p^\mu)=p^{3n-1}\#S_t(p^e,p^\mu)
\end{equation}
for $t<e/2$ and $t\le e-\max\{\mu,1\}$.

We now write
\[T(p^t)=\#\{\x\bmod{p^t}:\, p^t|\nabla F(\x)\},\]
whence
\[\#\{\x\bmod{p^e}:\, p^t|\nabla F(\x)\}=p^{3n(e-t)}T(p^t)\]
for $t\le e$.  Then for any non-negative integer $\tau\le e$ we have
$$
M(p^e,p^\mu)=\sum_{0\le t<\tau}\#S_t(p^e,p^\mu)+
\#\left\{\x\bmod{p^e}:
\begin{array}{l}
\x\equiv\x_0\bmod{p^\mu},\\ 
p^e\mid F(\x),\\
p\nmid\N(\x),\,p^\tau |\nabla F(\x)
\end{array}
\right\},
$$
whence
\[\left|M(p^e,p^\mu)-
\hspace{-0.1cm}
\sum_{0\le t<\tau}
\hspace{-0.1cm}
\#S_t(p^e,p^\mu)\right|\le
\#\{\x\bmod{p^e}:\,p^\tau |\nabla F(\x)\}=p^{3n(e-\tau)}T(p^\tau).\]
Similarly we have
\[\left|M(p^{e+1},p^\mu)-\sum_{0\le t<\tau}\#S_t(p^{e+1},p^\mu)\right|\le
p^{3n(e+1-\tau)}T(p^\tau).\]
We therefore wish to use \eqref{compare} for every value $t<\tau$.  This will
require that $\tau-1<e/2$ and $\tau-1\le e-\max\{\mu,1\}$.  This condition
is  equivalent to requiring that $\tau-1\le (e-1)/2$ and $\tau-1\le
e-\max\{\mu,1\}$, or that $\tau\le (e+1)/2$ and $\tau\le e-\max\{\mu-1,0\}$.
Thus if 
\[\tau=\tau(e)=
\min\left\{e-\max\{\mu-1,0\}\,,\,\left[\frac{e+1}{2}\right]\right\}\]
then (\ref{compare}) yields
\[\left|M(p^{e+1},p^\mu)-p^{3n-1}M(p^e,p^\mu)\right|\le 
2p^{3n(e+1-\tau)}T(p^\tau).\]

We proceed to insert this bound into (\ref{f1}). If
$e=\alpha-1$ we find that $\tau(e)=[\alpha/2]$ providing that
$\alpha\ge\max\{2\mu-1,1\}$. In fact we have made the stronger assumption
\eqref{eq:fir}, and we deduce that
\begin{equation}\label{f0est}
f_0(p^\alpha,p^\mu)\ll p^{\alpha-3n[\alpha/2]}T(p^{[\alpha/2]}).
\end{equation}
We have therefore reduced our problem to one of providing a suitable
upper bound for $T(p^t)$. We now recall that
\[F(\x)=2^{-\kappa}\left(\tr_{L/\QQ}(\delta\bN_{K/L}(\b{p})\bN_{K/L}(\b{q}))
-2\bN_{K/\QQ}(\b{s})\right).\]
If we set
\[G(\b{p},\b{q})=\tr_{L/\QQ}(\delta\bN_{K/L}(\b{p})\bN_{K/L}(\b{q}))\]
it then follows that
\begin{equation}\label{Test}
T(p^t)\le T_1(p^t)T_2(p^t),
\end{equation}
where
\[T_1(p^t):=\#\{(\b{p},\b{q})\bmod{p^t}:\,\nabla_{\b{p}}G(\b{p},\b{q})\equiv
\nabla_{\b{q}}G(\b{p},\b{q})\equiv\b{0}\bmod{p^t}\}\]
and
\[T_2(p^t):=\#\{\b{s}\bmod{p^t}:\,\nabla\bN_{K/\QQ}(\b{s})
\equiv\b{0}\bmod{p^t}\}.\]
We begin by explaining our estimation of $T_2(p^t)$.  The treatment of
$T_1(p^t)$ will then be in the same spirit, but a little more complicated.
Over $\overline{\QQ}$ we have
\[\bN_{K/\QQ}(\b{s})=\prod_{i=1}^n L_i(\b{s})\]
for certain linearly independent linear forms 
\[L_i(\b{s})=\sum_{j=1}^n c_{ij}s_j\]
with coefficients in the normal closure $N$, say, of $K$.  Indeed our
original choice of the basis $\omega_1,\ldots,\omega_n$ ensures that
the $c_{ij}$ are algebraic integers.  We now have
\[\frac{\partial}{\partial s_j}\bN_{K/\QQ}(\b{s})=\sum_{i=1}^n c_{ij}
\mu_i,\]
where
\[\mu_i:=L_1(\b{s})\cdots L_{i-1}(\b{s})L_{i+1}(\b{s})\cdots L_n(\b{s}).\]
Let $\mathbf{C}$ denote the matrix $(c_{ij})_{i,j\le n}$ and write
$D=\det \mathbf{C}$, so
that $D^2$ is the discriminant of $K$. Suppose now that $p^t$ divides
$\nabla\N_{K/\QQ}(\b{s})$.  Then if $\bmu$ is the column vector with
elements $\mu_i$ we will have $\mathbf{C}\bmu\equiv\b{0}\bmod{p^t}$.  This
divisibility relation may be interpreted in the ring of integers for $N$.
We now pre-multiply by the matrix $\mathbf{C}^{{\rm adj}}$, whose entries are
algebraic integers, and use the fact that $\mathbf{C}^{{\rm
    adj}}\mathbf{C}=(\det \mathbf{C})I$ to
deduce that $D\bmu\equiv\b{0}\bmod{p^t}$.  We conclude that $p^t\mid
D\mu_i$ for each index $i$, where divisibility is again within the
ring of integers of $N$.  Then $p^{tn}\mid D^n\prod_i
\mu_i=D^n\bN_{K/\QQ}(\b{s})^{n-1}$, and hence $p^{\lceil tn/(n-1)\rceil}\mid
D^2\bN_{K/\QQ}(\b{s})$ whenever $p^t\mid\nabla \bN_{K/\QQ}(\b{s})$.  

Let
\[
\sigma=\max\left\{\lceil tn/(n-1)\rceil-v_p(D^2)\,,\,0\right\}.
\]
Since $t\le\sigma\le 2t$ we then have 
\begin{align*}
T_2(p^t)&=
p^{n(t-\sigma)} 
\#\{\b{s}\bmod{p^{\sigma}}:\,p^t\mid\nabla\bN_{K/\QQ}(\b{s})\}\\ 
&\le 
p^{n(t-\sigma)} 
\#\{   
\b{s}\in \NN^n\cap(0,p^\sigma ]^n :\,p^\sigma\mid\bN_{K/\QQ}(\b{s})\}. 
\end{align*}
Suppose that $(p)$ splits over $K$ as
$$
(p)=\mathfrak{p}_1^{e_1}\cdots \mathfrak{p}_r^{e_r}\quad \mbox{with}\quad
N_{K/\QQ}(\mathfrak{p}_i)=p^{f_i},
$$
for $1\le i\le r$.
Let $\alpha=\sum_{j=1}^n s_j\omega_j$, so that $\alpha$ is a non-zero element
of $\Kring$.  If we now set $v_i=v_{\mathfrak{p}_i}(\alpha)$ we deduce that
\[
p^\sigma\mid \prod_{i=1}^r  
N_{K/\QQ}(\mathfrak{p}_i)^{v_i},
\]
whence $\sigma\le\sum f_iv_i$.  It then follows that $\sigma\le\sum
f_ig_i$, where $g_i=\min\{v_i,\sigma\}$.  Thus for each element $\b{s}$
there are non-negative integers $g_i\le\sigma$ for which $\sum
f_ig_i\ge\sigma$ and
\[
\sum_{j=1}^n s_j\omega_j\in \mathfrak{p}_1^{g_1}\cdots \mathfrak{p}_r^{g_r}.
\]
This condition restricts $\b{s}$ to a lattice $\Lambda$ say, with
$p^\sigma\ZZ^n\subseteq\Lambda\subseteq\ZZ^n$, and with
$p^\sigma\mid\det(\Lambda)$.  The number of choices for $g_1,\ldots,g_r$
is at most $(\sigma+1)^r\le(\sigma+1)^n$.
It therefore follows that
\begin{equation}\label{T2est}
T_2(p^t)\le   
p^{n(t-\sigma)}
(\sigma+1)^n p^{n\sigma-\sigma}
\le (2t+1)^n p^{nt-\sigma}\ll 
(2t+1)^np^{nt-\lceil tn/(n-1)\rceil},
\end{equation}
since $D$ is fixed.

We now examine $T_1(p^t)$ in a similar way. The form $G(\b{p},\b{q})$ takes
the shape
\[\delta\prod_{i=1}^{n/2}L_i(\b{p})\prod_{i=1}^{n/2}L_i(\b{q})+
\delta^{\sigma}\prod_{i=n/2+1}^{n}L_i(\b{p})\prod_{i=n/2+1}^{n/2}L_i(\b{q}).\]
If we replace $\delta L_1(\b{p})$ by $L_1(\b{p})$ and 
$\delta^\sigma L_{n/2+i}(\b{p})$ by $L_{n/2+1}(\b{p})$ we find that 
\[\frac{\partial G(\b{p},\b{q})}{\partial p_j}=\sum_{i=1}^n c_{ij}\mu_i,\]
where
\[c_{ij}=\frac{\partial L_i(\x)}{\partial x_j}\]
and
\[\mu_i=L_1(\b{p})\cdots L_{i-1}(\b{p})L_{i+1}(\b{p})\cdots L_{n/2}(\b{p})
\prod_{i=1}^{n/2}L_i(\b{q})\]
for $1\le i\le n/2$ and
\[\mu_i=L_{n/2+1}(\b{p})\cdots L_{n/2+i-1}(\b{p})L_{n/2+i+1}(\b{p})
\cdots L_{n}(\b{p})
\prod_{i=n/2+1}^n L_i(\b{q})\]
for $n/2<i\le n$.
Thus $p^t\mid D\mu_i$ for $1\le i\le n$, where $D$ is the determinant
of the matrix $\mathbf{C}=(c_{ij})_{i,j\le n}$.  In particular $D$ is a non-zero
algebraic integer.  Taking the product for $i\le n$ yields $p^{tn}\mid
D^n \bN_{K/\QQ}(\b{p})^{n/2-1}\bN_{K/\QQ}(\b{q})^{n/2}$.  
By symmetry we also obtain
$p^{tn}\mid D^n \bN_{K/\QQ}(\b{q})^{n/2-1}\bN_{K/\QQ}(\b{p})^{n/2}$, and hence
\[p^{2tn}\mid D^{2n} \bN_{K/\QQ}(\b{p})^{n-1}\bN_{K/\QQ}(\b{q})^{n-1}.\]
It follows that $p^\sigma \mid\bN_{K/\QQ}(\b{p})\bN_{K/\QQ}(\b{q})$ with
\[\sigma=\max\left\{\lceil 2tn/(n-1)\rceil-v_p(D^4)\,,\,0\right\}.\]

We can now complete the argument as before.  This time there will be
two sets of non-negative exponents $g_i$ and $g_i'$, say, corresponding
to $\b{p}$ and $\b{q}$ respectively, and such that $\sum_i
f_i(g_i+g_i')\geq \sigma$.  Thus the factor $(\sigma+1)^n$ must be replaced by
$(\sigma+1)^{2n}$. For each such set of exponents we 
find that $(\b{p},\b{q})$ is
restricted to a sublattice of $\ZZ^{2n}$ of index at least $p^\sigma$, 
and we deduce as before that
\[T_1(p^t)\ll (4t+1)^{2n}p^{2nt-\lceil 2tn/(n-1)\rceil}.\]
In view of (\ref{Test}) and (\ref{T2est}) we deduce that
\[T(p^t)\ll 
(4t+1)^{3n}p^{3nt-\lceil tn/(n-1)\rceil-\lceil 2tn/(n-1)\rceil}
\ll (4t+1)^{3n}p^{3nt-3t-2},\]
since if $m\in\ZZ$ we have 
$\lceil\theta\rceil\ge m+1$ for any real number $\theta>m$.  It now
follows from (\ref{f0est}) that
\[f_0(p^\alpha,p^\mu)\ll (2\alpha+1)^{3n}p^{\alpha-3[\alpha/2]-2},\]
and this suffices for Lemma \ref{sscon} when 
\eqref{eq:fir} holds. 
Finally, it is clear in fact that we still have this estimate
if $2\le \alpha\le 2\mu-1$, since if $\mu>0$ then $p$ belongs to the
finite set of divisors of $M$. This therefore 
completes the proof of Lemma~\ref{sscon}.

\bigskip

Our final task in this section is to establish a 
lower bound for $\sigma_\infty$ to
complement the upper bound in 
\eqref{eq:upper-sig}.
For any $\b{c}\in \RR^2$ let $B(\b{c};\rho)\subset \RR^2$ denote the box
centered on $\b{c}$ with side length $2\rho$. 
The final part of Lemma \ref{omegaprops} implies that 
there exists $\rho\gg G^{-1}U^{n/2}$ such that 
$$
\omega(x_1,x_2)\gg G^{2-n}
$$ 
for every $(x_1,x_2)\in
B(\delta \bN_{K/L}(U\b{u}^{(\RR)});\rho)$.  Here we view 
$\delta \bN_{K/L}(U\b{u}^{(\RR)})$ as a vector $(c_1,c_2)$ in  
$\RR^2$.
In this way we deduce from \eqref{eq:defJ} and \eqref{eq:sig-inf}
that 
\begin{equation}\label{last}
\sigma_\infty\gg G^{2-n} \sum_{\substack{\w\in \cW\cap \ZZ^n}}
\sum_{\substack{\v\in \cV\cap\ZZ^{n}}}
J(\v;\w),
\end{equation}
where if $a_i=a_i(\v)$ and $b=b(\w)$ are given by \eqref{eq:aib} then
$$
J(\v;\w):=
\meas\left\{x\in \RR:   
\left(
a_2(\v)x, -a_1(\v)x+\frac{b(\w)}{
a_2(\v)}\right)\in 
B(\b{c};\rho)
\right\}.
$$
The minimum distance from the line $(a_2x,-a_1x+b/a_2)$ to the point
$(c_1,c_2)$ is equal to $|a_1c_1+a_2c_2-b|/\sqrt{a_1^2+a_2^2}$.
Suppose now that we have points $\v$ and $\w$ satisfying
\begin{equation}\label{smaller}
|\v-V\v^{(\RR)}|<\lambda G^{-1}V\qquad\mbox{and}\qquad
|\w-W\w^{(\RR)}|<\lambda G^{-1}W
\end{equation}
for some $\lambda\le 1$.  If we set
\[a_1^{(\RR)}=\tr_{L/\QQ}(\bN_{K/L}(\v^{(\RR)})),\qquad
a_2^{(\RR)}=\tr_{L/\QQ}(\tau \bN_{K/L}(\v^{(\RR)}))\]
and
\[b^{(\RR)}=2\bN_{K/\QQ}(\w^{(\RR)}),\]
then
\[a_1(\v)=V^{n/2}a_1^{(\RR)}+O(\lambda G^{-1}V^{n/2}),\qquad
a_2(\v)=V^{n/2}a_2^{(\RR)}+O(\lambda G^{-1}V^{n/2})\]
and
\[b(\w)=W^{n}b^{(\RR)}+O(\lambda G^{-1}W^{n}).\]
Moreover
\[c_1V^{n/2}a_1^{(\RR)}+c_2 V^{n/2}a_2^{(\RR)}-W^{n}b^{(\RR)}
=F(V\v^{(\RR)};W\w^{(\RR)};U\u^{(\RR)})=0.\]
Thus $|a_1c_1+a_2c_2-b|=O(\lambda G^{-1}W^n)$, while
$a_1^2+a_2^2\gg V^n$ since
$a_1^{(\RR)},a_2^{(\RR)}$ are not both zero. 
Since $\rho\gg G^{-1}U^{n/2}$, it
follows that
\[\frac{|a_1c_1+a_2c_2-b|}{\sqrt{a_1^2+a_2^2}}\le \frac{\rho}{2}\]
providing that we take $\lambda$ as a sufficiently small positive constant.

We now see that, for points $(\v,\w)$ satisfying (\ref{smaller}), the line
$(a_2x,-a_1x+b/a_2)$ meets the disc $B(\delta\bN_{K/L}(U\b{u}^{(\RR)});\rho)$
in a segment of length $\gg\rho$, so that
\[J(\v;\w)\gg \frac{\rho}{\max\{|a_1| ,|a_2|\}}\gg
\frac{\rho}{V^{n/2}}\gg G^{-1}H^{n/2}.\]
The number of available points $(\v,\w)$ is $\gg G^{-2n}(VW)^n$ and we
therefore conclude from (\ref{last}) that
\[\sigma_\infty\gg G^{1-3n}H^nV^{2n}.\]
This  should be compared with the upper bound 
\eqref{eq:upper-sig}.

\bigskip

Bringing together our lower bounds for $\mathfrak{S}$ and $\sigma_\infty$ 
in Lemma \ref{lem:onion}, we 
deduce that 
$$
\cl{M}\gg G^{1-3n} H^nV^{2n}
+ O_\eta\left(Q^6 H^{n/2} V^{2n+2\eta}\right),
$$
provided that $H^{n/2}\leq V$. 
This therefore leads to the following conclusion.

\begin{lemma}\label{lem:final-main'}
Let $G=\log V$. 
Assume that $H^{n/2}\leq V$ and 
$$
Q\leq H^{n/12}V^{-\eta/2}.
$$ 
Then we have 
$$
\cl{M}(G,H,V)\gg (\log V)^{1-3n} H^nV^{2n}.
$$
\end{lemma}

Recalling \eqref{eq:me} and our choice $G=\log V$, it 
is now time to select parameters $Q, H,V$ such that
$
\mathcal{E}(G,H,V)=o(\mathcal{M}(G,H,V)).
$
We will choose
$
Q= H^{(n-1)/12},
$
with which choice  Lemma~\ref{lem:final-main'}
implies that 
$\cl{M}(G,H,V)\gg (\log V)^{1-3n} H^nV^{2n}$, if 
$$
V^{6\eta}\leq H\leq V^{2/n}.
$$
In line with Lemma \ref{lem:suffice}  we let $V$ run through large integers congruent to $1$ modulo $M$. Next we choose 
$H_0=1+M[V^{1/(10 n^2)}]$, which is a positive  integer 
congruent to $1$ modulo $M$.
But then $H=H_0^2$ has order  $V^{1/(5n^2)}$ and so  
$Q=H^{(n-1)/12}$
satisfies the conditions 
of  Lemma \ref{lem:treat_E}. This implies that the required estimate
for $\mathcal{E}(G,H,V)$ holds and so  completes the proof of Theorem \ref{thmain-a}.


\begin{thebibliography}{9}
\bibitem{bartels-81}
H.-J. Bartels, Zur Arithmetik von Konjugationsklassen in 
algebraischen Gruppen. {\em 
J.\ Algebra}  {\bf 70}  (1981), 179--199.

\bibitem{bartels-81'}
H.-J. Bartels, Zur Arithmetik von
Diedergruppenerweiterungen. {\em Math.\ Ann.}  {\bf 256}  (1981),
465--473.   

\bibitem{brudern}
J. Br\"udern, 
Binary additive problems and the circle method, multiplicative
sequences and convergent sieves. {\em Analytic number theory}, 91--132,
Cambridge Univ. Press, Cambridge, 2009.  

\bibitem{ct-salb}
J.-L. Colliot-Th\'el\`ene and P. Salberger,
Arithmetic on some singular cubic hypersurfaces.
{\em Proc.\ London Math.\ Soc.}  {\bf 58} (1989),  519--549.

\bibitem{CTSb}  J.-L. Colliot-Th\'el\`ene and J.J. Sansuc,
La $R$-\'equivalence sur les tores. {\em Ann.\ Sci.\ \'Ecole
Norm.\ Sup.} {\bf 10} (1977), 175--229.

\bibitem{CTSa}  J.-L. Colliot-Th\'el\`ene and J.J. Sansuc,
Principal homogeneous spaces under flasque tori: applications.
{\em J.\ Algebra} {\bf 106} (1987), 148--205.

\bibitem{ct-skoro}  J.-L. Colliot-Th\'el\`ene and A.N. Skorobogatov,
Descent on fibrations over $\PP_k^1$ revisited.
{\em Math.\ Proc.\ Camb.\ Phil.\ Soc.} {\bf 128} (2000), 383--393.


\bibitem{CTHS}J.-L. Colliot-Th\'el\`ene, D. Harari and
  A.N. Skorobogatov, 
Valeurs d'un polyn\^ome \`a une variable repr\'esent\'es par une 
norme. {\em Number theory and algebraic geometry}, 69--89,
London Math.\ Soc.\ Lecture Note Ser.  {\bf 303}
Cambridge Univ.\ Press, Cambridge, 2003.

\bibitem{crelle}
J.-L. Colliot-Th\'el\`ene, J.-J. Sansuc and P. Swinnerton-Dyer, 
Intersections of two quadrics and Ch\^atelet surfaces, I.
{\em J.\ reine angew.\ Math.}  {\bf 373}  (1987), 37--107;
II. {\em J.\ reine angew.\ Math.}  {\bf 374}  (1987), 72--168.

\bibitem{CT-94}
J.-L. Colliot-Th\'el\`ene and P. Swinnerton-Dyer, 
Hasse principle and weak approximation for pencils of Severi--Brauer
and similar varieties. {\em J.\ reine angew.\ Math.}  {\bf 453}  (1994),
49--112.

\bibitem{98a}J.-L. Colliot-Th\'el\`ene, A.N. Skorobogatov and 
P. Swinnerton-Dyer, 
Rational points and zero-cycles on fibred varieties: Schinzel's hypothesis
and Salberger's device.
{\em J.\ reine angew.\ Math.} {\bf 495} (1998), 1--28.

\bibitem{FI}
E. Fouvry and H. Iwaniec, Gaussian primes.
{\em Acta Arith.} {\bf 79} (1997), 249--287.


\bibitem{gras}G. Gras, {\em Class field theory: from theory to
 practice}. Springer-Verlag, 2002.

\bibitem{GT}B. Green and T. Tao, Linear equations in primes.
{\em Annals of Math.} {\bf 171} (2010),  1753--1850.

\bibitem{GTZ}B. Green, T. Tao and T. Ziegler, 
An inverse theorem for the Gowers $U^{s+1}[N]$-norm.
{\em Submitted}, 2010.

\bibitem{gurak}
S. Gurak, On the Hasse norm principle.
{\em J.\ reine angew.\ Math.} {\bf 299/300} (1978), 16--27.


\bibitem{RMI}D.R. Heath-Brown, 
The ternary Goldbach problem. {\em
Rev.\ Mat.\ Iberoamericana} {\bf 1} (1985), 45--59.

\bibitem{HBS} D.R. Heath-Brown and A.N. Skorobogatov, Rational
solutions of certain equations involving norms. 
{\em Acta Math.} {\bf 189} (2002), 161--177.

\bibitem{hooley}C. Hooley,
On the number of points on a complete intersection over a finite field.
{\em J. Number Theory} {\bf 38} (1991), 338--358. 

\bibitem{hux}
M.N. Huxley, The large sieve inequality for algebraic 
number fields. {\em Mathematika}  {\bf 15} (1968), 178--187.

\bibitem{IK}
H. Iwaniec and E. Kowalski, {\em Analytic number theory}. 
American Math.\ Soc.\ Colloq.\ Pub.\ {\bf 53}, American Math.\ Soc., 2004. 


\bibitem{linnik}
J.V. Linnik, {\em 
The dispersion method in binary additive problems}. 
American Math.\ Soc.\, Providence, R.I., 1963. 


\bibitem{sansuc} 
J.-J. Sansuc, 
Groupe de Brauer et arithm\'etique des groupes alg\'ebriques lin\'eaires sur un corps de nombres, {\em J.\ reine angew.\ Math.}  {\bf 327}  (1981), 12--80.

\bibitem{vosk} 
V.E. Voskresenski$\breve{\i}$, 
{\em Algebraic groups and their birational invariants} (translated from 
Russian  by B. Kunyavski$\breve{\i}$). 
Translations of Math. Monographs {\bf 179}, American Math.\ Soc.\, 1998. 

\end{thebibliography}
\end{document}